\documentclass[english]{article}
\usepackage{custom_tex}
\usepackage[algo2e,linesnumbered,ruled]{algorithm2e}
\begin{document}
\definecolor{ed}{RGB}{225,0,0}
\definecolor{yue}{RGB}{0,100,100}
\newcommand{\ed}[1]{\textcolor{ed}{[ED: #1]}}
\newcommand{\yue}[1]{\textcolor{yue}{[Yue: #1]}}

\title{WONDER: Weighted one-shot distributed ridge regression in high dimensions\footnote{A previous version of the manuscript had the  title "One-shot distributed ridge regression in high dimensions".}}
\date{\today}
\author{Edgar Dobriban\footnote{Wharton Statistics Department, University of Pennsylvania. 3730 Walnut Street, Philadelphia, PA, USA. E-mail: \texttt{dobriban@wharton.upenn.edu}.} \, and Yue Sheng\footnote{Graduate Group in Applied Mathematics and Computational Science, Department of Mathematics, University of Pennsylvania. E-mail: \texttt{yuesheng@sas.upenn.edu}.}}


\maketitle
\begin{abstract}
In many areas, practitioners need to analyze large datasets that challenge conventional single-machine computing. To scale up data analysis, distributed and parallel computing approaches are increasingly needed. 

Here we study a fundamental and highly important problem in this area: How to do ridge regression in a distributed computing environment? Ridge regression is an extremely popular method for supervised learning, and has several optimality properties, thus it is important to study. We study one-shot methods that construct weighted combinations of ridge regression estimators computed on each machine. By analyzing the mean squared error in a high dimensional random-effects model where each predictor has a small effect, we discover several new phenomena. 

{\bf Infinite-worker limit}: The distributed estimator works well for very large numbers of machines, a phenomenon we call ``infinite-worker limit". 

{\bf Optimal weights:} The optimal weights for combining local estimators sum to more than unity, due to the downward bias of ridge. Thus, all averaging methods are suboptimal. 

We also propose a new Weighted ONe-shot DistributEd Ridge regression (WONDER) algorithm. We test WONDER in simulation studies and using the Million Song Dataset as an example. There it can save at least 100x in computation time, while nearly preserving test accuracy.
\end{abstract}

\setcounter{tocdepth}{1}



\section{Introduction}

Computers have changed all aspects of our world. Importantly, computing has made data analysis more convenient than ever before. However, computers also pose limitations and challenges for data science.  For instance, hardware architecture is based on a model of a universal computer---a Turing machine---but in fact has physical limitations of storage, memory, processing speed, and communication bandwidth over a network. As large datasets become more and more common in all areas of human activity, we need to think carefully about working with these limitations. 

How can we design methods for data analysis (statistics and machine learning) that scale to large datasets? A general approach is \emph{distributed and parallel computing}. Roughly speaking, the data is divided up among computing units, which perform most of the computation locally, and synchronize by passing relatively short messages. While the idea is simple, a good implementation can be hard and nontrivial. Moreover, different problems have different inherent needs in terms of local computation and global communication resources. For instance, in statistical problems with high levels of noise, simple one-shot schemes like averaging estimators computed on local datasets can sometimes work well. 

In this paper, we study a fundamental problem in this area. We are interested in linear regression, which is arguably one of the most important problems in statistics and machine learning. A popular method for this model is \emph{ridge regression} (aka Tikhonov regularization), which regularizes the estimates using a quadratic penalty to improve estimation and prediction accuracy. We aim to understand how to do ridge regression in a distributed computing environment. We are also interested in the important \emph{high-dimensional} setting, where the number of features can be very large. In fact our approach allows the dimension and sample size to have any ratio. We also work in a random-effects model where each predictor has a small effect on the outcome, which is the model for which ridge regression is best suited.

We consider the simplest and most fundamental method, which performs ridge regression locally on each dataset housed on the individual machines or other computing units, sends the estimates to a global datacenter (or parameter server), and then constructs a final one-shot estimator by taking a linear combination of the local estimates. As mentioned, such methods are sometimes near-optimal, and it is therefore well-justified to study them. We will later give several additional justifications for our work. 

However, in contrast to existing work, we introduce a completely new mathematical approach to the problem, which has never been used for studying distributed ridge regression before. Specifically, we leverage and further develop sophisticated recent techniques from random matrix theory and free probability theory in our analysis. This enables us to make important contributions, that were simply unattainable using more ``traditional" mathematical approaches.

To give a sense of our results, we provide a brief discussion here. We have a dataset consisting of $n$ datapoints, for instance 1000 heart disease patients. Each datapoint has an outcome $y_j$, such as blood pressure, and features $x_j$, such as age, height, electronic health records, lab results, and genetic variables. Our goal is to predict the outcome of interest (i.e., blood pressure) for new patients based on their features, and to estimate the relationship of the outcome to the features.

The samples are distributed across several sites, for instance patients from different countries are housed in different data centers. We will refer to the sites as ``machines", though they may actually be other computing entities, such as entire computer networks or data centers. In many important settings, it can be impossible to share the data across the different sites, for instance due to logistical or privacy reasons. 

Therefore, we assume that each site has a subset of the samples. Our approach is to train ridge regression on this local data. As usual, we can arrange the local dataset (say on the $i$-th machine) into a feature matrix $X_i$, where each row contains a sample (i.e., datapoint), and an outcome vector $Y_i$ where each entry is an outcome. We compute the local ridge regression estimates $$\hat \beta_i = (X_i^\top X_i+\lambda_iI_p)^{-1}X_i^\top Y_i,$$ where $\lambda_i$ are some regularization parameters. We then aggregate them by a weighted combination, constructing the final \emph{one-shot distributed ridge estimator} (where $k$ is the number of sites)
$$\hbeta_{dist}=\sum_{i=1}^kw_i\hbeta_i.$$

The important questions here are: 
\benum 
\item How does this work? 
\item How to tune the parameters? (such as the regularization parameters and weights)
\eenum

Question (1) is of interest because we wish to know when one-shot methods are a good approach, and when they are not. For this we need to understand the performance as a function of the key problem parameters, such as the signal strength, sample size, and dimension. For question (2), the challenge is posed by the constraints of the distributed computing environment, where standard methods for parameter tuning such as cross-validation may be expensive.

In this work we are able to make several crucial contributions to these questions. We work in an asymptotic setting where $n, p$ grow to infinity at the same rate, which effectively gives good results for any $n,p$. We study a linear-random effects model, where each regressor has a small random effect on the outcome. This is a good model for the applications where ridge regression is used, because ridge does not assume sparsity, and has optimality properties in certain dense random effects models. Importantly, this analysis does \emph{not} assume any sparsity in a high-dimensional setting. Sparsity has been one of the biggest driving forces in statistics and machine learning in the last 20 years. Our work is in a different line of work, and shows that meaningful results are available without sparsity.

\begin{figure}
\centering
\begin{subfigure}{.5\textwidth}
  \centering
  \includegraphics[scale=0.5]{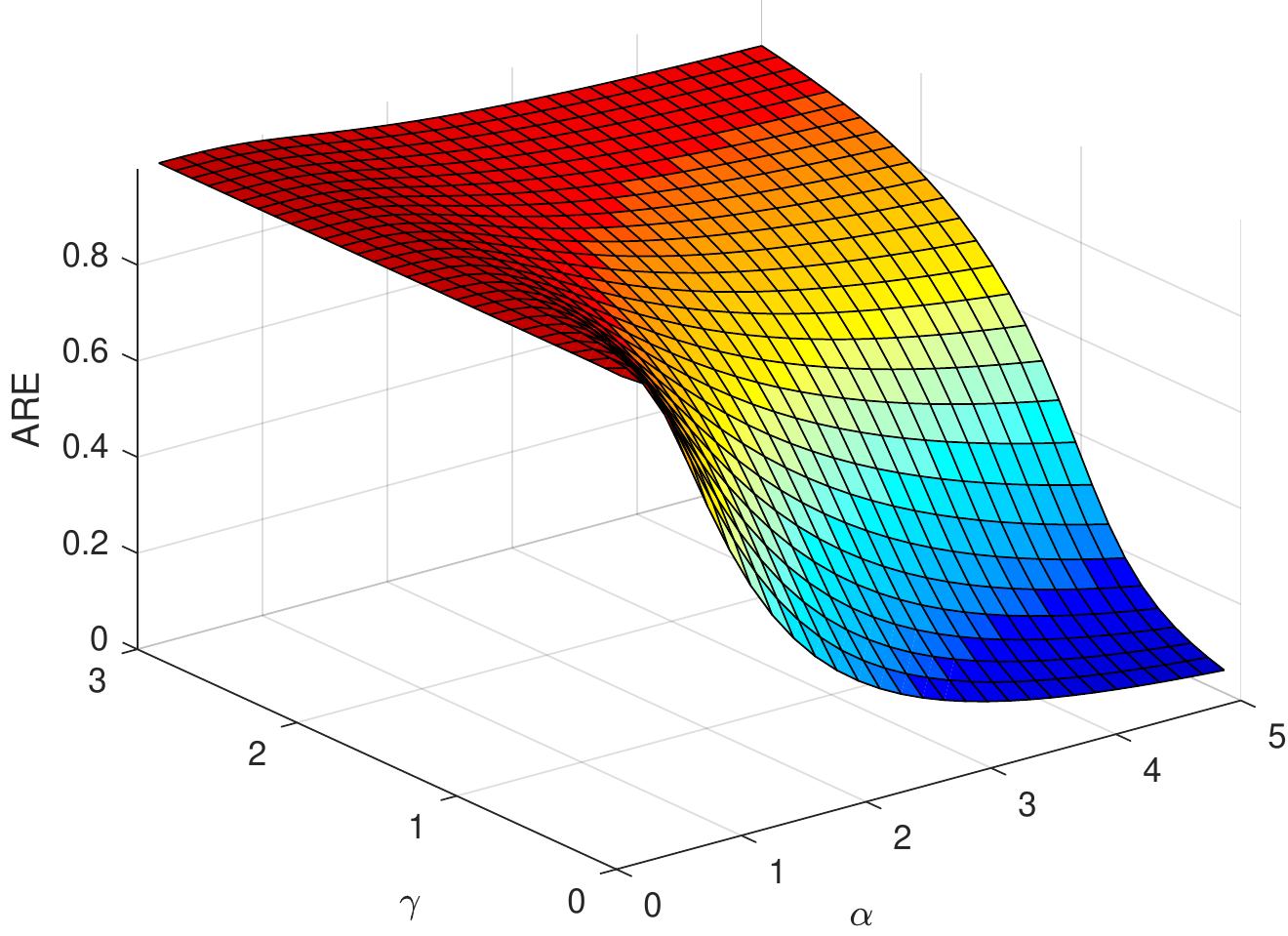}
  \caption{Surface}
\end{subfigure}%
\begin{subfigure}{.5\textwidth}
  \centering
  \includegraphics[scale=0.5]{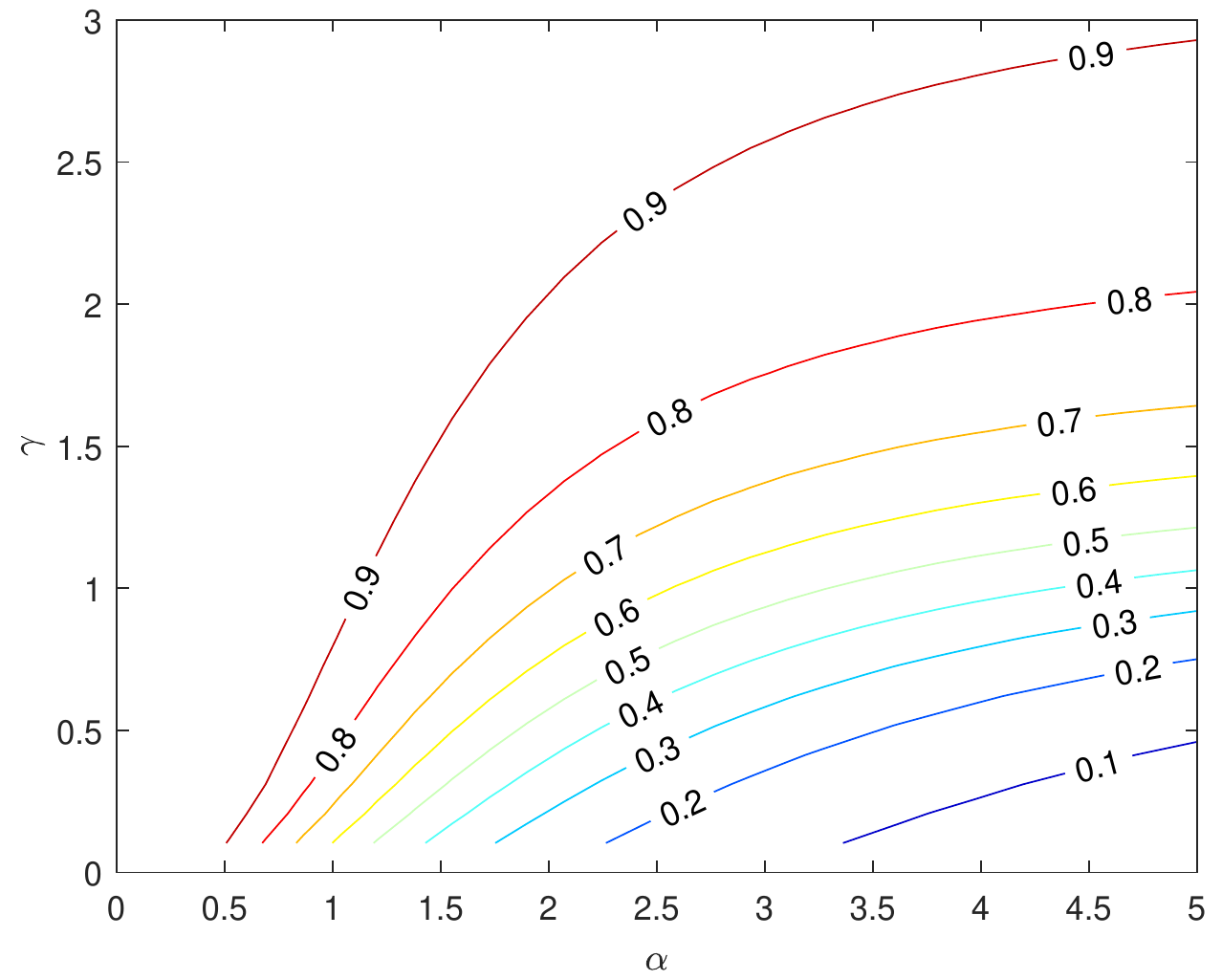}
  \caption{Contour}
\end{subfigure}
\caption{Efficiency loss due to one-shot distributed learning. This plot shows the relative mean squared error of centralized ridge regression compared to optimally weighted one-shot distributed ridge regression. This quantity is at most unity, and the larger, the ``better" distributed ridge works. Specifically, the model is asymptotic, and we show the dependence of the Asymptotic Relative Efficiency (ARE) on the aspect ratio $\gamma = \lim p/n$ (where $n$ is sample size and $p$ is dimension) and on the signal strengh $\alpha = \sqrt{\E\|\beta\|^2}$, in the \emph{infinite-worker limit} when we distribute our data over many machines. We show (a) surface and (b) contour plots of the ARE. See the text for details.}
  \label{arefig}
\end{figure}

We find the limiting mean squared error of the one-shot distributed ridge estimator. This enables us to characterize the optimal weights and tuning parameters, as well as the \emph{relative efficiency} compared to centralized ridge regression, meaning the ratio of the risk of usual ridge to the distributed estimator. This can precisely pinpoint the computation-accuracy tradeoff achieved via one-shot distributed estimation. See Figure \ref{arefig} for an illustration.

As a consequence of our detailed and precise risk analysis, we make several qualitative discoveries that we find quite striking: 

\benum 
\item {\bf Efficiency depends strongly on signal strength.} The statistical efficiency of the one-shot distributed ridge estimator depends strongly on signal strength. The efficiency is generally high (meaning distributed ridge regression works well) when the signal strength is low.

\item {\bf Infinite-worker limit.} The one-shot distributed estimator does not lose all efficiency compared to the ridge estimator even in the limit of \emph{infinitely many  machines}. Somewhat surprisingly, this suggests that simple one-shot weighted combination methods for distributed ridge regression can work well even for very large numbers of machines. The statement that this can be achieved by communication-efficient methods is nontrivial. This finding is clearly important from a practical perspective.

\item {\bf Decoupling.} When the features are uncorrelated, the problem of choosing the optimal regularization parameters \emph{decouples} over the different machines. We can choose them in a locally optimal way, and they are also globally optimal. We emphasize that this is a very delicate result, and is not true in general for correlated features. Moreover, this discovery is also important in practice, because it gives conditions under which we can choose the regularization parameters separately for each machine, thus saving valuable computational resources.

\item {\bf Optimal weights do not sum to unity.} Our work uncovers unexpected properties of the optimal weights. Naively, one may think that the weights need to sum to unity, meaning that we need a weighted average. However, it turns out the optimal weights sum to more than unity, because of the negative bias of the ridge estimator. This means that \emph{any type of averaging method is suboptimal}. We characterize the optimal weights and under certain conditions find their explicit analytic form.
\eenum

Based on these results, we propose a new Weighted ONe-shot DistributEd Ridge regression algorithm (WONDER). We also confirm these results in detailed simulation studies and on an empirical data example, using the Million Song Dataset. Here WONDER can be used over 100-way splits of the data with 5\% loss of prediction accuracy.

We also emphasize that some aspects of our work can help practitioners directly (e.g., our new algorithm), while others are developed for deepening our understanding of the nature of the problem. We discuss the practical implications of our work in Section \ref{rc}.

The paper is structured as follows. We discuss some related work in Section \ref{relw}. We start with finite sample results in Section \ref{fin_sam}. We provide asymptotic results for features with an arbitrary covariance structure in Section \ref{gen_cov}. We consider the special case of an identity covariance in Section \ref{id_cov}. In Section \ref{algo} we provide an explicit algorithm for optimally weighted one-shot distributed ridge. We also study in detail the properties of the estimation error, relative efficiency (including minimax properties in Section \ref{mmx}), tuning parameters (and decoupling), as well as optimal weights, including answers to the questions above.  We provide numerical simulations throughout the paper, and additional ones in Section \ref{expr}, along with an example using an empirical dataset. 
The code for our paper is available at \url{github.com/dobriban/dist_ridge}.




\subsection{Related work}
\label{relw}
Here we discuss some related work. Historically, distributed and parallel computation has first been studied in computer science and optimization \citep[see e.g.,][]{bertsekas1989parallel,lynch1996distributed,blelloch2010parallel,boyd2011distributed,rauber2013parallel,koutris2018algorithmic}. However, the problems studied there are quite different from the ones that we are interested in. Those works often focus on problems where correct answers are required within numerical precision, e.g., 16 bits of accuracy. However, when we have noisy datasets, such as in statistics and machine learning, numerical precision is neither needed nor usually possible. We may only hope for 3-4 bits of accuracy, and thus the problems are different.


The area of distributed statistics and machine learning has attracted increasing attention only relatively recently, see for instance \cite{mcdonald2009efficient,zhang2012communication,zhang2013communication, li2013statistical, zhang2013divide, duchi2014optimality, chen2014split, mackey2011divide, zhang2015divide, braverman2016communication, jordan2016communication, rosenblatt2016optimality, smith2016cocoa, banerjee2016divide, zhao2016partially, xu2016optimal,fan2017distributed, lin2017distributed, lee2017communication, volgushev2017distributed, shang2017computational, battey2018distributed, zhu2018distributed, xichen2018Quantile, xichen2018fone, xichen2018svm,shi2018massive,duan2018distributed,liu2018many,cai2020distributed}, and the references therein. See \cite{huo2018aggregated} for a review. We can only discuss the most closely related papers due to space limitations. 

\cite{zhang2013communication} study the MSE of averaged estimation in empirical risk minimization. Later \cite{zhang2015divide} study divide and conquer \emph{kernel ridge regression}, showing that the partition-based estimator achieves the statistical minimax rate over all estimators, when the number of machines is not too large. These results are very general, however they are not as explicit or precise as our results. In addition they consider fixed dimensions, whereas we study increasing dimensions under random effects models. \cite{lin2017distributed} improve the above results, removing certain eigenvalue assumptions on the kernel, and sharpening the rate. 

\cite{guo2017learning} study regularization kernel networks, and propose a debiasing scheme that can improve the behavior of distributed estimators. This work is also in the same framework as those above (general kernel, fixed dimension). \cite{xu2016optimal} propose a distributed General Cross-Validation method to choose the regularization parameter. 

\cite{rosenblatt2016optimality} consider averaging in distributed learning in fixed and high-dimensional M-estimation, without studying regularization. \cite{lee2017communication} study sparse linear regression, showing that averaging debiased lasso estimators can achieve the optimal estimation rate if the number of machines is not too large. A related work is \cite{battey2018distributed}, which also includes hypothesis testing under more general sparse models. These last two works are on a different problem (sparse regression), whereas we study ridge regression in random-effects models.

\section{Finite sample results}
\label{fin_sam}

We start our study of distributed ridge regression by a finite sample analysis of estimation error in linear models. Consider the standard linear model
\beq
Y = X\beta+\ep.
\label{lrm}
\eeq
Here $Y\in\mathbb R^n$ is the $n$-dimensional continuous outcome vector of $n$ independent samples (e.g., the blood pressure level of $n$ patients, or the amount of time spent on an activity by $n$ internet users), $X$ is the $n\times p$ design matrix containing the values of $p$ features for each sample (e.g., demographical and genetic variables of each patient). Moreover, $\beta=(\beta_1,\dots,\beta_p)^\top\in\mathbb R^p$ is the $p$-dimensional vector of unknown regression coefficients. 

Our goals are to predict the outcome variable for future samples, and also to estimate the regression coefficients. The outcome vector is affected by the random noise $\ep=(\ep_1,\dots,\ep_n)^\top\in\mathbb R^n$.  We assume that the coordinates of $\ep$ are independent random variables with mean zero and variance $\sigma^2$.

The \emph{ridge regression} (or Tikhonov regularization) estimator is one of the most popular methods for estimation and prediction in linear models. Recall that the ridge estimator of $\beta$ is 
$$\hbeta(\lambda) = (X^\top X+n\lambda I_p)^{-1} X^\top Y,$$
where $\lambda$ is a tuning parameter. This estimator has many justifications. It shrinks the coefficients of the usual ordinary least squares estimator, which can lead to improved estimation and prediction. When the entries of $\beta$ and $\ep$ are iid Gaussian, and for suitable $\lambda$, it is the posterior mean of $\beta$ given the outcomes, and hence is a Bayes optimal estimator for any quadratic loss function, including estimation and prediction error. 

Suppose now that we are in a distributed computation setting. The samples are distributed across $k$ different sites or machines. For instance, the data of users from a particular country may be stored in a separate datacenter.  This may happen due to memory or storage limitations of individual data storage facilities, or may be required by data usage agreements. As mentioned, for simplicity we call the sites ``machines".

 We can write the partitioned data as 
$$X = 
\begin{bmatrix}
    X_1 \\
    \hdots \\
    X_k
\end{bmatrix}
,\,\,
Y = 
\begin{bmatrix}
    Y_1 \\
    \hdots \\
    Y_k
\end{bmatrix}.
$$
Thus the $i$-th machine contains $n_i$ samples whose features are stored in the $n_i\times p$ matrix $X_i$  and also the corresponding $n_i\times 1$ outcome vector $Y_i$. 

Since the ridge regression estimator is a widely used gold standard method, we would like to understand how we can approximate it in a distributed setting. Specifically, we will focus on one-shot \emph{weighting} methods, where we perform ridge regression locally on each subset of the data, and then aggregate the regression coefficients by a weighted sum. There are several reasons to consider weighting methods:
\benum 
\item This is a practical method with \emph{minimal communication cost}. When communication is expensive, it is imperative to develop methods that minimize communication cost. In this case, one-shot weighting methods are attractive, and so it is important to understand how they work. In a well-known course on scalable machine learning, Alex Smola calls such methods ``idiot-proof" \citep{SML}, meaning that they are straightforward to implement (unlike some of the more sophisticated methods).

\item Averaging (which is a special case of one-shot weighting) has already been studied in several works on distributed ridge regression (e.g., \cite{zhang2015divide,lin2017distributed}), and much more broadly in distributed learning, see the related work section for details. Such methods are \emph{known to be rate-optimal} under certain conditions. 
\item However, in our setting, we are able to discover several \emph{new phenomena} about one-shot weighting. For instance, we can quantify in a much more nuanced way the accuracy loss compared to centralized ridge regression.

\item Weighting may serve as a useful \emph{initialization to iterative methods}. In practical distributed learning problems, iterative optimization algorithms such as distributed gradient descent or ADMM \citep{boyd2011distributed} may be used. However, there are examples where the first step of the iterative method has \emph{worse} performance than a simple averaging \citep{pourshafeie2018caring}. Therefore, we can imagine hybrid or warm start methods that use weighting as an initialization. This also suggests that studying one-shot weighting is important.
\eenum

Therefore, we define \emph{local} ridge estimators for each dataset $X_i,Y_i$, with regularization parameter $\lambda_i$ as
$$\hat \beta_i(\lambda_i) = (X_i^\top X_i+n_i\lambda_iI_p)^{-1}X_i^\top Y_i.$$
We consider combining the local ridge estimators at a central server via a one-step weighted summation. We will find the optimally weighted one-shot distributed estimator 
$$\hbeta_{dist}(w)=\sum_{i=1}^kw_i\hbeta_i.$$
Note that, unlike ordinary least squares (OLS), the local ridge estimators are always well-defined, i.e. $n_i$ can be smaller than $p$. Also, for the distributed OLS estimator averaging local OLS solutions, it is natural to require $\sum_iw_i=1$, because this ensures unbiasedness \citep{dobriban2018Distributed}. However, the ridge estimators are biased, so it is not clear if we should put any constraints on the weights. In fact we will find that the optimal weights typically do not sum to unity. These features distinguish our work from prior art, and lead to some surprising consequences. 

Throughout the paper, we will frequently use the notations $\hSigma=n^{-1}X^\top X$ and $\hSigma_i=n_i^{-1}X_i^\top X_i$. A stepping stone to our analysis is the following key result.

\begin{theorem}[Finite sample risk and efficiency of optimally weighted distributed ridge for fixed regularization parameters]
\label{re_ridge}
Consider the distributed ridge regression problem described above. Suppose we have a dataset with $n$ datapoints (samples), each with an outcome and $p$ features. The dataset is distributed across $k$ sites. Each site has a subset $X_i,Y_i$ of the data, with the $n_i\times p$ matrix $X_i$ of features of $n_i$ samples, and the corresponding outcomes $Y_i$. We compute the local ridge regression estimator $\hat \beta_i(\lambda_i) =(X_i^\top X_i+n_i\lambda_iI_p)^{-1}X_i^\top Y_i$ with fixed regularization parameters $\lambda_i>0$ on each dataset. We send the local estimates to a central location, and combine them via a weighted sum, i.e., $\hbeta_{dist}(w)=\sum_{i=1}^kw_i\hbeta_i$.

Under the linear regression model \eqref{lrm}, the optimal weights that minimize the mean squared error of the distributed estimator are 
$$w^\ast=(A+R)^{-1}v,$$
where the quantities $v,A,R$ are defined below.
\benum
\item  $v$ is a $k$-dimensional vector with $i$-th coordinate $\beta^\top Q_i\beta$, and $Q_i$ are the $p\times p$ matrices $Q_i=(\hSigma_i+\lambda_iI_p)^{-1}\hSigma_i$, 
\item $A$ is a $k\times k$ matrix with $(i,j)$-th entry $\beta^\top Q_iQ_j\beta$, and
\item $R$ is a $k\times k$ diagonal matrix with $i$-th diagonal entry $n_i^{-1}\sigma^2\tr[(\hSigma_i+\lambda_iI_p)^{-2}\hSigma_i]$.
\eenum

The mean squared error of the \emph{optimally weighted} distributed ridge regression estimator $\hbeta_{dist}$ with $k$ sites equals
\beqs
\textnormal{MSE}_{dist}^*(k)
=\E\|\hbeta_{dist}(w^*)-\beta\|^2 
= \|\beta\|^2-v^\top(A+R)^{-1}v,
\eeqs
\end{theorem}

See Section \ref{pf:re_ridge} for the proof. The argument proceeds via a direct calculation, recognizing that \emph{finding the optimal weights for combining the local estimators $\hbeta_i$ can be viewed as a $k$-parameter regression problem} of $\beta$ on $\hbeta_i$, for $i=1,\ldots, k$. 

This result quantifies the mean squared error of the optimally weighted distributed ridge estimator \emph{for fixed regularization parameters $\lambda_i$}. Later we will study how to choose the regularization parameters optimally. The result also gives an exact formula for the optimal weights. 
However, the optimal weights depend on the unknown regression coefficients $\beta$, and are thus not directly usable in practice. Instead, our approach is to make stronger assumptions on $\beta$ under which we can develop estimators for the weights.


{\bf Computational efficiency.} We take a short detour here to discuss computational efficiency. Here by computational efficiency we mean the total time consumption. Computing one ridge regression estimator $(X^\top X + \lambda I_p)^{-1}X^\top Y$ for a fixed regularization parameter $\lambda$ and $n\times p$ design matrix $X$ can be done in time $O(np \min(n,p))$ by first computing the SVD of $X$. This automatically gives the ridge estimator for all values of $\lambda$. 

How much time can we save by distributing the data? Suppose first that $n\ge p$, in which case the total time consumption is $O(np^2)$. Computing ridge locally on the $i$-th machine takes $O(n_ip \min(n_i,p))$ time. Suppose next that we distribute equally to $k$ of machines, and we also have $n_i = n/k \ge p$. Then the time consumption is reduced to $O((n/k) p^2) = O(np^2/k)$. In this case we can say that the total time consumption decreases \emph{proportionally to the number of machines}. This shows the benefit of parallel data processing. 

On the other extreme, if $n \le p$, then $n_i=n/k\le p$, the total time consumption is reduced from $O(n^2 p)$ to $O((n/k)^2 p) = O(n^2p/k^2)$. This shows that the total time consumption decreases \emph{quadratically} in the number of machines (albeit of course the constant is much worse). If we are in an intermediate case where $n\ge p$ and  $n_i=n/k\le p$, then the time decreases at a rate between linear and quadratic.

\subsection{Addressing reader concerns}
At this stage, our readers may have several concerns about our approach. We address some concerns in turn below. 
\begin{enumerate}
\item Does it make sense to average ridge estimators, which can be biased?

A possible concern is that we are working with biased estimators. Would it make sense to debias them first, before weighting? A similar approach has been used for sparse regression, with the debiased Lasso estimators \citep{lee2017communication,battey2018distributed}. However, our results \emph{allow the regularization parameters to be arbitrarily close to zero}, which leads to least squares estimators, with an inverse or pseudoinverse $(X_i^\top X_i)^\dagger$. These are the ``natural" debiasing estimators for ridge regression. For OLS, these are exactly unbiased, while for pseudoinverse, they are approximately so. Hence our approach allows nearly unbiased estimators, and  we automatically discover when this is the optimal method. 

\item Is it possible to improve the weighted sum of local ridge estimators $\hbeta_i$ in trivial ways? 

One-shot weighting is merely a heuristic. If it were possible to improve it in a simple way, then it would make sense to study those methods instead of weighting. However, we are not aware of such methods. For instance, one possibility is to try and add the constant vector into the regression on the global parameter server, because this may help reduce the bias. In simulation studies, we have observed that this approach does not usually lead to a perceptible decrease in MSE.  Specifically we have found that under the simulation setting common throughout the paper, the MSEs with and without a constant term are close (see Section \ref{add_const} for details).

\end{enumerate}

\section{Asymptotics under linear random-effects models}
\label{gen_cov}

The finite sample results obtained so far can be hard to interpret, and do not allow us to directly understand the performance of the optimal one-shot distributed estimator. Therefore, we will consider an asymptotic setting that leads to more insightful results.

Recall that our basic linear model is $Y=X\beta+\ep$, where the error $\ep$ is random. Next, we also assume that a \emph{random-effects model} holds. We assume $\beta$ is random---independently of $\ep$---with coordinates that are themselves independent random variables with mean zero and variance $p^{-1}\sigma^2\alpha^2$. Thus, each feature contributes a small random amount to the outcome. Ridge regression is designed to work well in such a setting, and has several optimality properties in variants of this model. The parameters are now $\theta=(\sigma^2, \alpha^2)$: the \emph{noise level} $\sigma^2$ and the \emph{signal-to-noise ratio} $\alpha^2$ respectively. This parametrization is standard and widely used (e.g. \cite{searle2009variance,dicker2017variance,dobriban2018high}). 

To get more insight into the performance of ridge regression in a distributed environment, we will take an asymptotic approach. Notice from Theorem \ref{re_ridge} that the mean squared error depends on the data only through simple functionals of the sample covariance matrices $\hSigma$ and $\hSigma_i$, such as
$$\beta^\top (\hSigma_i+\lambda_iI_p)^{-1}\hSigma_i \beta,\,\,\,\,\,\,
\beta^\top(\hSigma_i+\lambda_iI_p)^{-1}\hSigma_i 
(\hSigma_j+\lambda_jI_p)^{-1}\hSigma_j \beta,\,\,\,\,\,\,
\tr[(\hSigma_i+\lambda_iI_p)^{-2}\hSigma_i]
.$$
When the coordinates of $\beta$ are iid, the means of the quadratic functionals become proportional to the \emph{traces} of functions of the sample covariance matrices. This motivates us to adopt models from \emph{asymptotic random matrix theory}, where the asymptotics of such quantities are a central topic.

We begin by introducing some key concepts from random matrix theory (RMT) which will be used in our analysis. We will focus on ''Marchenko-Pastur'' (MP) type sample covariance matrices, which are fundamental and popular in statistics (see e.g., \cite{bai2009spectral, anderson1958introduction, paul2014random,yao2015large}).  A key concept is the spectral distribution, which for a $p\times p$ symmetric matrix $A$ is the distribution $F_A$ that places equal mass on all eigenvalues $\lambda_i(A)$ of $\Sigma$. This has cumulative distribution function (CDF) $F_A(x)=p^{-1}\sum_{i=1}^p\mathbf{1}(\lambda_i(A)\leq x)$. A central result in the area is the Marchenko-Pastur theorem, which states that eigenvalue distributions of sample covariance matrices converge \citep{marchenko1967distribution,bai2009spectral}. We state the required assumptions below: 

\begin{assu}
\label{mp}
Consider the following conditions:
\benum
\item The $n\times p$ design matrix $X$ is generated as $X=Z\Sigma^{1/2}$ for an $n\times p$ matrix $Z$ with i.i.d. entries (viewed as coming from an infinite array), satisfying $\E[Z_{ij}]=0$ and $\E[Z_{ij}^2]=1$, and a deterministic $p\times p$ positive semidefinite population covariance matrix $\Sigma$.

\item The sample size $n$ grows to infinity proportionally with the dimension $p$, i.e. $n,p\to\infty$ and $p/n\to\gamma\in(0,\infty)$.

\item The sequence of spectral distributions $F_{\Sigma}:=F_{\Sigma,n,p}$ of $\Sigma:=\Sigma_{n,p}$ converges weakly to a limiting distribution $H$ supported on $[0,\infty)$, called the population spectral distribution.

\eenum
\end{assu}

Then, the Marchenko-Pastur theorem states that with probability $1$, the spectral distribution $F_{\hSigma}$ of the sample covariance matrix $\hSigma$ also converges weakly (in distribution) to a limiting distribution $F_\gamma:=F_\gamma(H)$ supported on $[0,\infty)$ \citep{marchenko1967distribution,bai2009spectral}. The limiting distribution is determined uniquely by a fixed-point equation for its \emph{Stieltjes transform}, which is defined for any distribution $G$ supported on $[0,\infty)$ as
\beqs
m_G(z):=\int_0^\infty\frac{1}{t-z}dG(t),~~~z\in\mathbb C\setminus\mathbb R^+.
\eeqs
With this notation, the Stieltjes transform of the spectral measure of $\hSigma$ satisfies
\beqs
m_{\hSigma}(z)=p^{-1}\tr[(\hSigma-zI_p)^{-1}]\to_{a.s.} m_{F_\gamma}(z),~~~z\in\mathbb C\setminus\mathbb R^+,
\eeqs
where $m_{F_\gamma}(z)$ is the Stieltjes transform of $F$. In addition, we denote by $m'(z)$ the derivative of the Stieltjes transform. Then, it is also known that
\beqs
p^{-1}\tr[(\hSigma-zI_p)^{-2}]\to_{a.s.}m_{F_\gamma}'(z).
\eeqs

The results stated above can be expressed in a different, and perhaps slightly more modern language, using \emph{deterministic equivalents} \citep{serdobolskii2007multiparametric,hachem2007deterministic, couillet2011deterministic, dobriban2018Distributed}. For instance, the Marchenko-Pastur law is a consequence of the following result. For any $z$ where it is well-defined, consider the resolvent $(\hSigma-zI_p)^{-1}$. This random matrix is \emph{equivalent} to a deterministic matrix $(x_p\Sigma-zI_p)^{-1}$ for a certain scalar $x_p=x(\Sigma,n,p,z)$, and we write $$
(\hSigma-zI_p)^{-1}\asymp (x_p\Sigma-zI_p)^{-1}.
$$
Here two sequences of $n\times n$ matrices $A_n, B_n$ (not necessarily symmetric) of growing dimensions are \emph{equivalent}, and we write 
$$A_n \asymp B_n$$ if 
$$\lim_{n\to\infty}\tr\left[C_n(A_n-B_n)\right]=0$$ 
almost surely, for any sequence $C_n$ of $n\times n$ deterministic matrices (not necessarily symmetric) with bounded trace norm, i.e., such that $\lim\sup\|C_n\|_{tr}<\infty$  \citep{dobriban2018Distributed}. Informally, any linear combination of the entries of $A_n$ can be approximated by the entries of $B_n$. This also can be viewed as a kind of \emph{weak convergence} in the matrix space equipped with an inner product (trace). From this, it also follows that the traces of the two matrices are equivalent, from which we can recover the MP law.

In \cite{dobriban2018Distributed}, we collected some useful properties of the calculus of deterministic equivalents. In this work, we use those properties extensively. We also develop and use a new \emph{differentiation rule for the calculus of deterministic equivalents} (see Section \ref{pf:diffcalculus}).

We are now ready to study the asymptotics of the risk. We express the limits of interest in \emph{two equivalent forms}, one in terms of \emph{population quantities} (such as the limiting spectral distribution $H$ of $\Sigma$), and one in terms of \emph{sample quantities} (such as the limiting spectral distribution $F_\gamma$ of $\hSigma$). Moreover, we will denote by $T$ a random variable distributed according to $H$, so that $\E_H g(T)$ denotes the mean of $g(T)$ when $T$ is a random variable distributed according to the limit spectral distribution $H$.

The key to obtaining the results based on population quantities is that the quadratic forms involving $\beta$ have asymptotic equivalents that only depend on $\alpha^2,\sigma^2$, based on the concentration of quadratic forms. Specifically, we have 
$$\beta^\top A \beta \approx \sigma^2\alpha^2/p\cdot\tr(A)$$ 
for suitable matrices $A$ (see the proof of Theorem \ref{ARE} for details). The key to the results based on sample quantities is the MP law and the calculus of deterministic equivalents.

\begin{theorem}[Asymptotics for distributed ridge, arbitrary regularization]
\label{ARE}
In the linear random-effects model under Assumption \ref{mp}, suppose in addition that the eigenvalues of $\Sigma$ are uniformly bounded away from zero and infinity, and that the entries of $Z$ have a finite $8+c$-th moment for some $c>0$. Suppose moreover that the local sample sizes $n_i$ grow proportionally to $p$, so that $p/n_i \to\gamma_i>0$. 

Then the optimal weights for distributed ridge regression, and its MSE, converge to definite limits. Recall from Theorem \ref{re_ridge} that we have the formulas $w^\ast=(A+R)^{-1}v$ and  MSE$^*_{dist}=\|\beta\|^2-v^\top(A+R)^{-1}v$ for the optimal finite sample weights and risk, and thus it is enough to find the limit of $v, A$ and $R$. These have the following limits:
\benum
\item With probability one, we have the convergence $v\to V\in\mathbb R^k$. The $i$-th coordinate of the limit $V$ has the following two equivalent forms, in terms of population and sample quantities, respectively:
\beqs V_i=\sigma^2\alpha^2\E_H\frac{x_iT}{x_iT+\lambda_i}=\sigma^2\alpha^2(1-\lambda_im_{F_{\gamma_i}}(-\lambda_i)). \eeqs
Recall that $H$ is the limiting population spectral distribution of $\Sigma$, and $T$ is a random variable distributed according to $H$. Among the empirical quantities, $F_{\gamma_i}$ is the limiting empirical spectral distribution of $\hSigma_i$ and $x_i:=x_i(H,\lambda_i,\gamma_i)>0$ is the unique solution of the fixed point equation
\beqs 1-x_i=\gamma_i\left[1-\lambda_i\int_0^\infty\frac{dH(t)}{x_it+\lambda_i}\right]=\gamma_i\left[1-\E_H\frac{\lambda_i}{x_iT+\lambda_i}\right]. \eeqs
It is part of the theorem's claim that there is such an $x_i$.

\item With probability one, $A\to \mathcal{A}\in\mathbb R^{k\times k}$. For $i\neq j$, the $(i,j)$-th entry of $\mathcal{A}$ is, in terms of the population spectral distribution $H$,
\beqs \mathcal{A}_{ij}=\sigma^2\alpha^2\E_H\frac{x_ix_jT^2}{(x_iT+\lambda_i)(x_jT+\lambda_j)}. \eeqs
The $i$-th diagonal entry of $\mathcal{A}$ is, in terms of population and sample quantities, respectively,
\begin{align*} 
\mathcal{A}_{ii}
&=\sigma^2\alpha^2\left[1-\E_H\frac{2\lambda_ix_iT+\lambda_i^2}{(x_iT+\lambda_i)^2}+\frac{\lambda_i^2\gamma_i x_i \left(\E_H\frac{T}{(x_iT+\lambda_i)^2}\right)^2}{1 +\gamma_i \lambda_i\E_H\frac{T}{(x_iT+\lambda_i)^2}}\right]\\
&=\sigma^2\alpha^2\left[1-2\lambda_im_{F_{\gamma_i}}(-\lambda_i)+\lambda_i^2m'_{F_{\gamma_i}}(-\lambda_i)\right].
\end{align*}

\item With probability one, the diagonal matrix $R$ converges, $R\to\mathcal{R}\in\mathbb R^{k\times k}$, where of course $\mathcal{R}$ is also diagonal. The $i$-th diagonal entry of $\mathcal{R}$ is, in terms of population and sample quantities, respectively,
\beqs \mathcal{R}_{ii}=\sigma^2\left[\frac{x_i\E_H\frac{T}{(x_iT+\lambda_i)^2}}{1+\lambda_i\gamma_i\E_H\frac{T}{(x_iT+\lambda_i)^2}}
 \right]=\sigma^2\left[\gamma_im_{F_{\gamma_i}}(-\lambda_i)-\gamma_i\lambda_im'_{F_{\gamma_i}}(-\lambda_i)\right]. \eeqs

\eenum
The limiting weights and MSE are then 
$$\mathcal W_k^*=(\mathcal{A}+\mathcal{R})^{-1}V$$
and
$$\mathcal M_k=\sigma^2\alpha^2-V^\top(\mathcal{A}+\mathcal{R})^{-1}V.$$
\end{theorem}

See Section \ref{pf:ARE} for the proof. The statement may look complicated, but the formulas simplify considerably in the uncorrelated case $\Sigma=I_p$, on which we will focus later. Moreover, these limiting formulas are also fundamental for developing consistent estimators for the optimal weights. To develop an algorithm for the practically common general covariance case, the following theorem is crucial. 

\begin{theorem}[Asymptotics for distributed ridge when the samples are equally distributed]
\label{ARE_equal}
Consider the assumptions and the notations of Theorem \ref{ARE}. We further assume the samples are equally distributed across the local machines, i.e. $n_1=n_2=\cdots=n_k=n/k$ and $\gamma_1=\gamma_2=\cdots=\gamma_k=k\gamma$. We use the same tuning parameter $\lambda$ for each local estimator. Then the limiting optimal weights $\mathcal{W}_k^*$ and the limiting MSE $\mathcal{M}_k$ have the following forms:
\begin{equation*}
\mathcal{W}_k^*=(1,1,\dots,1)^\top\cdot\frac{\sigma^2\alpha^2(1-\lambda m)}{\mathcal{F}+k\mathcal{G}}~~\textnormal{and}~~\mathcal{M}_k=\sigma^2\alpha^2-\frac{\sigma^4\alpha^4(1-\lambda m)^2k}{\mathcal{F}+k\mathcal{G}}.
\end{equation*}
Here $\mathcal{F}$ and $\mathcal{G}$ are defined as follows:
\begin{equation*}
\mathcal{F}=\sigma^2\alpha^2\frac{k\gamma\lambda^2(m-\lambda m')^2}{1-k\gamma+k\gamma\lambda m'}+\sigma^2k\gamma(m-\lambda m')
\end{equation*}
and
\begin{equation*}
\mathcal{G}=\sigma^2\alpha^2\left(1-2\lambda m+\lambda^2m'-\frac{k\gamma\lambda^2(m-\lambda m')^2}{1-k\gamma+k\gamma\lambda m'}\right)
\end{equation*}
where $m:=m_{F_{k\gamma}}(-\lambda)$ and $m':=-\frac{dm}{d\lambda}$.
\end{theorem}

See Section \ref{pf:ARE_equal} for the proof. Based on this theorem, we are able to develop an algorithm which works for arbitrary covariance structures. See Section \ref{algo} for the details.

Now we discuss the problem of estimating the optimal weights, which is crucial for developing practical methods. The results in Theorem \ref{ARE_equal} show that to estimate the weights consistently, if the tuning parameter $\lambda$ is known, we \emph{only need to estimate $\alpha^2,\sigma^2$} consistently. The reason is that we can use $\tr(\hSigma_i+\lambda I)^{-1}/p$ to approximate $m$, and use $\tr(\hSigma_i+\lambda I)^{-2}/p$ to approximate $m'$.

Estimating these two parameters is a well-known problem, and several approaches have been proposed, for instance restricted maximum likelihood (REML) estimators  \citep{jiang1996reml,searle2009variance, dicker2014variance, dicker2016variance, jiang2016high}, etc. We can use---for instance---results from \cite{dicker2017variance}, who showed that the Gaussian MLE is consistent and asymptotically efficient for $\theta=(\sigma^2,\alpha^2)$ even in the non-Gaussian setting of this paper (see Section \ref{signoise} for a summary).

\section{Special case: identity covariance}
\label{id_cov}

When the population covariance matrix is the identity, that is $\Sigma = I$, the results simplify considerably. In this case the features are nearly uncorrelated.  It is known that the limiting Stieltjes transform $m_{F_\gamma} :=m_\gamma$ of $\hSigma$ has the explicit form \citep{marchenko1967distribution}:
\beq m_\gamma(z) = \frac{(z+\gamma-1)+\sqrt{(z+\gamma-1)^2-4z \gamma}}
{-2z \gamma}. 
\label{st_form}
\eeq
As usual in the area, we use the principal branch of the square root of complex numbers.

\subsection{Properties of the estimation error and asymptotic relative efficiency}

We can use the closed form expression for the Stieltjes transform to get explicit formulas for the optimal weights. From Theorem \ref{ARE}, we conclude the following simplified result:

\begin{theorem}[Asymptotics for isotropic population covariance, arbitrary regularization parameters]
\label{AREsimple}
In addition to the assumptions of Theorem \ref{ARE}, suppose that the population covariance matrix $\Sigma=I$. Then the limits of $v, A$ and $R$ have simple explicit forms:
\benum
\item The $i$-th coordinate of $V$ is:
\beqs V_i=\sigma^2\alpha^2[1-\lambda_im_{\gamma_i}(-\lambda_i)], \eeqs
where $m_{\gamma_i}(-\lambda_i)$ is the Stieltjes transform given above in equation \eqref{st_form}.

\item The entries of $\mathcal{A}$ are
$$
\mathcal{A}_{ij}=
\begin{cases}
\sigma^2\alpha^2[1- \lambda_im_{\gamma_i}(-\lambda_i)]
\cdot 
[1- \lambda_jm_{\gamma_j}(-\lambda_j)],\,\,& \text{for } i\neq j\\
\sigma^2\alpha^2\left[1-2\lambda_im_{\gamma_i}(-\lambda_i)
+\lambda_i^2 m'_{\gamma_i}(-\lambda_i)\right],\,\,& \text{for } i=j.
\end{cases}
$$

\item The $i$-th diagonal entry of $\mathcal{R}$ is 
\beqs \mathcal{R}_{ii}=
\sigma^2\gamma_i\left[m_{\gamma_i}(-\lambda_i)-\lambda_im'_{\gamma_i}(-\lambda_i)\right]. \eeqs

\eenum
The limiting optimal weights for combining the local ridge estimators are $\mathcal W_k^*=(\mathcal{A}+\mathcal{R})^{-1}V$, and MSE of the optimally weighted distributed estimator is
\beqs
\mathcal{M}_k
=\frac{\sigma^2\alpha^2}{1+\sum_{i=1}^k\frac{V_i^2}
{\sigma^2\alpha^2(\mathcal{R}_{ii} + \mathcal{A}_{ii})-V_i^2}}.
\eeqs

\end{theorem}
See Section \ref{pf:AREsimple} for the proof. This theorem shows the surprising fact that the limiting risk \emph{decouples} over the different machines. By this we mean that the limiting risk can be written in a simple form, involving a sum of terms depending on each machine, without any interaction. This seems like a major surprise.

To explain in more detail the decoupling phenomenon, let us study how the local risks are related to the distributed risks. Define $V = V( \gamma, \lambda)$ to be the limiting scalar $V \in \R$ defined above, in the special case $k=1$. Explicitly, this is the limit of the quantity $\beta^\top Q \beta$, where $Q=(\hSigma+\lambda I_p)^{-1}\hSigma$, as given in Theorem \ref{re_ridge} applied for $k=1$. Let $D$ be the scalar expression $D(\gamma, \lambda) = \sigma^2\alpha^2(\mathcal{R} + \mathcal{A})-V$ when $k=1$. With these notations, the risk $\mathcal{M}_1$ of ridge regression when computed on the entire dataset equals
\begin{align*}
\mathcal{M}_1(\gamma, \lambda)
= \frac{\sigma^2\alpha^2}{1+\frac{V( \gamma, \lambda)}{D(\gamma, \lambda)}}.
\end{align*}
Moreover, the risk of optimally weighted one-shot distributed ridge over $k$ subsets, with arbitrary regularization parameters $\lambda_i$, equals
\begin{align*}
\mathcal{M}_k(\gamma_1,\ldots, \gamma_k, \lambda_1,\ldots,\lambda_k)
= \frac{\sigma^2\alpha^2}
{1+\sum_{i=1}^k\frac{V_i^2( \gamma_i, \lambda_i)}
{D_i(\gamma_i, \lambda_i)}}.
\end{align*}
Then one can check that we have the following equations connecting the risk computed on the entire dataset and the distributed risk:
\begin{align*}
\frac{\sigma^2\alpha^2}
{\mathcal{M}_k( \gamma_1,\ldots, \gamma_k, \lambda_1,\ldots,\lambda_k)} -1
&= 
\sum_{i=1}^k\frac{\sigma^2\alpha^2}{\mathcal{M}_1( \gamma_i, \lambda_i)} -k, \\ 
\mathcal{M}_k( \gamma_1,\ldots, \gamma_k, \lambda_1,\ldots,\lambda_k)
&= \frac{1}
{\sum_{i=1}^k\frac{1}{\mathcal{M}_1( \gamma_i, \lambda_i)} 
+\frac{1-k}{\sigma^2\alpha^2}}.
\end{align*}
These equations are precisely what we mean by \emph{decoupling}. The distributed risk can be written as a function of the type $1/(\sum_i1/x_i+b)$ of the distributed risks. Therefore, there are no ``interactions" between the different risk functions. Similar expressions have been obtained for linear regression \citep{dobriban2018Distributed}.

Next, we discuss in more depth why the limiting risk decouples. Mathematically, the key reason is that when $\Sigma=I$, the limit of $A_{ij}$ for $i\neq j$ decouples into a product of two terms. Therefore, the distributed risk function involves a quadratic form with zero \emph{off-diagonal} terms. This is not the case for general population covariance $\Sigma$. We provide an explanation via free probability theory in Section \ref{free_prob_exp}.

An important consequence of the decoupling is that \emph{we can optimize the individual risks over the tuning parameters $\lambda_i$ separately}.

\begin{proposition}[Optimal regularization (tuning) parameters, and risk]
\label{parameter}
Under the assumptions of Theorem \ref{AREsimple}, the optimal regularization (tuning) parameters $\lambda_i$ that minimize the local MSEs also minimize the distributed risk $\mathcal{M}_k$. They have the form
\beqs \lambda_i=\frac{\gamma_i}{\alpha^2},~~i=1,2,\dots,k. \eeqs
Moreover, the risk $\mathcal{M}_k$ of distributed ridge regression with optimally tuned regularization parameters is
\beqs \mathcal{M}_k=\frac{\sigma^2\alpha^2}
{1+\sum_{i=1}^k\left[\frac{\alpha^2}{\gamma_im_{\gamma_i}(-\gamma_i/\alpha^2)}-1\right]},
\eeqs
\end{proposition}
See Section \ref{pf:parameter} for the proof. 

The main goal of our paper is to study the behavior of the one-shot distributed ridge estimator and compare it with the centralized estimator. It is helpful to first understand the properties of the \emph{optimal risk} function $\phi(\gamma):=\gamma m_\gamma(-\gamma/\alpha^2)$. The optimal risk function equals the optimally tuned global risk $\mathcal{M}_1$ up to a factor $\sigma^2$. It has the explicit form
$$\phi(\gamma)=\gamma m_{\gamma}(-\gamma/\alpha^2) 
=\frac{-\gamma/\alpha^2+\gamma-1+\sqrt{(-\gamma/\alpha^2+\gamma-1)^2+4\gamma^2/\alpha^2}}{2\gamma/\alpha^2}.$$

\begin{proposition}[Properties of the optimal risk function]
\label{optimalrisk}
The optimal risk function $\phi(\gamma)$ has the following properties:
\benum
\item {\bf Monotonicity:} $\phi(\gamma)$ is an increasing function of $\gamma\in[0,\infty)$ with $\lim_{\gamma\rightarrow0_+}\phi(\gamma)=0$ and $\lim_{\gamma\rightarrow+\infty}\phi(\gamma)=\alpha^2$.

\item {\bf Concavity:} When $\alpha\leq1,\phi(\gamma)$ is a concave function of $\gamma\in[0,\infty)$. When $\alpha>1,$ $\phi(\gamma)$ is convex for small $\gamma$ (close to $0$), and concave for large $\gamma$.

\eenum

\end{proposition}

\begin{figure}
\begin{subfigure}{.5\textwidth}
  \centering
\includegraphics[scale=0.4]
{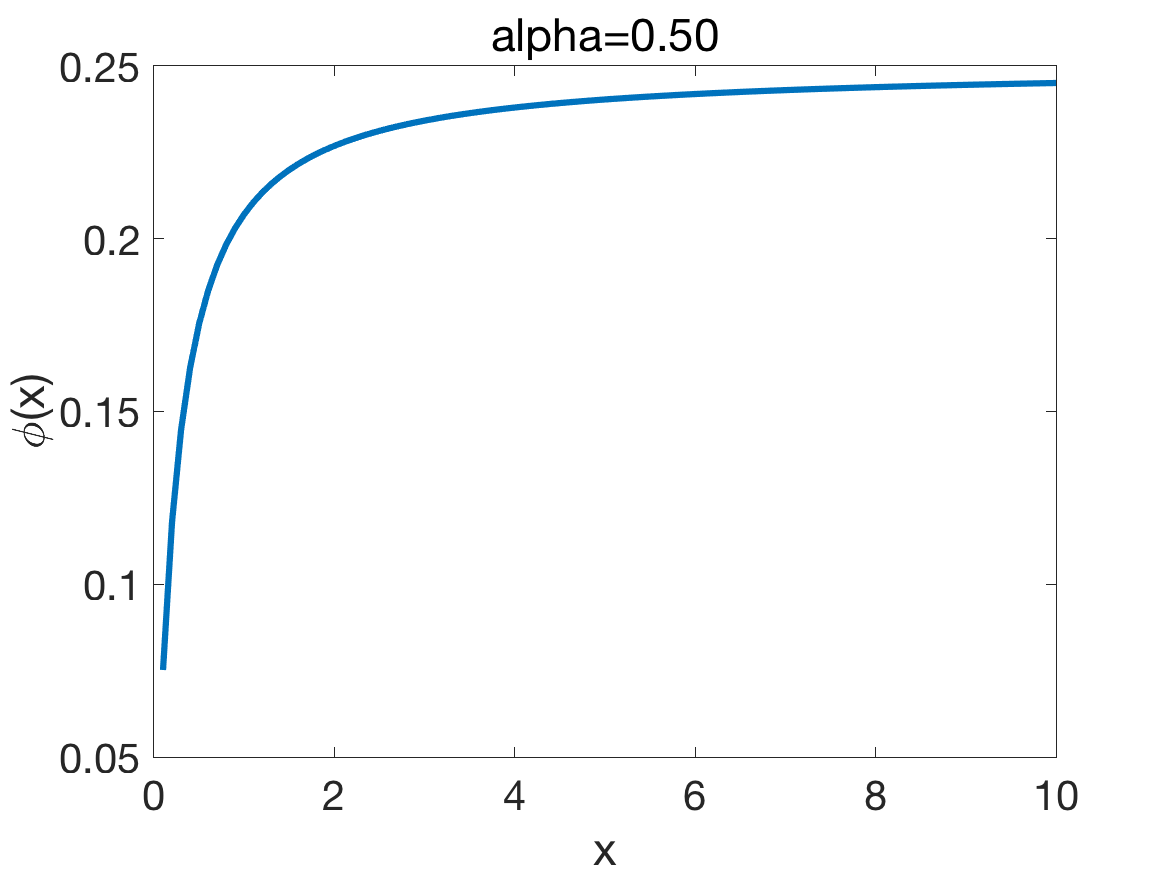}
\end{subfigure}
\begin{subfigure}{.5\textwidth}
  \centering
\includegraphics[scale=0.4]
{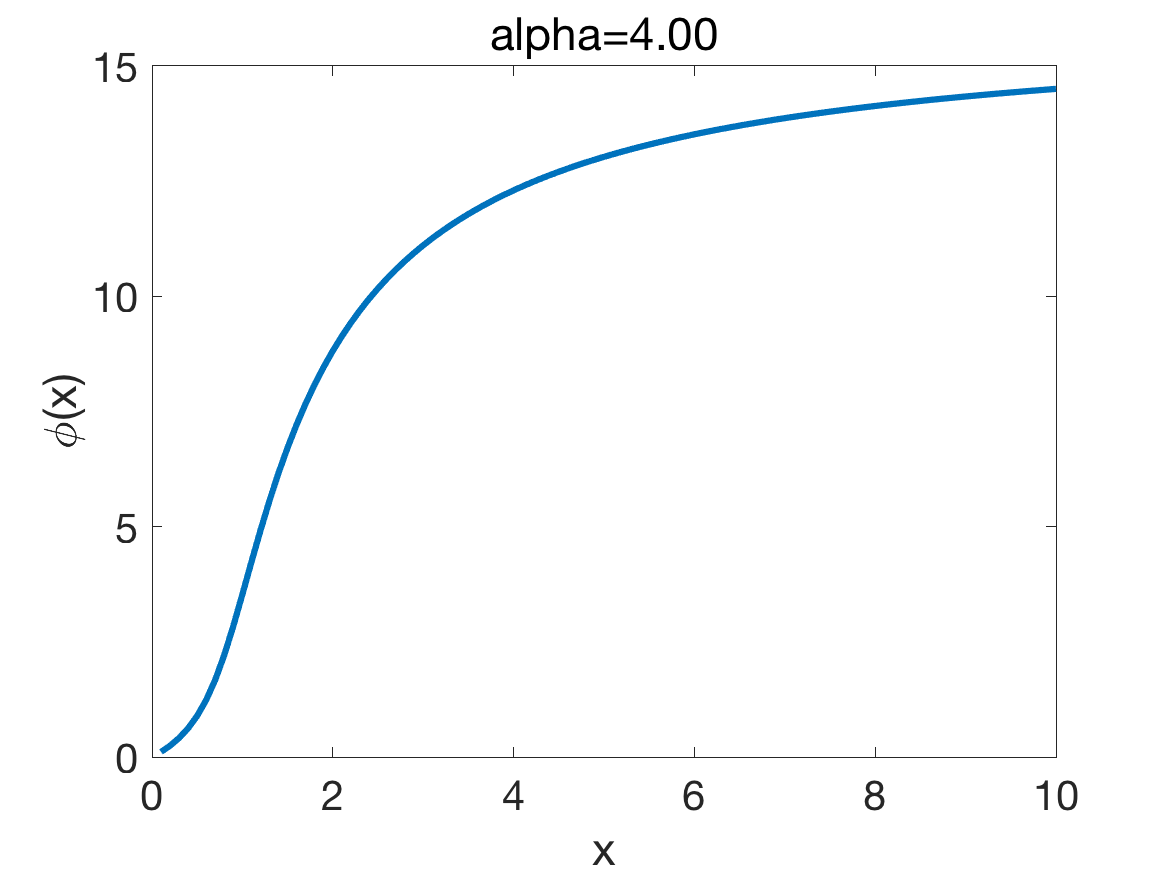}
\end{subfigure}
\caption{Plots of the optimal risk function $\phi$ as a function of the aspect ratio $\gamma$ (denoted by $x$ in the plots), for different signal strength parameters $\alpha$.}
\label{phi}
\end{figure}

See Section \ref{pf:optimalrisk} for the proof. See also Figure \ref{phi} for plots of $\phi$ for different $\alpha$, which show its monotonicity and convexity properties. The aspect ratio $\gamma$ characterizes the dimensionality of the problem. It makes sense that $\phi(\gamma)$ is increasing, since the regression problem should become more difficult as the dimension increases. For the second property, the concavity of the function means that it grows very fast to approach its limit. When the signal-to-noise ratio $\alpha^2$ is small, the risk is concave, so it grows fast with the dimension. But when the signal-to-noise ratio becomes large, the risk will grow much slower at the beginning. Here the phase transition happens at $\alpha^2=1$. This gives insight into the effect of the signal-to-noise ratio on the regression problem.

To compare the distributed and centralized estimators, we will study their (asymptotic) relative efficiency (ARE), which is the (limit of the) ratio of their mean squared errors. Here we assume each estimator is optimally tuned. This quantity, which is at most unity, captures the loss of efficiency due to the distributed setting. An ARE close to $1$ is ``good", while an ARE close to $0$ is ``bad". From the results above, it follows that the ARE has the form
\beqs
ARE=\frac{\mathcal{M}_1}{\mathcal{M}_k}=\frac{\gamma m_{\gamma}(-\gamma/\alpha^2)}{\alpha^2}\left[ 1+\sum_{i=1}^k\left( \frac{\alpha^2}{\gamma_im_{\gamma_i}(-\gamma_i/\alpha^2)}-1  \right)\right ] \le 1.
\eeqs

We have the following properties of the ARE.
\begin{theorem}[Properties of the asymptotic relative efficiency (ARE)]
\label{AREprop}
The asymptotic relative efficiency (ARE) has the following properties:
\benum

\item {\bf Worst case is equally distributed data}: For fixed $k,\alpha^2$ and $\gamma$, the ARE attains its minimum when the samples are \emph{equally distributed} across $k$ machines, i.e. $\gamma_1=\gamma_2=\cdots=\gamma_k=k\gamma$. We denote the minimal value by $\psi(k,\gamma,\alpha^2)$. That is
\beqs \min_{\gamma_1,\ldots,\gamma_k}ARE=\psi(k,\gamma,\alpha^2):=\frac{\gamma m_{\gamma}(-\gamma/\alpha^2)}{\alpha^2}\left(1-k+\frac{\alpha^2}{\gamma m_{k\gamma}(-k\gamma/\alpha^2)} \right ). \eeqs

\item {\bf Adding more machines leads to efficiency loss}: For fixed $\alpha^2$ and $\gamma$, $\psi(k,\gamma,\alpha^2)$ is a decreasing function on $k\in[1,\infty)$ with 
$\lim_{k\to1_+}\psi(k,\gamma,\alpha^2)=1$ and \emph{infinite-worker limit}

$$\lim_{k\to\infty}\psi(k,\gamma,\alpha^2)=h(\alpha^2,\gamma) <1.$$ 

Here we can view $\psi$ as a continuous function of $k$ for convenience, although originally it is only well-defined for $k\in\mathbb N$. We emphasize that the infinite-worker limit tells us how much efficiency we have for a very large number of machines. It is a nontrivial result that this quantity is strictly positive.

\item {\bf Form of the infinite-worker limit}: As a function of $\alpha^2$ and $\gamma$, $h(\alpha^2,\gamma)$ has the explicit form
\beqs
h(\alpha^2,\gamma)=\frac{-\gamma/\alpha^2+\gamma-1+\sqrt{(-\gamma/\alpha^2+\gamma-1)^2+4\gamma^2/\alpha^2}}{2\gamma}\left(1+\frac{\alpha^2}{\gamma(1+\alpha^2)}\right).
\eeqs

\item {\bf Edge cases of the infinite-worker limit}: For fixed $\alpha^2$, $h(\alpha^2,\gamma)$ is an increasing function of $\gamma\in[0,\infty)$ with limit
\beqs
\lim_{\gamma\to 0}h(\alpha^2,\gamma)=\frac{1}{1+\alpha^2},~~\lim_{\gamma\to\infty}h(\alpha^2,\gamma)=1.
\eeqs
On the other hand, for fixed $\gamma$, $h(\alpha^2,\gamma)$ is a decreasing function of $\alpha^2\in[0,\infty)$ with limit
\beqs
\lim_{\alpha^2\to0}h(\alpha^2, \gamma)=1,~~
\lim_{\alpha^2\rightarrow\infty}h(\alpha^2, \gamma)
=
\begin{cases}
1-\frac{1}{\gamma^2},~~\gamma>1,\\
0,~~0<\gamma\leq1.
\end{cases}
\eeqs

\eenum

\end{theorem}

\begin{figure}
\begin{subfigure}{.5\textwidth}
  \centering
\includegraphics[scale=0.4]
{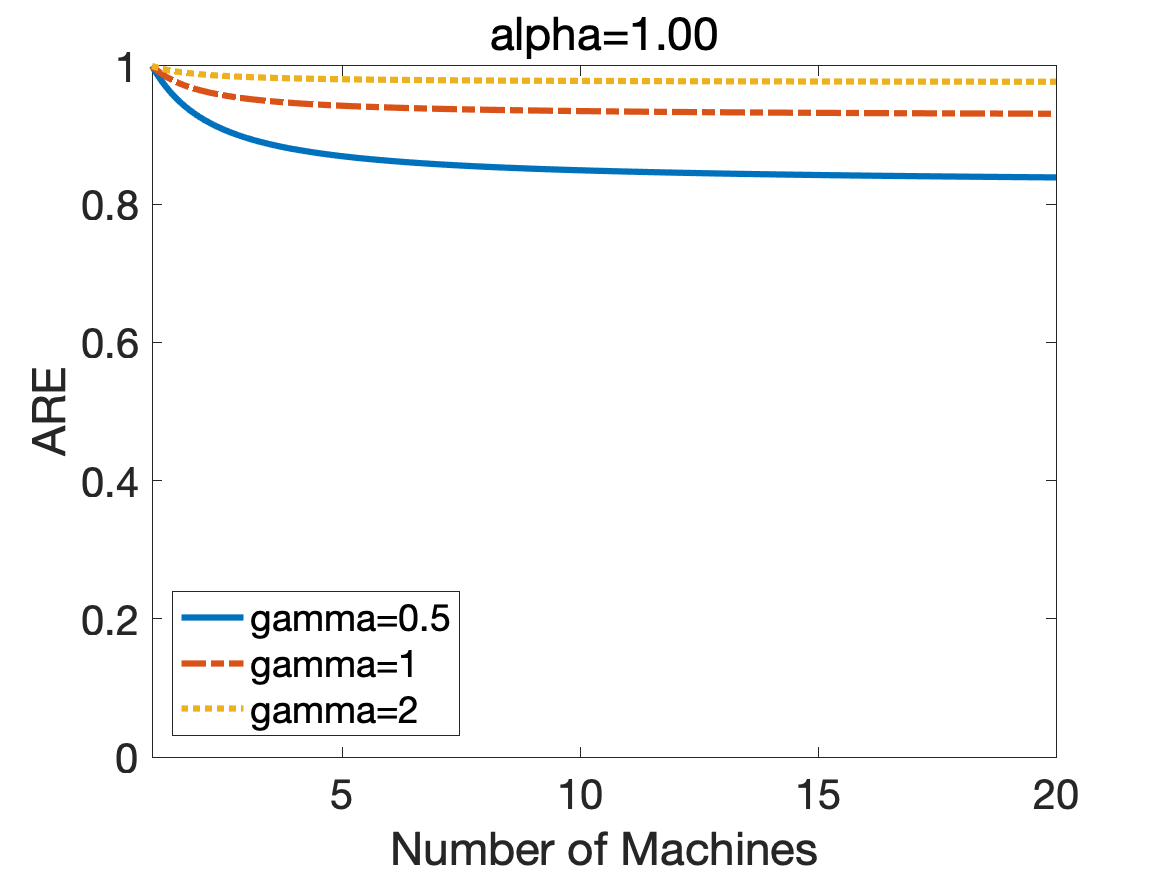}
\end{subfigure}
\begin{subfigure}{.5\textwidth}
  \centering
\includegraphics[scale=0.4]
{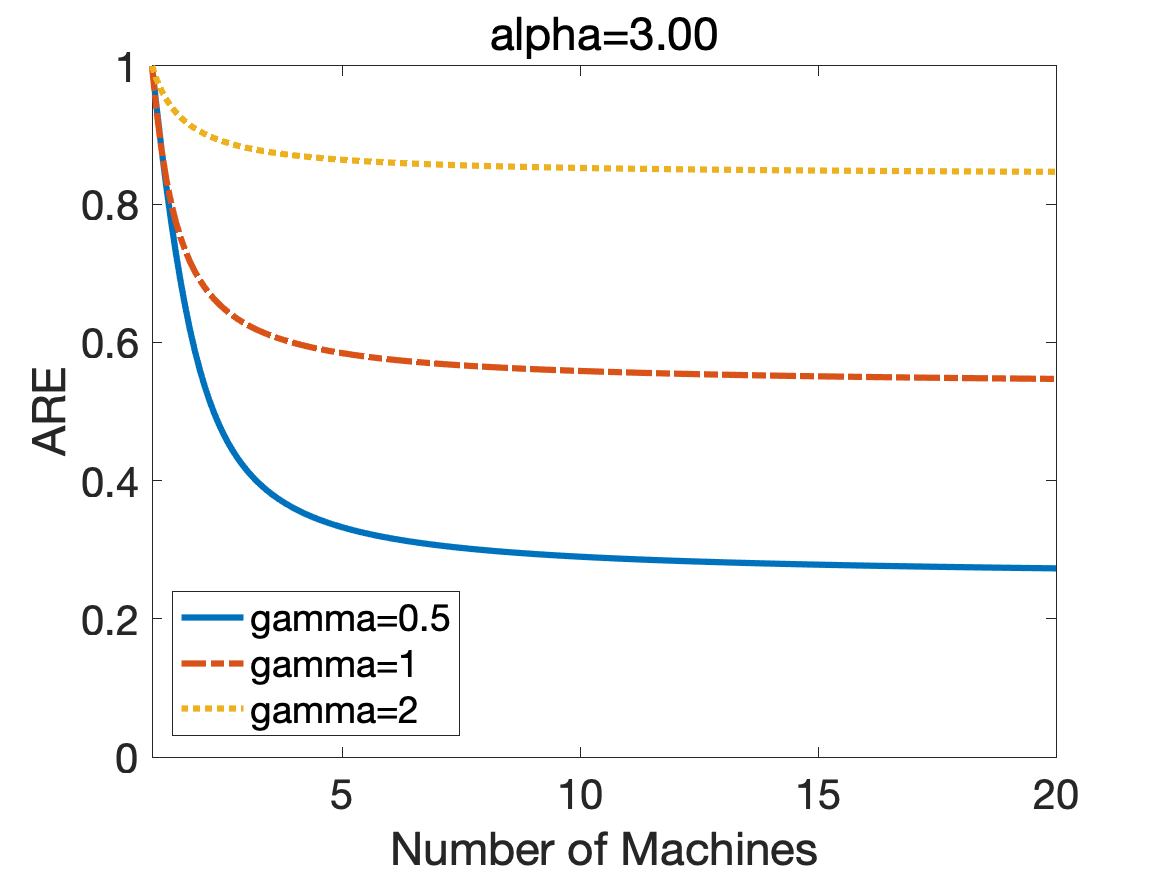}
\end{subfigure}
\caption{Plots of the asymptotic relative efficiency $\psi$ when the datasets are evenly distributed, for different $\alpha$ and $\gamma$. See Theorem \ref{AREprop} for the properties of the ARE.}
\label{ridgeARE}
\end{figure}

See Section \ref{pf:AREprop} for the proof. See Figure \ref{ridgeARE} for some plots of the evenly distributed ARE $\psi$ for various $\alpha$ and $\gamma$ and Figure \ref{arefig} for the surface and contour plots of $h(\alpha^2,\gamma)$. The efficiency loss tends to be larger (ARE is smaller) when the signal-to-noise ratio $\alpha^2$ is larger. The plots confirm the theoretical result that the efficiency always decreases with the number of machines. Relatively speaking, the distributed problem becomes easier and easier as the dimension increases, compared to the aggregated problem (i.e., the ARE increases in $\gamma$ for fixed parameters). This can be viewed as a blessing of dimensionality. 

We also observe a nontrivial \emph{infinite-worker limit}. Even in the limit of many machines, distributed ridge \emph{does not lose all efficiency}. This is in contrast to doing linear regression on each machine, where all efficiency is lost when the local sample sizes are less than the dimension \citep{dobriban2018Distributed}. This is one of the few results in the distributed learning literature where one-step weighting gives nontrivial results for \emph{arbitrary large} $k$, i.e., we can take $k \to \infty$ and we still obtain nontrivial results. We find this quite remarkable.

Overall, the ARE is generally large, \emph{except} when $\gamma$ is small and $\alpha$ is large. This is a setting with strong signal and relatively low dimension, which is also the ``easiest" setting from a statistical point of view. In this case, perhaps we should use other techniques for distributed estimation, such as iterative methods.

\subsection{Properties of the optimal weights}

Next, we study properties of the optimal weights. This is important, because choosing them is a crucial practical question. The literature on distributed regression typically considers simple averages of local estimators, for which  $\hbeta_{dist}=k^{-1}\sum_{i=1}^k\hbeta_i$ (see, e.g. \cite{zhang2015divide, lee2017communication, battey2018distributed}). In contrast, we will find that the optimal weights \emph{do not sum up to unity}.

Formally, we have the following properties of the optimal weights.
\begin{theorem}[Properties of the optimal weights]
\label{weights}
The asymptotically optimal weights $\mathcal W_k^*=(\mathcal{A}+\mathcal{R})^{-1}V$ have the following properties:
\benum
\item  {\bf Form of the optimal weights}: The $i$-th coordinate of $\mathcal W_k$ is:
$$
\mathcal W_{k,i}=\left(\frac{\alpha^2}{\gamma_im_{\gamma_i}(-\gamma_i/\alpha^2)}\right)\cdot \left(\frac{1}{1+\sum_{i=1}^k\left[\frac{\alpha^2}{\gamma_im_{\gamma_i}(-\gamma_i/\alpha^2)}-1\right]}\right),
$$
and the sum of the limiting weights is always greater than or equal to one: $\sum_{i=1}^k \mathcal W_{k,i} \ge 1$. When $k\ge 2$, the sum is strictly greater than one: 
$$\sum_{i=1}^k \mathcal W_{k,i} > 1.$$

\item {\bf Evenly distributed optimal weights}: When the samples are evenly distributed, so that all limiting aspect ratios $\gamma_i$ are equal, $\gamma_i=k\gamma$, then all $\mathcal W_{k,i}$ equal the \emph{optimal weight function} $\mathcal{W}(k,\gamma,\alpha^2)$, which has the form
$$
\mathcal{W}(k,\gamma,\alpha^2)=\frac{\alpha^2}{\alpha^2k+(1-k)k\gamma\cdot m_{k\gamma}(-k\gamma/\alpha^2)}.
$$
This can also be written in terms of the \emph{optimal risk function} $\phi(\gamma,\alpha^2)$ defined above as
$$
\mathcal{W}(k,\gamma,\alpha^2)=\frac{\alpha^2}{\alpha^2k-(k-1)\phi(k\gamma,\alpha^2)}.
$$

\item {\bf Limiting cases}: For fixed $k$ and $\alpha^2$, the optimal weight function $\mathcal{W}(k,\gamma,\alpha^2)$ is an increasing function of $\gamma\in[0,\infty)$ with $\lim_{\gamma\to0_+}\mathcal{W}(\gamma)=1/k$ and $\lim_{\gamma\to\infty}\mathcal{W}(\gamma)=1$.

\eenum
\end{theorem}

\begin{figure}
\begin{subfigure}{.5\textwidth}
  \centering
\includegraphics[scale=0.5]
{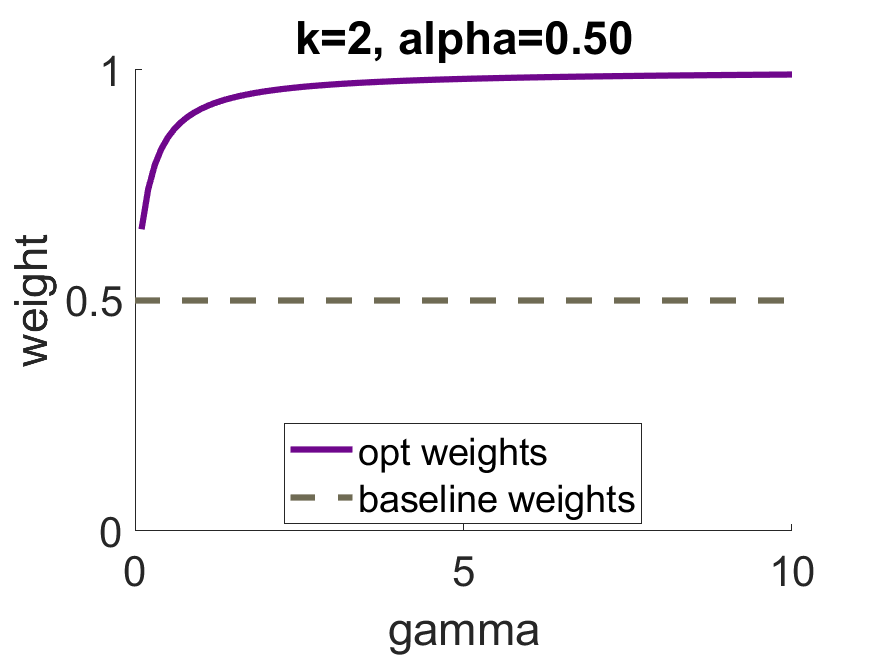}
\end{subfigure}
\begin{subfigure}{.5\textwidth}
  \centering
\includegraphics[scale=0.5]
{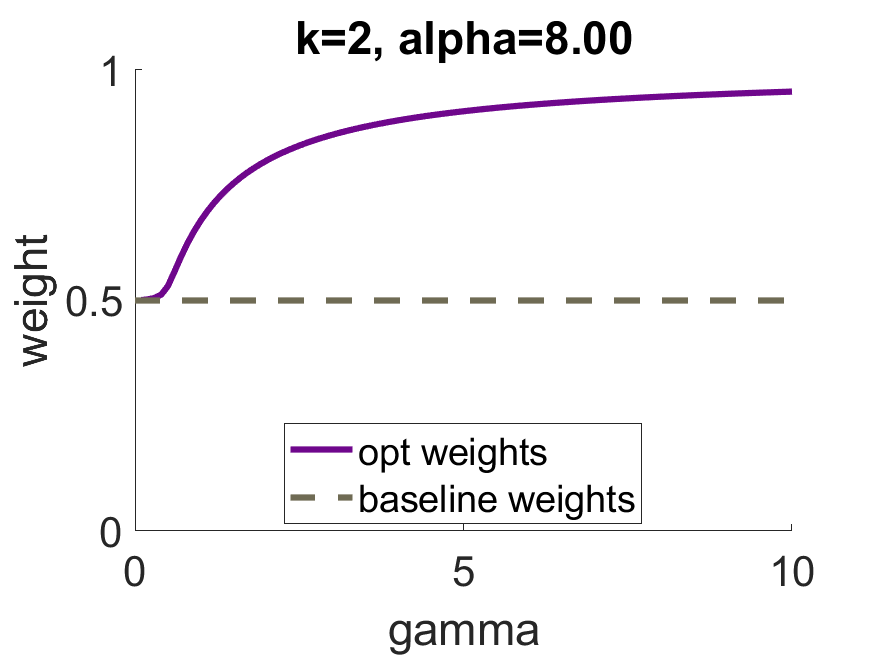}
\end{subfigure}
\caption{Plots of optimal weights for different $\alpha$.}
\label{optw}
\end{figure}

\begin{figure}
\centering
\begin{subfigure}{.45\textwidth}
  \centering
  \includegraphics[scale=0.5]{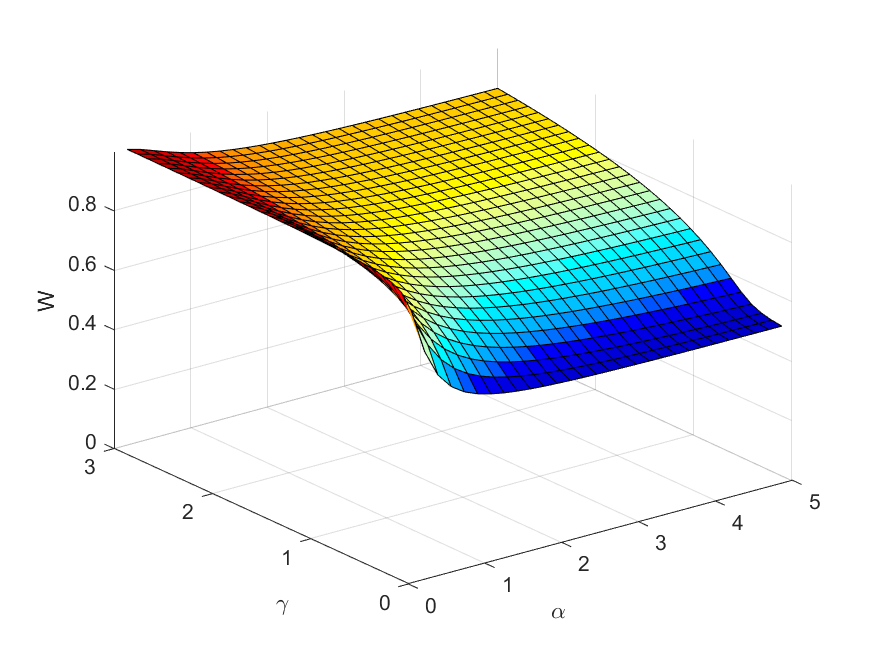}
  \caption{Surface}
\end{subfigure}
\begin{subfigure}{.45\textwidth}
  \centering
  \includegraphics[scale=0.5]{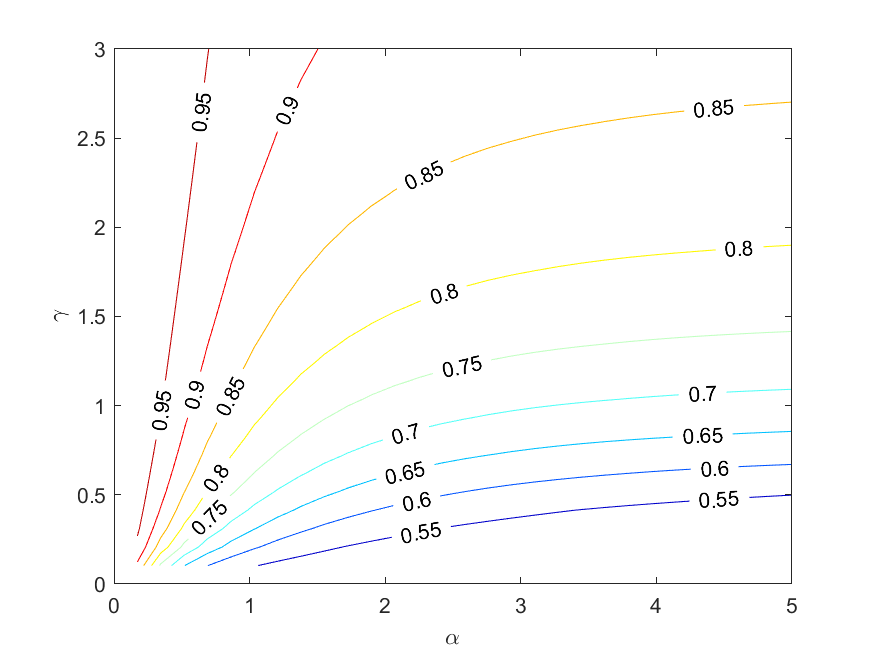}
  \caption{Contour}
\end{subfigure}
\caption{Surface and contour plots of the optimal weights.}
  \label{optw2}
\end{figure}

See Section \ref{pf:weights} for the proof.
See Figures \ref{optw} and \ref{optw2} for some plots of the optimal weight function with $k=2$. We can see that the optimal weights are usually large, and always greater than $1/k$. When the signal-to-noise ratio $\alpha^2$ is small, the weight function is concave and increases fast to approach one. In the low dimensional setting where $\gamma\to0$, the weights tend to the uniform average $1/k$. Hence in this setting we recover the classical uniform averaging methods, which makes sense, because ridge regression with optimal regularization parameter tends to linear regression in this regime.

\begin{figure}
\begin{subfigure}{.5\textwidth}
  \centering
\includegraphics[scale=0.5]
{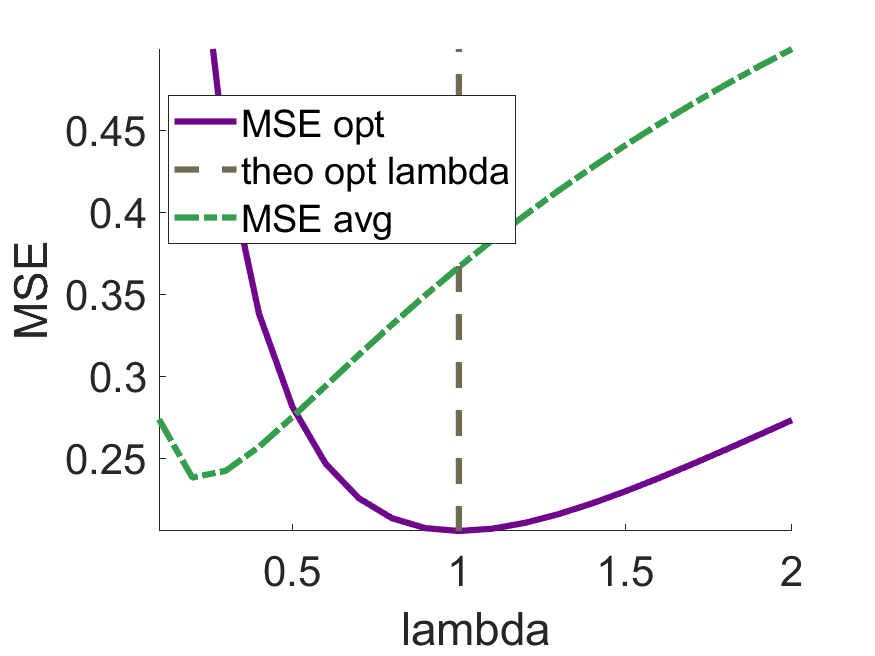}
\end{subfigure}
\begin{subfigure}{.5\textwidth}
  \centering
\includegraphics[scale=0.5]
{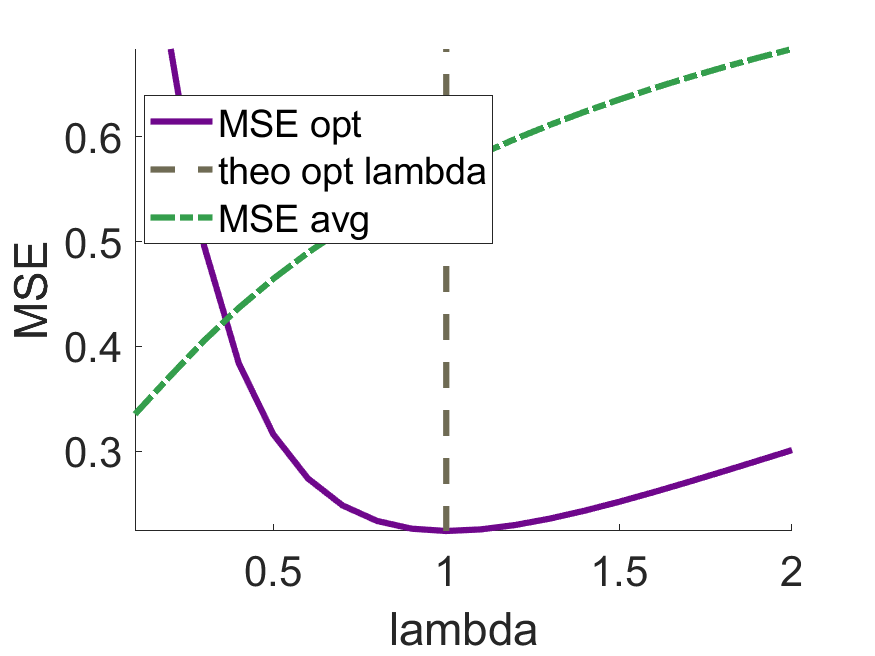}
\end{subfigure}
\caption{Distributed risk as a function of the regularization parameter. We plot both the risk with optimal weights (MSE opt) and the risk obtained from sub-optimal averaging (MSE avg). We set $\alpha=1$, $\gamma=0.17$ and $k=5,10$.}
\label{drl}
\end{figure}

How much does optimal weighting help? It is both interesting and important to know this, especially compared to naive uniform weighting, because it allows us to compare our proposed weighting method to the ``baseline". See Figure \ref{drl}. We have plotted the risk of distributed ridge regression for both the optimally weighted version and the simple average, as a function of the regularization parameter. We observe that \emph{optimal weighting can lead to a 30-40\% decrease in the risk}. Therefore, our proposed weighting scheme can lead to a substantial benefit.

Why are the weights large, and why do they sum to a quantity greater than one? The short intuitive answer is that ridge regression is \emph{negatively} (or \emph{downward}) biased, and so we must \emph{counter the effect of bias by upweighting}. This also can be viewed as a way of \emph{debiasing}. In different contexts, it is already well known that debiasing can play a kew role in distributed learning (\cite{lee2017communication, battey2018distributed}). We provide a slightly more detailed intuitive explanation in Section \ref{int_exp}. 

\subsection{Out-of-sample prediction}

So far, we have discussed the estimation problem. In real applications, out-of-sample prediction is also of interest. We consider a test datapoint $(x_t, y_t)$, generated from the same model $y_t=x_t^\top\beta+\ep_t$, where $x_t,\ep_t$ are independent of $X,\ep$. We want to use $x_t^\top\hbeta$ to predict $y_t$, and the out-of-sample prediction error is defined as $\E(y_t-x_t^\top\hbeta)^2$. Then we have the following proposition.

\begin{proposition}[Out-of-sample prediction error (test error) and relative efficiency]
\label{oe_ridge}
Under the conditions of Theorem \ref{AREsimple}, the limiting out-of-sample prediction error of the optimal distributed estimator $\hbeta_{dist}$ is 
$$
\mathcal{O}_k=\sigma^2+\mathcal{M}_k.
$$
Thus, the asymptotic out-of-sample relative efficiency, meaning the ratio of prediction errors, is
$$
OE=\frac{\mathcal{O}_1}{\mathcal{O}_k}=\frac{\mathcal{M}_1+\sigma^2}{\mathcal{M}_k+\sigma^2},
$$
and the efficiency for prediction is higher than for estimation
$OE\geq ARE.$ 
Furthermore, when the samples are equally distributed, the relative efficiency has the form
$$
\Psi(k,\gamma,\alpha^2)=\frac{1+\gamma m_\gamma(-\gamma/\alpha^2)}{1+\frac{\alpha^2\gamma m_{k\gamma}(-k\gamma/\alpha^2)}{\alpha^2+(1-k)\gamma m_{k\gamma}(-k\gamma/\alpha^2)}},
$$
and the corresponding infinite-worker limit (taking $k\to\infty$) is
$$
\mathcal{H}(\alpha^2,\gamma)=\frac{1+\gamma m_\gamma(-\gamma/\alpha^2)}{1+\frac{\gamma\alpha^2(1+\alpha^2)}{\alpha^2+\gamma(1+\alpha^2)}}.
$$

\end{proposition}

\begin{figure}
\begin{subfigure}{.5\textwidth}
  \centering
\includegraphics[scale=0.4]
{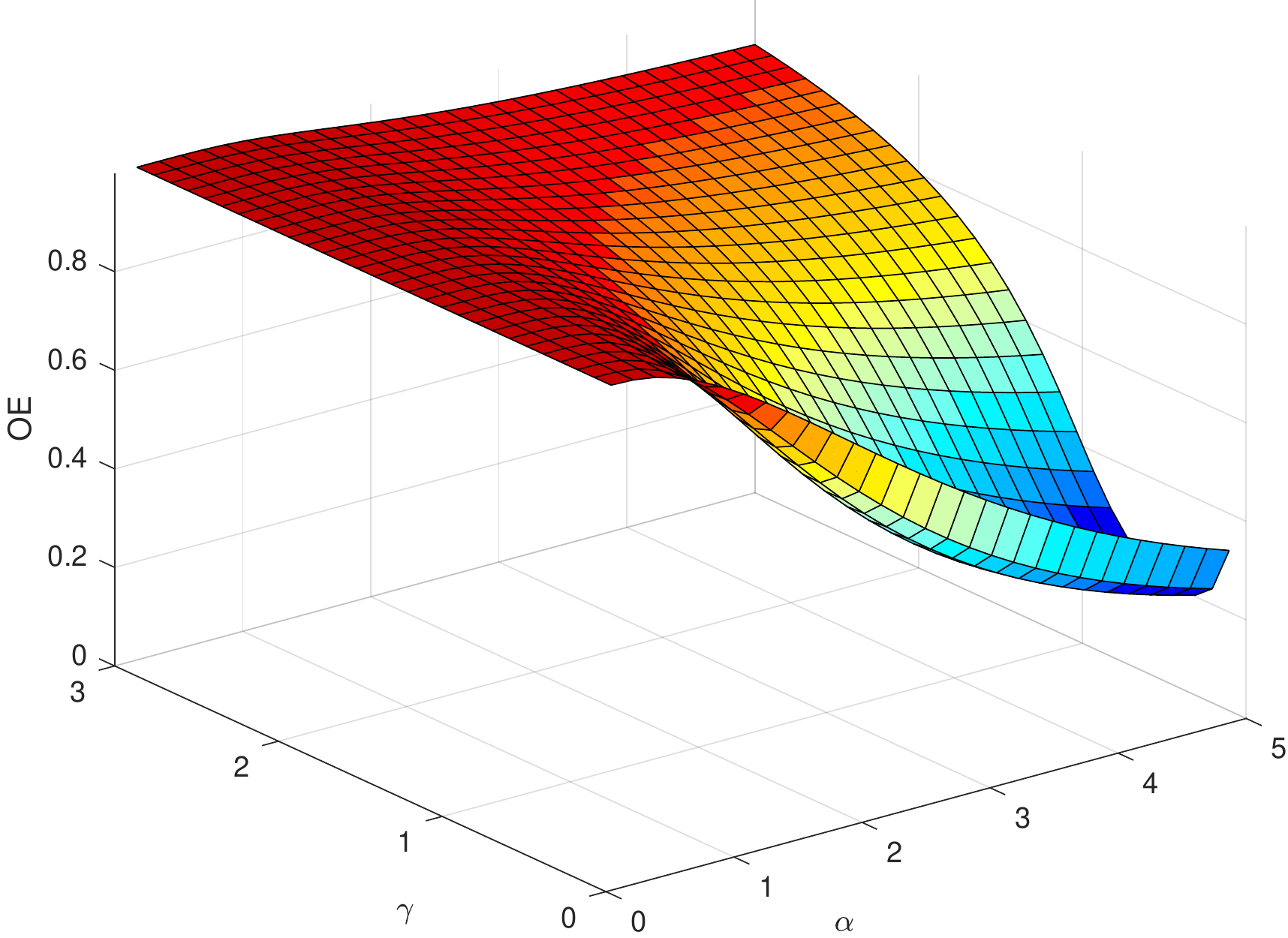}
\end{subfigure}
\begin{subfigure}{.5\textwidth}
  \centering
\includegraphics[scale=0.4]
{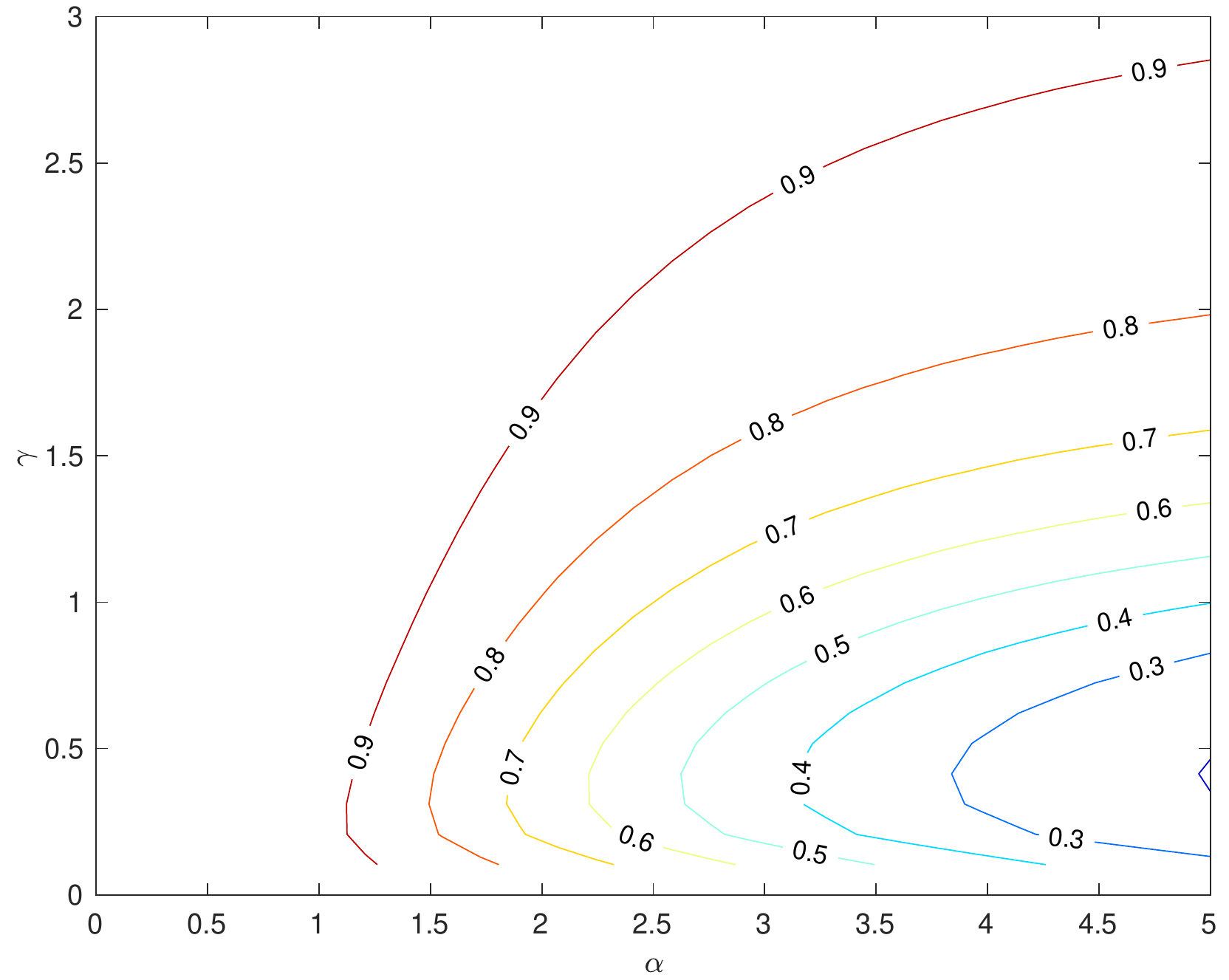}
\end{subfigure}
\caption{Limit of OE: (a) surface and (b) contour plots of $\mathcal{H}(\alpha^2,\gamma)$.}
\label{oefig}
\end{figure}

See Section \ref{pf:oe} for the proof and Figure \ref{oefig} for some plots. This proposition implies that, for the identity covariance case, the efficiency loss of the distributed estimator in terms of the test error is always less than the loss in terms of the estimation error. When the signal-to-noise ratio $\alpha^2$ is small, the relative efficiency is always very large and close to $1$. This observation can be an encouragement to use our distributed methods for out-of-sample prediction.

\subsection{Choosing the regularization parameter}

Previous work found that, under certain conditions, the regularization parameters on the individual machines should be chosen as if they had the all samples \citep{zhang2015divide}. Our findings are consistent with these results.  However, the reasons behind our findings are very different from prior work. The intuition for the previous results is that the \emph{variance} of distributed estimators averages out, while the \emph{bias} does not do so. Therefore, the regularization parameters should be chosen such that the local bias is lower than for locally optimal tuning. This means that we should use smaller regularization parameters locally.

In our case, we find that for isotropic covariance, the optimal risk \emph{decouples across machines}. Hence, the regularization parameters on the machines can be chosen optimally for each machine. Moreover, in our asymptotics the locally optimal choice is a \emph{constant multiple} of the globally optimal choice, namely the multiple in front of the identity matrix in the local ridge estimator $(X_i^\top X_i+n_i\lambda_iI_p)^{-1}X_i^\top Y_i$ should be $\lambda_i = p/(n_i\alpha^2)$ whereas the globally optimal $\lambda$ is $\lambda = p/(n\alpha^2)$. 

Roughly speaking, this derivation reaches the same conclusion as prior work about the choice of regularization parameters, namely that the regularization parameters on the machines should be chosen as if they had the all samples. However, we emphasize that our results are very different, because the optimal weighting procedure has weights summing to \emph{greater than unity}. Moreover, we also consider the proportional-limit case, and the conclusion for regularization parameters only applies to the isotropic case.

\subsection{Implications and practical relevance}
\label{rc}
We discuss some of the implications of our results. Our finite-sample results show that the optimal way to weight the estimators depends on functionals of the unknown parameter $\beta$, while the asymptotic results in general depend on the eigenvalues of $\hSigma$ (or $\Sigma$). These are unavailable in practice, and hence these results can typically not be used on real datasets. However, since our results are precise and accurate (they capture the \emph{truth} about the problem), we interpret this as saying that \emph{the problem is hard in general}. Meaning that optimal weighting for ridge regression is a challenging statistical problem. In practice that means that we may be content with uniform weighting. It remains to be investigated in future work how much we should up-adjust those equal weights.

The optimal weights become usable in the case of spherical data, when $\Sigma=I$ (or, more accurately, the limiting spectral distribution of $\Sigma$ is the point mass at unity). In practice, we can get closer to this assumption by using some form of \emph{whitening} on the data, for instance by scaling all variables to the same scale, by estimating $\Sigma$ over restricted classes, such as assuming block-covariance structures. Alternatively, we can use correlation screening, where we remove features with high correlation. At this stage, all these approaches are heuristic, but we include them to explain how our results can be relevant in practice. It is a topic of future research to make these ideas more concrete. In the algorithm we proposed in Section \ref{algo}, we use grid search to find a good tuning parameter under general covariance structures.

On the theoretical side, our results can also be interpreted as a form of \emph{reduction} between statistical problems. \emph{If} we can estimate the quadratic functionals of the unknown regression parameter involved in our weights, \emph{then} we can do optimally weighted ridge regression. In this sense, we reduce distributed ridge regression to the estimation of those quadratic functionals. We think that in the challenging and novel setting of distributed learning, such reductions can be both interesting and potentially useful. 

An important question is ``Should we use distributed linear or ridge regression?". If we have $n_i \ge p$ and linear regression is defined on each local machine, then we can use either distributed linear \citep{dobriban2018Distributed} or ridge regression. Linear regression has the advantage that the optimal weights are easy to find. Therefore, if we cannot reasonably reduce to the case $\Sigma=I$, it seems we should use linear regression.

\subsection{Minimax optimality of the optimal distributed estimator}
\label{mmx}
To deepen our understanding of the distributed problem, we next show that the optimal distributed ridge estimator is asymptotically rate-minimax. Suppose without loss of generality that the noise level $\sigma^2=1$, and let $\mathbb S^{p-1}(\alpha)=\{\beta\in\mathbb R^p; ||\beta||=\alpha\}$ denote the sphere of radius $\alpha\geq 0$ in $\mathbb R^p$ centered at the origin. Then the minimax risk for estimating $\beta$ over the sphere $\mathbb S^{p-1}(\alpha)$ is
$$
r(\alpha)=\inf_{\hbeta}\sup_{\beta\in\mathbb S^{p-1}(\alpha)}R(\hbeta,\beta)=\inf_{\hbeta}\sup_{\beta\in\mathbb S^{p-1}(\alpha)}\E_\beta||\hbeta-\beta||^2,
$$
where the expectation is over both $X$ and $\ep$. This problem has been well studied by \cite{dicker2014ridge}, who reduced it to the following Bayes problem. 
Let $\pi$ be the uniform measure on $\mathbb S^{p-1}(\alpha)$. Then the Bayes risk with respect to $\pi$ is
$$
r_B(\alpha)=\inf_{\hbeta}\int_{\mathbb S^{p-1}(\alpha)}R(\hbeta,\beta)d\pi(\beta)=\inf_{\hbeta}\E_{\pi}||\hbeta-\beta||^2.
$$
The Bayes estimator is the posterior mean
$
\hbeta_{\mathbb S^{p-1}(\alpha)}=\E_{\pi}(\beta|y,X).
$
So the corresponding Bayes risk is 
$
r_B(\alpha)=\E_{\pi}||\hbeta_{\mathbb S^{p-1}(\alpha)}-\beta||^2.
$
Then, the Bayes estimator also minimizes the original minimax risk and $r(\alpha)=r_B(\alpha)$ \citep{dicker2014ridge}. 

Recall that the ridge estimator with optimally tuned regularization parameter is
$$
\hbeta_r(\alpha)=(X^\top X+\frac{p}{\alpha^2}I_p)^{-1}X^\top Y,
$$
which can be interpreted as the posterior mean of $\beta$ under the normal prior assumption $\beta\sim \N(0,\alpha^2/pI_p)$. When $p$ is very large, the normal distribution $\N(0,\alpha^2/pI_p)$ is very close to the uniform distribution on $\mathbb S^{p-1}(\alpha)$, so we would expect that $\hbeta_{\mathbb S^{p-1}(\alpha)}\approx\hbeta_r(\alpha)$. With this intuition, \cite{dicker2014ridge} further showed that, as $p,n\to\infty,p/n\to\gamma\in(0,\infty)$, for any $\beta\in\mathbb S^{p-1}(\alpha)$
$$
\lim_{n,p\to\infty} 
\left[R(\hbeta_{\mathbb S^{p-1}(\alpha)},\beta)-R(\hbeta_r(\alpha),\beta)\right]=0.
$$
So the global ridge estimator is asymptotically exact minimax. 

We call an estimator is \emph{asymptotically rate-minimax} if asymptotically its risk is at most a constant times the minimax risk. For our distributed problem, we have the following result:
\begin{theorem}[Minimaxity of the optimal distributed estimator]
\label{minimax}
For fixed signal strength $\alpha^2$, the optimally weighted distributed ridge estimator is asymptotically rate minimax. Specifically, its risk $\mathcal{M}_k$ is less than the risk $\mathcal{M}_1$ of the global ridge estimator multiplied by a constant $C=1+\alpha^2$ which only depends on the signal strength $\alpha^2$, and not on the aspect ratio $\gamma=\lim p/n$ and number of machines $k$. Specifically
$$
\mathcal{M}_k\leq (1+\alpha^2)\mathcal{M}_1.
$$
Moreover, for fixed aspect ratio $\gamma>1$, the distributed risk $\mathcal{M}_k$ is less than the global risk $\mathcal{M}_1$ times a constant $C'=\gamma^2/(\gamma^2-1)$ which is independent of $\alpha^2$ and $k$, i.e.
$$
\mathcal{M}_k\leq \frac{\gamma^2}{\gamma^2-1}\mathcal{M}_1.
$$
Therefore, in either case, the optimally weighted distributed ridge estimator is asymptotically rate minimax.
\end{theorem}

See Section \ref{pf:minimax} for the proof. The minimax optimality result is nontrivial, and does not hold for some simpler estimators. For instance, for the null estimator $\hbeta_{null}=0$, the corresponding ARE can be written in terms of the optimal risk function $\phi(\gamma)$ as
$$
\lim_{n,p\to\infty}\frac{R(\hbeta_r(\alpha), \beta)}{R(\hbeta_{null}, \beta)}=\frac{\phi(\gamma)}{\alpha^2}=\frac{\gamma m_\gamma(-\gamma/\alpha^2)}{\alpha^2}.
$$
When $\gamma\to\infty$, we know that $\gamma/\alpha^2m_\gamma(-\gamma/\alpha^2)\to1$, so that even \emph{the null estimator is asymptotically exact minimax}. In this regime, exact minimaxity is a weak result. When $\gamma\to 0$ however, we have $\gamma/\alpha^2m_\gamma(-\gamma/\alpha^2)\to0$ for any $\alpha$, and so the null estimator does not perform well (has zero efficiency). However, the distributed estimator is still asymptotically rate-minimax. 

\section{WONDER: Algorithms for weighted one-shot distributed ridge regression}
\label{algo}

So far, most of our results on distributed ridge regression are purely theoretical. In practice, it would be very helpful to have an implementable algorithm. In fact, our theory for distributed ridge regression allows us to develop an efficient algorithm which works for designs $X$ with arbitrary covariance structures $\Sigma$.

Recall that we have $n$ samples distributed across $k$ machines. For simplicity, let us assume the samples are equally distributed. On the $i$-th machine, we compute a local ridge estimator $\hbeta_i$, local estimators $\hsigma^2_i$, $\hat{\alpha}^2_i$ of the signal-to-noise ratio and the noise level. From Theorem \ref{ARE_equal}, we know that the other quantities needed to find the optimal weights are $m, m'$ and $\lambda$. For $m$ and $m'$, by the definition of the Stieltjes transform, they can be approximated by
\begin{equation*}
\frac{\tr(\hSigma_i+\lambda I)^{-1}}{p}\approx m(-\lambda)
~~\textnormal{and}~~\frac{\tr(\hSigma_i+\lambda I)^{-2}}{p}\approx m'(-\lambda).
\end{equation*}
Here we only need to use local data. The remaining question is: how do we choose the tuning parameter $\lambda$? One way may be grid search. From the theory for the isotropic design, a proper initial guess would be $\lambda=kp/(n\alpha^2)$. Then we can search around this initial guess to find a good parameter with small prediction error.

We assume the data are already mean-centered, which can be performed exactly in one additional round of communication, or approximately by centering the individual datasets. 

Now we have all the quantities we need for our Weighted ONe-shot DistributEd Ridge regression algorithm (WONDER). We send them to a global machine or data center, and aggregate them to compute a weighted ridge estimator. See Algorithm \ref{alg:distalgo1} for more details. WONDER is communication efficient as the local machines only need to send the local ridge estimator $\hbeta_i$ and some scalars to the global datacenter.

\begin{algorithm2e}[h]
    \SetKwInOut{Input}{Input}
    \SetKwInOut{Output}{Output}
    
    \Input{Data matrices ($n_i \times p$) and outcomes ($n_i \times 1$), $(X_i, Y_i)$ distributed across $k$ sites}
    \Output{Distributed ridge estimator $\hbeta_{dist}$ of regression coefficients $\beta$}
    \BlankLine
    \For{$i \leftarrow 1$ \KwTo $k$}{
    Compute the MLE $\htheta_i=(\hsigma^2_i, \hat{\alpha}^2_i)$ locally on $i$-th machine\;
    Send $\htheta_i$ to the global data center\;
    }
    At the data center, combine $\htheta_i$ to get a global estimator $\htheta=(\hsigma^2, \hat{\alpha}^2)=k^{-1}\sum_{i=1}^k\htheta_i$ and send it back to the local machines\;
    Choose a set of tuning parameters $\mathcal{S}$ around the initial guess $\lambda_0=kp/(n\hat{\alpha}^2)$\;
    \For{$\lambda\in\mathcal{S}$}{
    \For{$i \leftarrow 1$ \KwTo $k$}{
    Compute the local ridge estimator $\hbeta_i(\lambda)=(X_i^\top X_i+n_i\lambda I_p)^{-1}X_i^\top Y_i$\;
    Compute the weight $\omega_i$ for the $i$-th local estimator by using the formulas from Theorem \ref{ARE_equal}:
    $$
    \omega_i(\lambda)=\frac{\hsigma^2\hat{\alpha}^2(1-\lambda m)}{\mathcal{F}+k\mathcal{G}}
    $$
    where we use $\tr(X_i^\top X_i/n_i+\lambda I)^{-1}/p$ to approximate $m$, and use $\tr(X_i^\top X_i/n_i+\lambda I)^{-2}/p$ to approximate $m'$\;
    Send $\hbeta_i(\lambda)$ and $\omega_i(\lambda)$ to the global data center;
    }
    Evaluate the performance of the distributed ridge estimator $\hbeta_{dist}(\lambda)=\sum_{i=1}^k\omega_i(\lambda)\hbeta_i(\lambda)$ on validation sets\;
    }
    Select the best tuning parameter $\lambda^*$ and output the corresponding distributed ridge estimator $\hbeta_{dist}(\lambda^*)=\sum_{i=1}^k\omega_i(\lambda^*)\hbeta_i(\lambda^*)$.
    \caption{WONDER: Weighted ONe-shot DistributEd Ridge regression algorithm, general design}
    \label{alg:distalgo1}
\end{algorithm2e}

\begin{algorithm2e}[h]
    \SetKwInOut{Input}{Input}
    \SetKwInOut{Output}{Output}
    
    \Input{Data matrices ($n_i \times p$) and outcomes ($n_i \times 1$), $(X_i, Y_i)$ distributed across $k$ sites}
    \Output{Distributed ridge estimator $\hbeta_{dist}$ of regression coefficients $\beta$}
    \BlankLine
    \For{$i \leftarrow 1$ \KwTo $k$}{
    Compute the MLE $\htheta_i=(\hsigma^2_i, \hat{\alpha}^2_i)$ locally on $i$-th machine\;
    Set local aspect ratio $\gamma_i=p/n_i$\;
    Set regularization parameter $\lambda_i=\gamma_i/\hat{\alpha}^2_i$\;
    Compute the local ridge estimator $\hbeta_i(\lambda_i)=(X_i^\top X_i+n_i\lambda_iI_p)^{-1}X_i^\top Y_i$\;
    Send $\htheta_i, \gamma_i$ and $\hbeta_i$ to the global data center.
    }
    At the data center, combine $\htheta_i$ to get a global estimator $\htheta=(\hsigma^2, \hat{\alpha}^2)$, by $\htheta=k^{-1}\sum_{i=1}^k\htheta_i$\;
    Evaluate the optimal risk functions for $i=1,2,\dots,k$
    $$
    \phi(\gamma_i)=\gamma_i m_{\gamma_i}(-\gamma_i/\hat{\alpha}^2)=\frac{-\gamma_i/\hat{\alpha}^2+\gamma_i-1+\sqrt{(-\gamma_i/\hat{\alpha}^2+\gamma_i-1)^2+4\gamma_i^2/\hat{\alpha}^2}}{2\gamma_i/\hat{\alpha}^2};
    $$\\
    Compute the optimal weights $\omega$, where the $i$-th coordinate of $\omega$ is 
    $$
    \omega_i=\left(\frac{\hat{\alpha}^2}{\phi(\gamma_i)}\right)\cdot \left(\frac{1}{1+\sum_{i=1}^k\left[\frac{\hat{\alpha}^2}{\phi(\gamma_i)}-1\right]}\right);
    $$\\
    Output the distributed ridge estimator $\hbeta_{dist}=\sum_{i=1}^k\omega_i\hbeta_i$.
    \caption{WONDER: Weighted ONe-shot DistributEd Ridge regression algorithm, isotropic design}
    \label{alg:distalgo2}
\end{algorithm2e}

For identity covariance, our results lead to a much simpler WONDER algorithm which requires even less communication and computation. See Algorithm \ref{alg:distalgo2}.

In the above WONDER algorithms, we combine the local estimators of the noise level and signal strength $\htheta_i$ to find a global estimator $\htheta$. A simple method is to take the average: $\htheta=k^{-1}\sum_{i=1}^k\htheta_i$. Another option is to use inverse-variance weighting, based on the asymptotic variance of the MLE (which then of course has to be estimated). 

Based on the results so far, it follows that our WONDER algorithm can consistently estimate the limiting optimal weights, and moreover it has asymptotically optimal mean squared error among all weighted distributed ridge estimators, at least for the identity covariance case. We omit the details.

\section{Experimental results}
\label{expr}
We present some numerical results in addition to the ones already shown in the paper.


\subsection{Finite-sample comparison of relative efficiency for isotropic covariance}

\begin{figure}
\begin{subfigure}{.5\textwidth}
  \centering
\includegraphics[scale=0.14]
{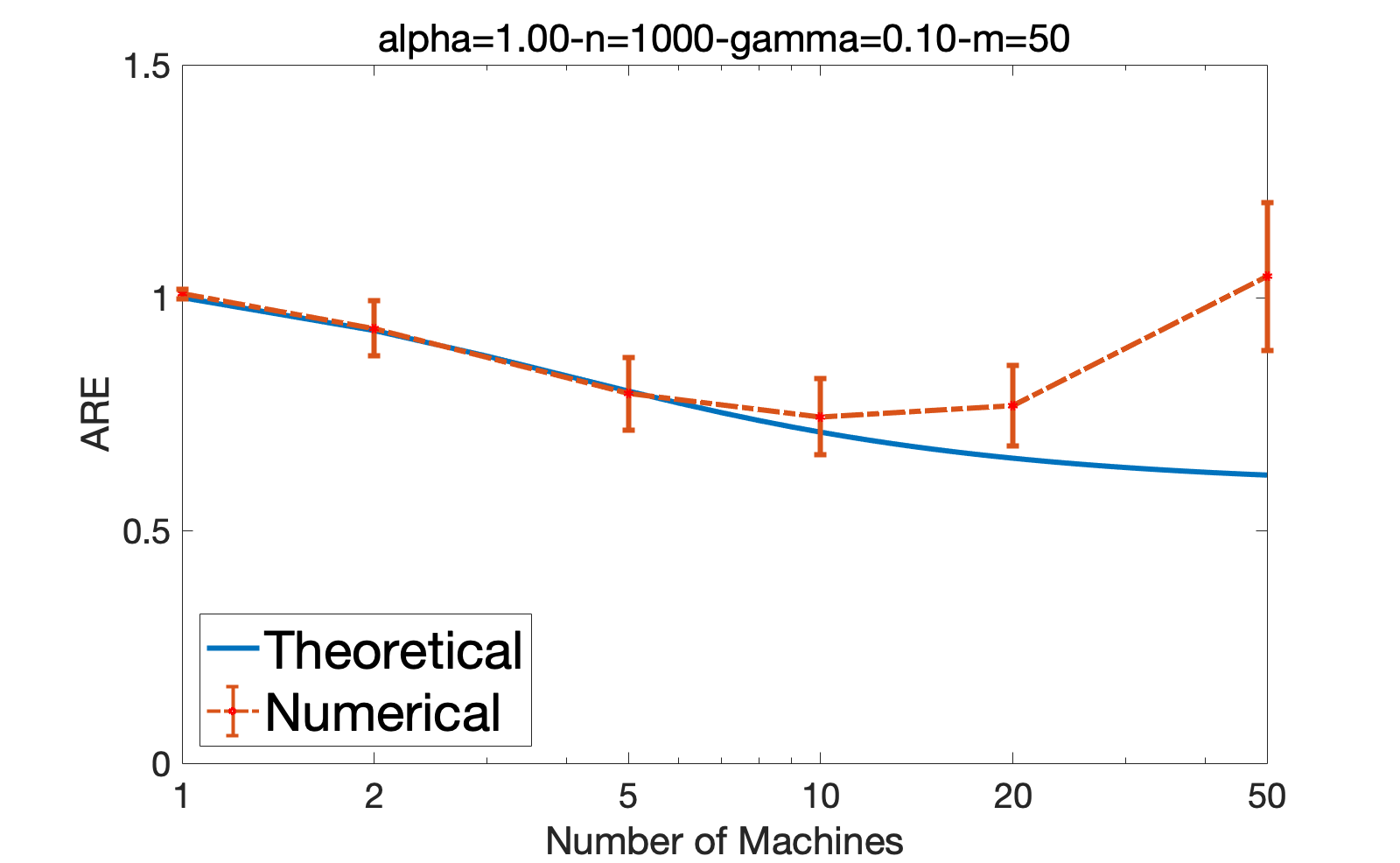}
\end{subfigure}
\begin{subfigure}{.5\textwidth}
  \centering
\includegraphics[scale=0.14]
{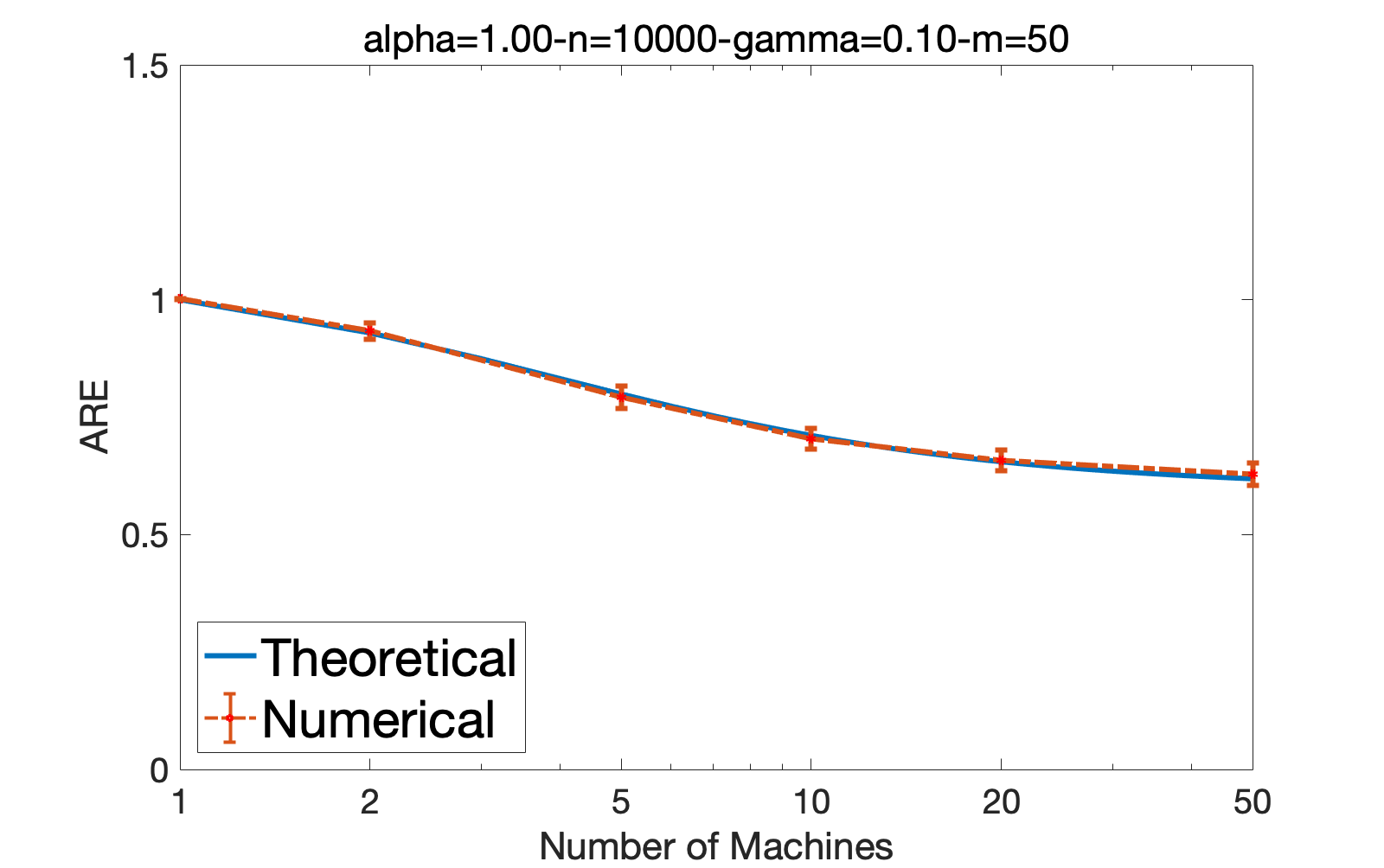}
\end{subfigure}
\caption{Realized relative efficiency in a regression simulation.}
\label{ActualARE}
\end{figure}

Figure \ref{ActualARE} shows a comparison of the theoretical formulas for ARE and realized relative efficiency in a regression simulation. Here the regression model is $Y=X\beta+\ep$, where $X$ is $n\times p$ with i.i.d. standard normal entries, $\beta$ is a $p$-dimensional random vector with i.i.d. mean $0$, variance $\alpha^2/p$ normal entries, and $\ep$ also has i.i.d standard normal entries. For each $k=1,2,5,10,20,50$, we split the data equally into $k$ groups and perform ridge regression on each group. For each group, we choose the same tuning parameter $\lambda_i=p/(n_i\alpha^2)$. For the global regression on the entire dataset, we choose the tuning parameter $\lambda=p/(n\alpha^2)$ optimally. 

We show the results of the expression for the realized relative efficiency $\|\hbeta-\beta\|^2/\|\hbeta_{dist}-\beta\|^2$ compared to the theoretical ARE. We generate 100 independent copies of $\ep$, perform regression, recording the realized relative efficiency $||\hbeta-\beta||^2/||\hbeta_{dist}-\beta||^2$, as well as its overall Monte Carlo mean. For the first plot, we take $n=1000,$ $p=100$, and $\alpha=\sigma=1$. 

As we can see in the plot, the theoretical formula is accurate only for a small number of machines. It turns out that this is due to finite-sample effects. In the second plot, we set $n=10000,$ $p=1000$ and $\alpha=\sigma=1$ such that the aspect ratio $\gamma=p/n$ is the same as before. In that case the theoretical formula becomes very accurate.

\subsection{Choosing the regularization for general covariance}

How can we choose the optimal regularization parameters when the predictors have a general covariance structure $\Sigma$? In this case, our theoretical results do not give an explicit expression for the optimal regularization parameters. In practice, one can use techniques like cross-validation to do selections. Here we present simulation results to shed light on the important question of how to choose them.

\begin{figure}
\begin{subfigure}{.5\textwidth}
  \centering
\includegraphics[scale=0.5]
{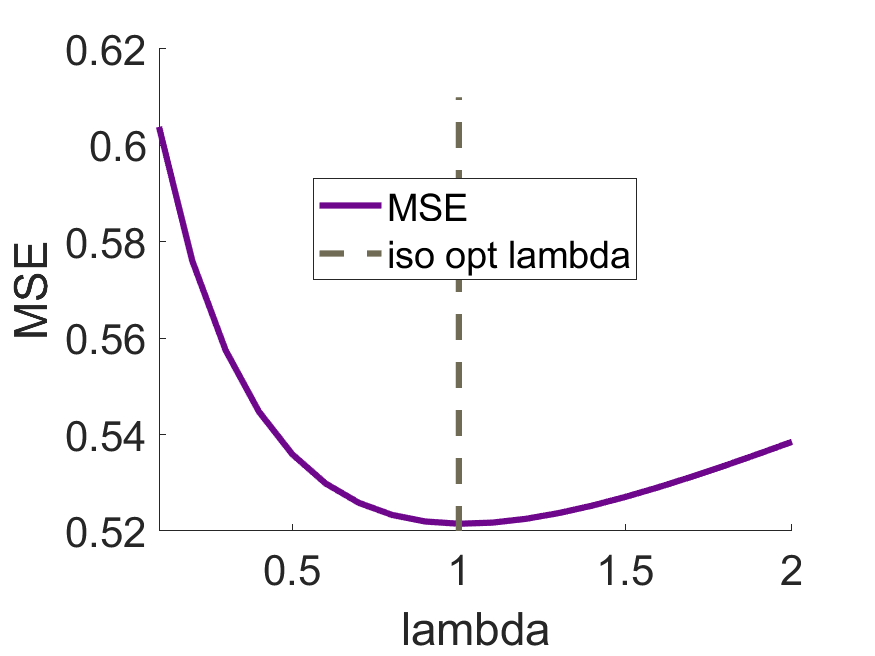}
\end{subfigure}
\begin{subfigure}{.5\textwidth}
  \centering
\includegraphics[scale=0.5]
{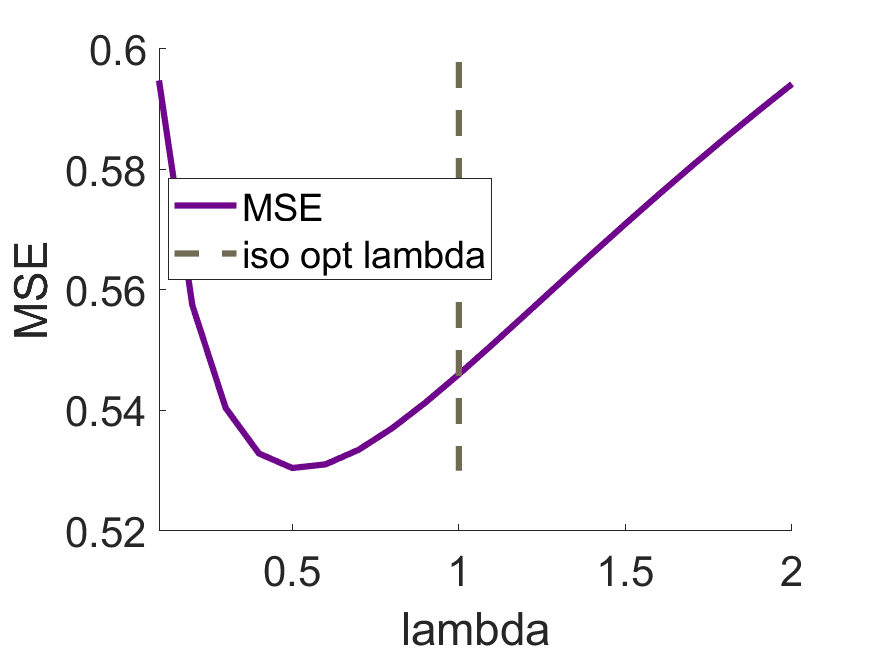}
\end{subfigure}
\begin{subfigure}{.5\textwidth}
  \centering
\includegraphics[scale=0.5]
{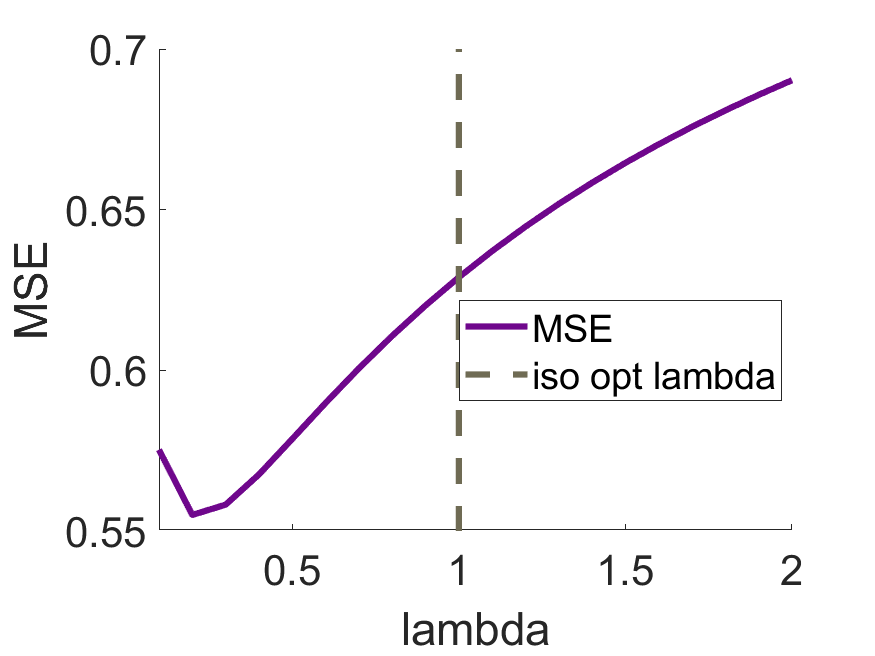}
\end{subfigure}
\begin{subfigure}{.5\textwidth}
  \centering
\includegraphics[scale=0.5]
{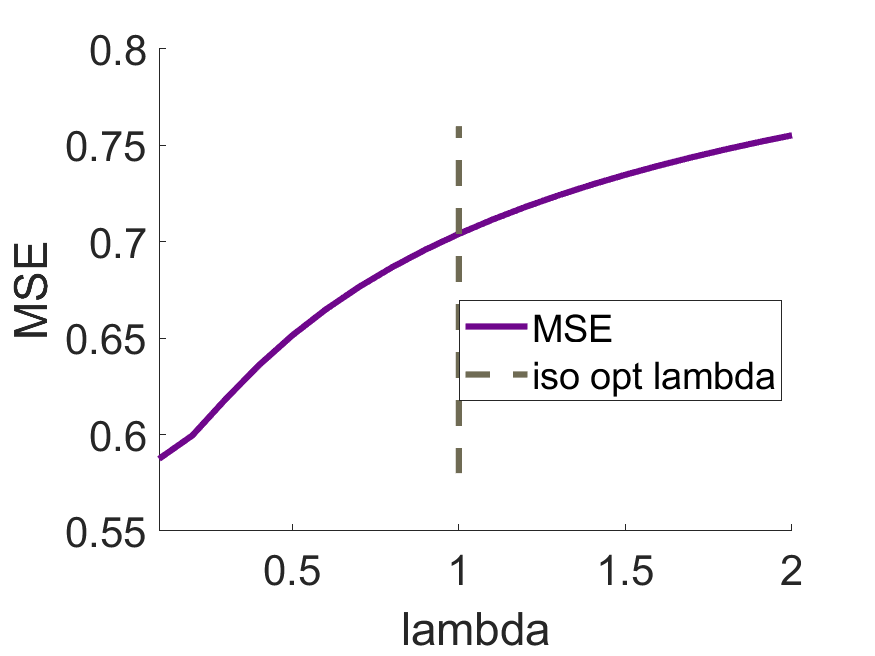}
\end{subfigure}
\caption{Distributed risk as a function of the regularization parameter. We plot the risk of the optimally weighted distributed estimator for an AR-1 covariance structure. We set $\alpha=1$, $\gamma=0.17$ and $k=1,2,5,10$.}
\label{gc}
\end{figure}

We use a similar simulation setup as in the previous sections, except we generate the datapoints independently from an autoregressive model of order one (AR-1), i.e., each datapoint $x_i$ is generated as $x_i \sim \N(0,\Sigma)$, where $\Sigma_{ij} = \rho^{|i-j|}$, and $\rho$ is the autocorrelation parameter. We choose $\rho = 0.9$. We also choose $n=3000$, $p=500$, and report the results of a simulation where we average over $n_{mc} = 20$ independent realizations of $\beta$. Figure \ref{gc} shows the optimal distributed risk $M^*(k)$ as a function of the local regularization parameter $\lambda$. We set all local regularization parameters to equal values, which is reasonable, since the local problems are exchangeable. We also parametrize the regularization parameters as multiples of the optimal parameter for the isotropic case (which equals $k\gamma/\alpha^2$).

We observe that for $k=1$, the optimal parameter is the same as in the isotropic case. This makes sense, because the optimal regularization parameter for one machine is always the same, regardless of the structure of the design. However for $k>1$, we observe that the regularization parameters are \emph{smaller} than the isotropic ones. This is an insight that has apparently not been available before. It is an interesting topic of future work to develop an intuitive understanding.

\subsection{Experiments on empirical data}


\begin{figure}
\centering
\includegraphics[scale=0.2]
{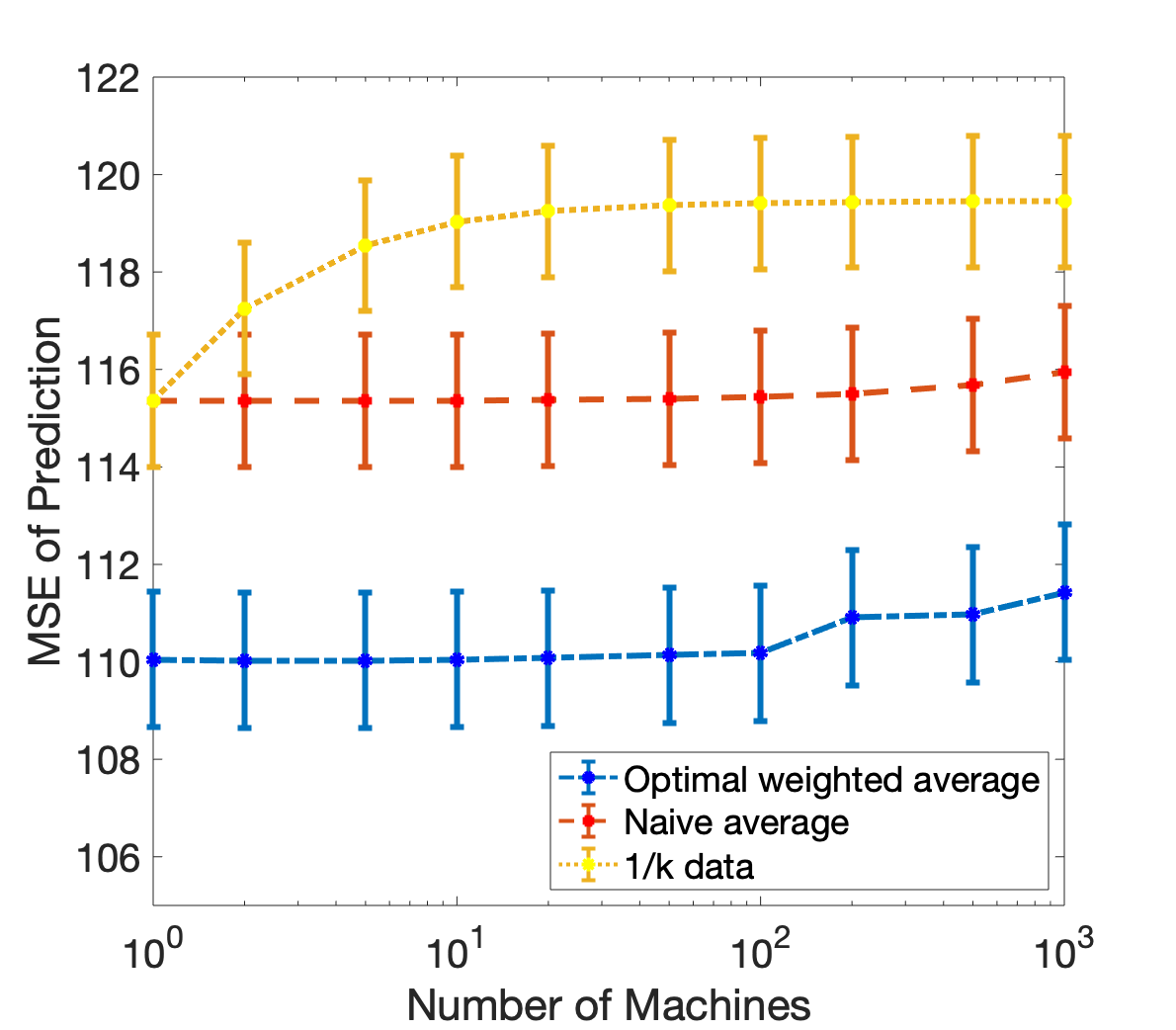}
\caption{Million Song Year Prediction Dataset (MSD). Optimal weighted average (WONDER), Naive average, and regression on $1/k$ fraction of data.}
\label{MSD}
\end{figure}

In this section, we present an empirical data example to examine the accuracy of our theoretical results. It is reasonable to compare the performance of different estimators in terms of the prediction (test) error. Figure \ref{MSD} shows a comparison of three estimators including our optimal weighted estimator on the Million Song Year Prediction Dataset (MSD) \citep{Bertin-Mahieux2011}. 

Specifically, we perform the following steps in our data analysis. We download the dataset from the UC Irvine Machine Learning Repository. The original dataset has $N=515,345$ samples and $p=91$ features. The dataset has already been divided into a training set and a test set. The training set consists of the first $463,715$ samples and the test set contains the rest. We attempt to predict the release year of a song. Before doing distributed regression, we first center and normalize both the design matrix $X$ and the outcome $Y$. Now we are ready to do ridge regression under the distributed setting. 

For each experiment, we randomly choose $n_{train}=10,000$ samples from the training set and $n_{test}=1,000$ samples from the test set. We construct the estimators on the training samples. Then we perform ridge regression in a distributed way to obtain our optimal weighted WONDER estimator as described in Algorithm \ref{alg:distalgo1}. We measure its performance on the test data by computing its MSE for prediction.  We choose the number of machines to be $k=1,2,5,10,20, 50,100,200,500,1000$, and we distribute the data evenly across the $k$ machines. Here we try different tuning parameters $\lambda$ around $kp/(n_{train}\cdot\hat{\alpha}^2)$, and use $\lambda=3kp/(n_{train}\cdot\hat{\alpha}^2)$ as our final parameter. (In practice, one may try a 1-D grid search to find the right scale.)

For comparison, we also consider two other estimators:
\benum

\item
The distributed estimator where we take the naive average (weight for each local estimator is simply $1/k$) and choose the local tuning parameter $\lambda=p/(n_{train}\cdot\hat{\alpha}^2)$. This formally agrees with the divide-and-conquer type estimator proposed in \cite{zhang2015divide}.
\item
The estimator using only a fraction $1/k$ of the data, which is just one of the local estimators. For this estimator, we choose the tuning parameter $\lambda=kp/(n_{train}\cdot\hat{\alpha}^2)$.
\eenum

We repeat the experiment for $T=100$ times, and report the average and $1/4$ standard deviation over all experiments on Figure \ref{MSD}. Each time we randomly collect new training and test sets. 

From Figure \ref{MSD}, we observe the following:
\benum

\item The WONDER estimator has smaller MSE than both the local estimator and the naive averaged estimator, which means optimal weighting can indeed help. 

\item It seems that data splitting does not have huge impact on the performance of the WONDER estimator. This phenomenon is compatible with our theory. Since the signal-to-noise ratio $\alpha^2$ is about $1.2$ for this data set, we are in a low SNR scenario. From Proposition \ref{oe_ridge} and Figure \ref{oefig}, we see that the performance of the distributed estimator is close to the global estimator in terms of the prediction error. 

\eenum

To conclude, in terms of computation-statistics tradeoff, this example suggests a very positive outlook on using distributed ridge regression via WONDER: The accuracy is affected very little even though the data is split up into 100 parts. Thus we save at least 100x in computation time, while we have nearly no loss in performance.

Finally, we mention that in Figure 4 of \cite{zhang2015divide}, the authors also compare the performance of the distributed estimator to the local estimator on the same Million Song data set. We notice that the MSE of prediction in their experiments is usually between $80$ and $90$, and variance is typically very small. In our experiments, both the MSE and variance are larger. The reason for this seems to be that they consider more general kernel ridge regression.

\section*{Acknowledgements}

The authors thank Yuekai Sun for discussions motivating our study, as well as John Duchi, Jason D. Lee, Xinran Li, Jonathan Rosenblatt, Feng Ruan, and Linjun Zhang for helpful discussions. They are grateful to Sifan Liu for thorough comments on an earlier version of the manuscript. They are also grateful to the associate editor and referees for valuable suggestions. ED was partially supported by NSF BIGDATA grant IIS 1837992. 
 
\section{Appendix}

\subsection{Proof of Theorem \ref{re_ridge}}
\label{pf:re_ridge}

We can calculate the MSE of the weighted sum as
\begin{align*}
M(w)= \E\left\|\sum w_i \hat\beta_i-\beta\right\|^2 &
= \E(\sum w_i \hat\beta_i-\beta)^\top (\sum w_j \hat\beta_j-\beta) \\
&= \sum_{ij} w_iw_j \cdot \E\hat\beta_i^\top\hat\beta_j
- 2\sum_i w_i \E\hat\beta_i^\top\beta
+ \|\beta\|^2.
\end{align*}
Let $\hB$ be the $p\times k$ matrix defined as $\hB = [\hbeta_1,\ldots,\hbeta_k]$. Then we can write the above MSE as
\begin{align*}
M(w)= w^\top \E \hB^\top \hB w - 2 \E \beta^\top\hB w
+ \|\beta\|^2.
\end{align*}
Let also 
\beqs B  = \E \hB = [\E\hbeta_1,\ldots,\E\hbeta_k].\eeqs
Since the local estimators are independent, we can write
\beqs
M(w)= w^\top (B^\top B+R) w - 2 \beta^\top B w + \|\beta\|^2,
\eeqs
where $R$ is a diagonal matrix with entries 
\beqs 
R_i = \E \|\hbeta_i\|^2 - \|\E\hbeta_i\|^2= \E \|\hbeta_i - \E\hbeta_i\|^2. 
\eeqs
The objective function $M(w)$ can be viewed as corresponding to a $k$-parameter linear regression problem, with unknown parameters $w_i$, design matrix $B$ and outcome vector $\beta$. Specifically, we regress $\beta$ on $\E\hB=  \E[\hbeta_1,\ldots,\hbeta_k]$. Therefore, the optimal weights are
\beqs w^* = (B^\top B+R)^{-1} B^\top\beta,
\eeqs 
and the optimal risk equals
\beqs 
M^* = M(w^*) = \beta^\top\left[I_p -B (B^\top B+R)^{-1} B^\top \right]\beta.
\eeqs
Now, to find $B=\E \hB$, we need $\E\hat\beta_i$. The expectation of the ridge regression estimator for the full dataset is
\beqs \E \hbeta(\lambda) = \E(X^\top X+n\lambda I_p)^{-1} X^\top Y = (X^\top X+n\lambda I_p)^{-1} X^\top X\beta. \eeqs
Letting $\hSigma = n^{-1} X^\top X$, this equals $\E \hbeta(\lambda) = (\hSigma+\lambda I_p)^{-1} \hSigma\beta$. Similarly, 
\beqs \E \hbeta_i(\lambda_i) 
= (X_i^\top X_i+n_i\lambda_i I_p)^{-1} X_i^\top X_i\beta.\eeqs
Let $Q_i = Q_i(\lambda_i) = (X_i^\top X_i+n_i\lambda_i I_p)^{-1} X_i^\top X_i$ be those matrices and let $\hSigma_i = n^{-1} X_i^\top X_i$. Then the above equals $Q_i = (\hSigma_i+\lambda_i I_p)^{-1} \hSigma_i$, and 
\beqs B  = [Q_1\beta; \ldots; Q_k\beta].\eeqs
Therefore, $B^\top B$ has entries $\beta^\top Q_iQ_j\beta$, while $B^\top \beta$ has entries $\beta^\top Q_i \beta$. Moreover, 
$$R_i 
= \E \|\hbeta_i - \E\hbeta_i\|^2 
= \E \| (X_i^\top X_i+n_i\lambda_iI_p)^{-1}X_i^\top \ep_i\|^2
=  \sigma^2\tr[(X_i^\top X_i+n_i\lambda_iI_p)^{-2}X_i^\top X_i].$$
We can also write this as $R_i = n_i^{-1}\sigma^2\tr[(\hSigma_i+\lambda_i I_p)^{-2} \hSigma_i]$. To conclude the optimal risk, we have
\beqs M ^*(k) = \|\beta\|^2- v^\top (A+R)^{-1} v, \eeqs
where 
\begin{align*}
v  &= B^\top \beta = vec[\beta^\top Q_i \beta],\\
A& = mat[\beta^\top Q_iQ_j\beta],\\
R& = \diag\left[n_i^{-1}\sigma^2\tr[(\hSigma_i+\lambda_i I_p)^{-2} \hSigma_i]\right],\\
Q_i& = (\hSigma_i+\lambda_i I_p)^{-1} \hSigma_i.
\end{align*}
Here we used the vectorization and to-matrix operators $vec,mat$. For the global MSE, we only need to consider the special case where $k=1$, which gives us 
$$\E||\hbeta-\beta||^2=M ^*(1) = \|\beta\|^2-\frac{(\beta^\top Q\beta)^2}{\beta^\top Q^2\beta+\sigma^2\tr[(X^\top X+n\lambda I_p)^{-2}X^\top X]},$$ 
where $Q=(\hSigma+\lambda I_p)^{-1}\hSigma$. This finishes the argument.

\subsection{Adding a constant to the regression}
\label{add_const}
We show below the details of the derivation of optimal weights for ridge regression when we also add a constant to the (biased) local estimators. In our calculation from Theorem \ref{re_ridge}, we need to change some details as follows: 

We need to define a new matrix  $\hB = [\hbeta_1,\ldots,\hbeta_k, p^{-1/2} 1_p]$ and new weights $w = [w; w_{k+1}]$. Clearly, we still have that 

$$B  = [\E\hbeta_1,\ldots,\E\hbeta_k, p^{-1/2} 1_p] = [Q_1\beta; \ldots; Q_k\beta, p^{-1/2} 1_p].$$

The new matrix $R$ is now diagonal with all entries as before, and the lower right corner entry is $R_{k+1} = 0$.

We consider the same regression problem as before, except we add an intercept into the matrix $B$ as above. The same algebraic form of the optimal weights and risk holds, with the new definitions above. The optimal risk is now 

\beqs M ^*(k) = \|\beta\|^2- v^\top (A+R)^{-1} v \eeqs
where 
\begin{align*}
v  &= B^\top \beta = [vec[\beta^\top Q_i \beta]; p^{-1/2} 1_p^\top\beta]\\
A& = 
\left[
  \begin{array}{ c c }
     mx[\beta^\top Q_iQ_j\beta] & vec[p^{-1/2} 1_p^\top Q_i \beta] \\
     vec[p^{-1/2} 1_p^\top Q_i \beta]  & 1 
  \end{array} \right]\\
R& = \diag\left[n_i^{-1}\tr[(\hSigma_i+\lambda_i I_p)^{-2} \hSigma_i]; 0\right]\\
Q_i& = (\hSigma_i+\lambda_i I_p)^{-1} \hSigma_i
\end{align*}
In simulation studies, we have observed that this approach typically does not lead to a significant decrease in MSE.



\subsection{Differentiation rule for calculus of deterministic equivalents}
\label{pf:diffcalculus}

\begin{theorem}[Differentiation rule for the calculus of deterministic equivalents]
\label{diffcalculus}
Suppose $T=T_n$ and $S=S_n$ are two (deterministic or random) matrix sequences of growing dimensions such that $f(z,T_n)\asymp g(z,S_n)$, where the entries of $f$ and $g$ are analytic functions in $z\in D$ and $D$ is an open connected subset of $\mathbb C$. Suppose that for any sequence $C_n$ of deterministic matrices with bounded trace norm we have
$$
|\tr\left[C_n(f(z,T_n)-g(z,S_n))\right]|\leq M
$$
for every $n$ and $z\in D$. Then we have $f'(z,T_n)\asymp g'(z,S_n)$ for $z\in D$, where the derivatives are entry-wise with respect to $z$.
\end{theorem}

To prove this theorem, we need to introduce a lemma from complex analysis which is a consequence of the dominated convergence theorem and Cauchy's integral formula. 

\begin{lemma}[see Lemma 2.14 in \cite{bai2009spectral}]
\label{vitali}
Let $f_1, f_2,\dots$ be analytic on the domain $D$, satisfying $|f_n(z)|\leq M$ for every $n$ and $z\in D$. Suppose that there is an analytic function on $D$ such that $f_n(z)\to f(z)$ for all $z\in D$. Then it also holds that $f'_n(z)\to f'(z)$ for all $z\in D$.
\end{lemma}

The proof of theorem \ref{diffcalculus} is clear. Since $\tr\left[C_n(f(z,T_n)-g(z,S_n))\right]$ is a sequence of analytic functions on $D$ with uniform bound, then from the definition of the deterministic equivalence, we have $\tr\left[C_n(f(z,T_n)-g(z,S_n))\right]\to0$. By lemma \ref{vitali}, the derivative also converges to $0$ for all $z\in D$, which finishes the proof.

\subsection{Proof of Theorem \ref{ARE}}
\label{pf:ARE}

The first step is to use the well-known concentration of quadratic forms to reduce to trace functionals (See e.g. Lemma C.3 of \cite{dobriban2018high} which is based on Lemma B.26 of \cite{bai2009spectral}). Since $\beta$ is independent of the data $X$ with mean zero and finite variance, under the moment assumptions imposed in the theorem, we have 
\beqs \beta^\top Q_i \beta - \sigma^2\alpha^2/p \cdot \tr Q_i \to_{a.s.} 0,\eeqs
\beqs \beta^\top Q_iQ_j\beta - \sigma^2\alpha^2/p \cdot \tr Q_i Q_j\to_{a.s.} 0,\eeqs
\beqs \beta^\top Q_i^2 \beta - \sigma^2\alpha^2/p \cdot \tr Q_i^2 \to_{a.s.} 0.\eeqs

Let us compute the limits of $v, A$ and $R$ respectively.
\benum
\item Limit of $v$: First of all, we have already known that 
\beqs \beta^\top Q_i \beta - \sigma^2\alpha^2/p \cdot \tr Q_i \to_{a.s.} 0,\eeqs
so it is sufficient to consider the limit of $\tr Q_i/p$.
Since

\beqs \tr Q_i/p = 1 - \lambda_i \tr[(\hSigma_i+\lambda_i I_p)^{-1}]/p.\eeqs

assuming that the spectral distribution of $\hSigma_i$ converges almost surely to $F_{\gamma_i}$, we thus have

\beqs \tr Q_i/p \to_{a.s.}  1 - \lambda_i \E_{F_{\gamma_i}}(T+\lambda_i)^{-1} = 1- \lambda_i m_{F_{\gamma_i}}(-\lambda_i).\eeqs

Above we have introduced the Stieltjes transform $m_{F_{\gamma_i}}$ as a limiting object. So, 

\beqs \beta^\top Q_i \beta \to_{a.s.} \sigma^2\alpha^2[1- \lambda_i m_{F_{\gamma_i}}(-\lambda_i)].\eeqs

For the form in terms of the population spectral distribution $H$, if $p/n\to\gamma$ and the spectral distribution of $\Sigma$ converges to $H$, we have by the general Marchenko-Pastur (MP) theorem of Rubio and Mestre \citep{rubio2011spectral}, that 
\begin{align*}
(\hSigma+ \lambda I)^{-1} &\asymp (x_p \Sigma + \lambda I)^{-1},
\end{align*}

where  $x_p$ is the unique positive solution of the fixed point equation 
\beqs 1-x_p=\frac{x_p}{n}\tr\left[\Sigma(x_p\Sigma+\lambda I)^{-1}\right]. \eeqs
When $n,p\to\infty$, $x_p\to x$ and $x$ satisfies the equation
\beqs 1 -x = \gamma \left[1 - \lambda\int_0^\infty \frac{dH(t)}{xt+\lambda}\right].\eeqs 
We remark that the assumptions made in the theorem suffice for using the Rubio-Mestre result. Moreover, we only use a special case of their result, similar to \cite{dobriban2018Distributed}.
Hence from the calculus of deterministic equivalents \citep{dobriban2018Distributed}, we can take the traces of the matrices in question to obtain
\beqs \tr Q_i/p 
= 1 - \lambda_i \tr[(\hSigma_i+\lambda_i I_p)^{-1}]/p
\asymp 1 - \lambda_i \tr[(x_i \Sigma + \lambda_i I)^{-1}]/p
\to_{a.s.} \E_H \frac{x_iT}{x_iT+\lambda_i}, \eeqs
where  $x_i=x(H,\gamma_i,-\lambda_i)$ is the unique solution of 
\beqs 1 -x_i = \gamma_i \left[1 - \lambda_i\int_0^\infty \frac{dH(t)}{x_it+\lambda_i}\right].\eeqs

\item Limit of $A$: Let us consider the cases $i\neq j$ and $i=j$ separately.
\benum
\item $i\neq j$: We begin by 
\beqs \beta^\top Q_iQ_j\beta - \sigma^2\alpha^2/p \cdot \tr Q_i Q_j\to_{a.s.} 0.\eeqs
Based on the above expression for $Q_i$, we have 
\beqs Q_iQ_j = I_p - \lambda_i (\hSigma_i+\lambda_i I_p)^{-1}
- \lambda_j (\hSigma_j+\lambda_j I_p)^{-1} 
+ \lambda_i\lambda_j  (\hSigma_i+\lambda_i I_p)^{-1}(\hSigma_j+\lambda_j I_p)^{-1}.\eeqs
So the key will be to find the limit of 
\beqs E_{ij}= p^{-1}\tr \{(\hSigma_i+\lambda_i I_p)^{-1}(\hSigma_j+\lambda_j I_p)^{-1}\}.\eeqs
 From the general MP theorem, since $p/n_i\to\gamma_i$, we have for all $i$,
\begin{align*}
(\hSigma_i+\lambda_i I_p)^{-1} &\asymp (x_{ip}\Sigma+\lambda_i I_p)^{-1}.
\end{align*}

Here $x_{ip}$ is the unique positive solution of the fixed point equation 
\beqs 1 -x_{ip} = \frac{x_{ip}}{n_i} \tr\left[\Sigma(x_{ip}\Sigma+\lambda_iI)^{-1}\right],\eeqs 
and $x_{ip}\to x_i$ as $n_i,p\to\infty$. By the product rule of the calculus of deterministic equivalents, we have for $i\neq j$
\begin{align*}
(\hSigma_i+\lambda_i I_p)^{-1}(\hSigma_j+\lambda_j I_p)^{-1} &
\asymp (x_{ip}\Sigma+\lambda_i I_p)^{-1}(x_{jp}\Sigma+\lambda_j I_p)^{-1}.
\end{align*}

Hence by the trace rule of deterministic equivalents,
\beqs E_{ij}\asymp 
p^{-1}\tr[(x_{ip}\Sigma+\lambda_i I_p)^{-1}(x_{jp}\Sigma+\lambda_j I_p)^{-1}]\eeqs
Moreover, since the spectral distribution of $\Sigma$ converges to $H$, we find for $i\neq j$
\beqs E_{ij} \to \E_H \frac{1}{(x_iT+\lambda_i)(x_jT+\lambda_j)}.\eeqs

Putting it together, 
\beqs Q_iQ_j \asymp I_p - \lambda_i (x_{ip}\Sigma+\lambda_i I_p)^{-1}
- \lambda_j (x_{jp}\Sigma+\lambda_j I_p)^{-1}
+ \lambda_i\lambda_j (x_{ip}\Sigma+\lambda_i I_p)^{-1}(x_{jp}\Sigma+\lambda_j I_p)^{-1}.\eeqs
So, again by the trace rule of deterministic equivalents, we have
\begin{align*}
p^{-1}\tr\{Q_iQ_j\} &
\to_{a.s.} 1-  \E_H \frac{\lambda_i }{x_iT+\lambda_i}
-  \E_H \frac{\lambda_j }{x_jT+\lambda_j}
+  \E_H \frac{\lambda_i\lambda_j }{(x_iT+\lambda_i)(x_jT+\lambda_j)}\\
&= x_ix_j\E_H \frac{T^2}{(x_iT+\lambda_i)(x_jT+\lambda_j)}.
\end{align*}
Therefore, for $i\neq j$
\beqs
A_{ij} \to \sigma^2\alpha^2\left[x_ix_j\E_H \frac{T^2}{(x_iT+\lambda_i)(x_jT+\lambda_j)}\right].
\eeqs

\item $i=j$: In this case, 
\beqs
\beta^\top Q_i^2\beta-\sigma^2\alpha^2/p\cdot\tr Q_i^2\to 0, 
\eeqs
where $Q_i^2=I_p-2\lambda_i(\hSigma_i+\lambda_iI_p)^{-1}+\lambda_i^2(\hSigma_i+\lambda_iI_p)^{-2}$. We can easily find the limit of $\tr Q_i^2/p$ in terms of empirical quantities, based on our knowledge of the convergence of Stieltjes transforms and its derivatives:
\beqs
\tr Q_i^2/p\to 1-2\lambda_im_{F_{\gamma_i}}(-\lambda_i)+\lambda_i^2m'_{F_{\gamma_i}}(-\lambda_i).
\eeqs
Therefore, for $i= j$
\beqs
A_{ii} \to \sigma^2\alpha^2[1-2\lambda_im_{F_{\gamma_i}}(-\lambda_i)+\lambda_i^2m'_{F_{\gamma_i}}(-\lambda_i)].
\eeqs

We can also express the limit of $A_{ii}$ in terms of the population spectral distribution $H$ by using Theorem \ref{diffcalculus}. For our purpose, let $T=\Sigma$, $S = \hSigma$, while
\begin{align*}
f(z,T) & = (x_p T-z I)^{-1},\\ 
g(z,S) & = (S-zI_p)^{-1}.
\end{align*}
From \cite{rubio2011spectral}, we know that for each $z\in D:=\mathbb C\setminus \mathbb R^+$, $f(z,\Sigma)\asymp g(z,\hSigma)$. $x_p$ is defined as 
$$
x_p=\frac{1}{n}\tr[(I+\frac{p}{n}e_pI)^{-1}]=\frac{1}{1+(p/n)e_p}=\frac{1}{1+\gamma_pe_p},
$$
and $e_p=e_p(z)$ is the Stieltjes transform of a certain positive measure on $\mathbb R^+$, obtained as the unique solution of the equation
$$
e_p=\frac{1}{p}\tr[\Sigma(x_p\Sigma-zI_p)^{-1}].
$$
It is well-known that $x_p(z), e_p(z)$ are both analytic functions on $D$. Then we can check that the conditions of theorem \ref{diffcalculus} hold in this case. First of all, for an invertible matrix $A$, $A^{-1}=(\det A)^{-1}A^*$, where $A^*$ is the adjugate matrix of $A$. Since $x_p$ is analytic, it is easy to verify that $\det(x_p\Sigma-zI_p), \det(\hSigma-zI_p)$ and all entries of $(x_p\Sigma-zI_p)^*, (\hSigma-zI_p)^*$ are analytic functions of $z$. So the entries of $f(z, \Sigma)$ and $g(z,\hSigma)$ are analytic in $D$. 

Next, we want to bound 
$$
\tr[C_n((x_p\Sigma-zI_p)^{-1}-(\hSigma-zI_p)^{-1})]\leq||C_n||_{\tr}\cdot||(x_p\Sigma-zI_p)^{-1}-(\hSigma-zI_p)^{-1}||_2.
$$
For a fixed $\delta>0$, let us define a domain $D_{\delta}:=\{z\in D:\textnormal{Re}z<-\delta\}\cup\{z\in D: |\textnormal{Im}z|>\delta\}$. Then, it is sufficient to find a uniform bound for $||(x_p\Sigma-zI_p)^{-1}-(\hSigma-zI_p)^{-1}||_2$ on $D_\delta$. In fact, we can bound $||(x_p\Sigma-zI_p)^{-1}||_2$ and $||(\hSigma-zI_p)^{-1}||_2$ separately. 
\benum
\item Bounding $||(\hSigma-zI_p)^{-1}||_2$:

$$||(\hSigma-zI_p)^{-1}||_2=\sigma_{\max}((\hSigma-zI_p)^{-1})=\max_i\frac{1}{|\hat{l}_i-z|},$$
where $\hat{l}_i$ is the $i$-th eigenvalue of $\hSigma$. Since $\hat{l}_i$ is always non-negative, we have
$$
\frac{1}{|\hat{l}_i-z|}=\frac{1}{|\hat{l}_i-\textnormal{Re}z-i\textnormal{Im}z|}=\frac{1}{\sqrt{(\hat{l}_i-\textnormal{Re}z)^2+(\textnormal{Im}z)^2}}\leq\frac{1}{\delta}.
$$

\item Bounding $||(x_p\Sigma-zI_p)^{-1}||_2$:

In this case, we need to use the properties of $e_p$ and $x_p$. Recall that $e_p$ is the Stieltjes transform of a certain measure on $\mathbb R^+$, i.e.
\begin{align*}
e_p(z)&=\int_0^{\infty}\frac{1}{t-z}d\mu(t)=\int_0^{\infty}\frac{1}{t-\textnormal{Re}z-i\textnormal{Im}z}d\mu(t)\\
&=\int_0^\infty\frac{t-\textnormal{Re}z}{(t-\textnormal{Re}z)^2+(\textnormal{Im}z)^2}d\mu(t)+i\int_0^\infty\frac{\textnormal{Im}z}{(t-\textnormal{Re}z)^2+(\textnormal{Im}z)^2}d\mu(t).
\end{align*}
So
\begin{align*}
x_p&=\frac{1}{1+\gamma_pe_p}=\frac{1}{1+\gamma_p\textnormal{Re}(e_p)+i\gamma_p\textnormal{Im}(e_p)}\\
&=\frac{1+\gamma_p\textnormal{Re}(e_p)}{(1+\gamma_p\textnormal{Re}(e_p))^2+(\gamma_p\textnormal{Im}(e_p))^2}-i\frac{\gamma_p\textnormal{Im}(e_p)}{(1+\gamma_p\textnormal{Re}(e_p))^2+(\gamma_p\textnormal{Im}(e_p))^2}.
\end{align*}
When $z\in D_\delta$, we can check that $\textnormal{Re}(x_p)>0$. Meanwhile, $\textnormal{Im}(x_p)$ and $\textnormal{Im}(z)$ have opposite signs.

Now, let us consider
$$||(x_p\Sigma-zI_p)^{-1}||_2=\sigma_{\max}((x_p\Sigma-zI_p)^{-1})=\max_k\frac{1}{|x_pl_k-z|},$$
where $l_k$ is the $k$-th eigenvalue of $\Sigma$. Since $l_k$ is non-negative, we have
\begin{align*}
\frac{1}{|x_pl_k-z|}&=\frac{1}{|l_k\textnormal{Re}(x_p)+il_k\textnormal{Im}(x_p)-\textnormal{Re}z-i\textnormal{Im}z|}\\
&=\frac{1}{\sqrt{(l_k\textnormal{Re}(x_p)-\textnormal{Re}z)^2+(l_k\textnormal{Im}(x_p)-\textnormal{Im}z)^2}}\\
&\leq\frac{1}{\delta}.
\end{align*}

\eenum

Finally, since $\delta$ is arbitrary, we can conclude that $f'(z, \Sigma)\asymp g'(z, \hSigma)$ for all $z\in D$.

Then let us compute the derivatives. For invertible $A=A(z)$, we have
$$
\frac{d(A^{-1})}{dz}=-A^{-1}\frac{d A}{d z}A^{-1},
$$
where the derivative is entry-wise. Thus
\begin{align*}
f'(z,T) &= -(x_p T-z I)^{-1}(x_p'T-I)(x_p T-zI_p)^{-1}=-(x_p T-zI_p)^{-2}(x_p'T-I),\\
g'(z,S) &= (S-zI_p)^{-2}. 
\end{align*}

Next, we need to calculate $x' = dx/dz$, where $x(z)$ is the limit of $x_p(z)$. In fact, by looking at the expression of $x_p(z)$, it is not hard to find that $x_p(z)$ is uniformly bounded on $D$. By using a similar argument, we have $x_p'\to x'$ on $D$. To find $x'$, let us start from the following fixed-point equation
\beqs 1 -x = \gamma \left[1 +z\E_H\frac{1}{xT-z}\right].\eeqs 
Take derivatives on both sides to get
\begin{align*}
-x' &= \gamma \left[z\E_H\frac{1}{xT-z}\right]'\\
-x' &= \gamma \left[\E_H\frac{1}{xT-z}
+ \E_H\frac{z-zTx'}{(xT-z)^2}\right]\\
x'\left[-1 + \gamma z\E_H\frac{T}{(xT-z)^2}\right] 
&= \gamma \E_H\frac{xT}{(xT-z)^2}\\
x'
&=\frac{\gamma \E_H\frac{xT}{(xT-z)^2}}
{-1 + \gamma z\E_H\frac{T}{(xT-z)^2}}.
\end{align*}
Therefore we obtain
\begin{align*}
(\hSigma-z I)^{-2} &\asymp  (x_p \Sigma-zI_p)^{-2}(I-x'_p\Sigma)\\
p^{-1} \tr(\hSigma-z I)^{-2} &\asymp  -x'_p p^{-1}\tr[\Sigma(x_p \Sigma -z I)^{-2}]+p^{-1}\tr[(x_p\Sigma -zI_p)^{-2}]\\
&\to \frac{\gamma \E_H\frac{xT}{(xT-z)^2}}
{1 - \gamma z\E_H\frac{T}{(xT-z)^2}}\E_H\frac{T}{(xT-z)^2}+\E_H\frac{1}{(xT-z)^2}\\
&= \frac{\gamma x \left(\E_H\frac{T}{(xT-z)^2}\right)^2}
{1 - \gamma z\E_H\frac{T}{(xT-z)^2}}+ \E_H\frac{1}{(xT-z)^2}.
\end{align*} 

Now, let $z=-\lambda$ and then we will have 
\begin{align*}
(\hSigma+\lambda I)^{-2} &\asymp  (x_p \Sigma+\lambda I)^{-2}(I-x'_p\Sigma)\\
p^{-1} \tr(\hSigma+\lambda I)^{-2} &\to \frac{\gamma x \left(\E_H\frac{T}{(xT+\lambda)^2}\right)^2}{1 +\gamma \lambda\E_H\frac{T}{(xT+\lambda)^2}}+ \E_H\frac{1}{(xT+\lambda)^2}.
\end{align*}
Finally, we can simply replace $\hSigma, \lambda, \gamma, x$ by $\hSigma_i, \lambda_i, \gamma_i, x_i$ to get the desired results.

\eenum

\item Limit of $R$: Recall that $R_i  = n_i^{-1}\sigma^2\tr[(\hSigma_i+\lambda_i I_p)^{-2} \hSigma_i]$. We note $p^{-1}\tr (\hSigma+\lambda I)^{-2} \to m_{F_\gamma}'(-\lambda)$ and $\hSigma(\hSigma+\lambda I)^{-2} 
= (\hSigma+\lambda I)^{-1} - \lambda (\hSigma+\lambda I)^{-2}$, so 
\beqs \frac{\tr [\hSigma(\hSigma+\lambda I)^{-2} ]}{n} \to \gamma[m_{F_\gamma}(-\lambda)-\lambda m_{F_\gamma}'(-\lambda)].\eeqs

Hence 
\beqs
R_{ii}  \to 
\sigma^2\left[\gamma_i[m_{F_{\gamma_i}}(-\lambda_i)-\lambda m_{F_{\gamma_i}}'(-\lambda_i)]\right].
\eeqs

Next, we find a limit in terms of population parameters
\begin{align*}
\hSigma(\hSigma+\lambda I)^{-2} 
&= (\hSigma+\lambda I)^{-1} - \lambda (\hSigma+\lambda I)^{-2}\\
&\asymp (x_p \Sigma + \lambda I)^{-1} - \lambda (x_p \Sigma+ \lambda I)^{-2}(I-x_p'\Sigma)\\
p^{-1}\tr \hSigma(\hSigma+\lambda I)^{-2} 
&\asymp p^{-1} \tr (x_p \Sigma + \lambda I)^{-1} 
- \lambda p^{-1} \tr \left[(I-x'_p\Sigma)(x_p \Sigma + \lambda I)^{-2}\right]\\
&\to \E_H\frac{1}{xT+\lambda}
- \lambda\frac{\gamma x \left(\E_H\frac{T}{(xT+\lambda)^2}\right)^2}
{1 + \gamma\lambda\E_H\frac{T}{(xT+\lambda)^2}}-\E_H\frac{\lambda}{(xT+\lambda)^2}\\
&=\E_H\frac{xT}{(xT+\lambda)^2}-\lambda\frac{\gamma x \left(\E_H\frac{T}{(xT+\lambda)^2}\right)^2}
{1 + \gamma\lambda\E_H\frac{T}{(xT+\lambda)^2}}\\
&=\frac{x\E_H\frac{T}{(xT+\lambda)^2}}{1+\lambda\gamma\E_H\frac{T}{(xT+\lambda)^2}},
\end{align*} 
where we used the differentiation rule of the calculus of deterministic equivalents. Hence we finally find the limit
\beqs
R_{ii}  \to 
\sigma^2\left[\frac{x_i\E_H\frac{T}{(x_iT+\lambda_i)^2}}{1+\lambda_i\gamma_i\E_H\frac{T}{(x_iT+\lambda_i)^2}}
 \right].
\eeqs

\eenum

\subsection{Proof of Theorem \ref{ARE_equal}}
\label{pf:ARE_equal}

Notice that, when the samples are equally distributed and we use the same tuning parameter $\lambda$ for all the local estimators, a direct consequence is that $x_i=x_j=x$ for all $i,j$, where $x$ is the unique solution of the following fixed point equation
\begin{equation*}
1-x=k\gamma\left[1-\lambda\int_0^\infty\frac{dH(t)}{xt+\lambda}\right]=k\gamma\left[1-\E_H\frac{\lambda}{xT+\lambda}\right]=k\gamma(1-\lambda m_{F_{k\gamma}}(-\lambda))=k\gamma(1-\lambda m).
\end{equation*}
In this case, we can express $\mathcal{A}_{ij}$ as
\begin{equation*}
\mathcal{A}_{ij}=\sigma^2\alpha^2\E_H\frac{(xT)^2}{(xT+\lambda)^2}=\sigma^2\alpha^2\int\frac{(xt)^2}{(xt+\lambda)^2}dH(t).
\end{equation*}
In order to express $\mathcal{A}_{ij}$ in terms of the sample quantities, we can start from the following equality
\begin{equation*}
\int\frac{1}{xt+\lambda}dH(t)=m.
\end{equation*}
Take derivatives with respect to $\lambda$, we have 
\begin{equation*}
\int\frac{x't+1}{(xt+\lambda)^2}dH(t)=m'.
\end{equation*}
Rearrange terms, we have
\begin{equation*}
\int\frac{x't+1}{(xt+\lambda)^2}dH(t)=\int\frac{(xt+\lambda-\lambda)\cdot\frac{x'}{x}+1}{(xt+\lambda)^2}dH(t)=\frac{x'}{x}m+\left(1-\frac{\lambda x'}{x}\right)\int\frac{1}{(xt+\lambda)^2}dH(t)=m'.
\end{equation*}
On the other hand, take derivatives with respect to $\lambda$ on the fixed point equation for $x$ gives us
\begin{equation*}
x'=k\gamma(m-\lambda m').
\end{equation*}
So 
\begin{equation*}
\int\frac{1}{(xt+\lambda)^2}dH(t)=\frac{xm'-x'm}{x-\lambda x'}=\frac{(1-k\gamma)m'+2k\gamma\lambda mm'-k\gamma m^2}{1-k\gamma+k\gamma\lambda^2m'}.
\end{equation*}
Then we have
\begin{align*}
\int\frac{(xt)^2}{(xt+\lambda)^2}dH(t)&=\int\frac{(xt+\lambda-\lambda)^2}{(xt+\lambda)^2}dH(t)\\
&=1-2\lambda m+\lambda^2\int\frac{1}{(xt+\lambda)^2}dH(t)\\
&=1-2\lambda m+\lambda^2m'-\frac{k\gamma\lambda^2(m-\lambda m')^2}{1-k\gamma+k\gamma\lambda m'}.
\end{align*}
Now we the expressions for $V, \mathcal{A},\mathcal{R}$, we also know $\mathcal{W}_k^*=(\mathcal{A}+\mathcal{R})^{-1}V$ and $\mathcal{M}_k=\sigma^2\alpha^2-V^\top(\mathcal{A}+\mathcal{R})^{-1}V$. By using the auxiliary functions $\mathcal{F}, \mathcal{G}$ defined in the theorem, we have
\begin{equation*}
\mathcal{M}_k=\sigma^2\alpha^2-\frac{1}{\mathcal{G}}V^\top ({\bf1\cdot 1}^\top+\diag(\mathcal{F}/\mathcal{G}))^{-1}V=\sigma^2\alpha^2-\frac{(\sigma^2\alpha^2(1-\lambda m))^2}{\mathcal{G}}{\bf 1}^\top ({\bf1\cdot 1}^\top+\diag(\mathcal{F}/\mathcal{G}))^{-1}{\bf 1},
\end{equation*}
where ${\bf 1}=(1,1,\cdots,1)^\top$ is the all-one vector. Then similar to the proof of Theorem \ref{ARE}, we can use the Sherman-Morrison formula to simply the expression, this leads to
\begin{align*}
\mathcal{M}_k&=\sigma^2\alpha^2-\frac{(\sigma^2\alpha^2(1-\lambda m))^2}{\mathcal{G}}{\bf 1}^\top\left(\diag(\mathcal{F}/\mathcal{G})^{-1}-\frac{\diag(\mathcal{F}/\mathcal{G})^{-1}{\bf 1\cdot 1^\top}\diag(\mathcal{F}/\mathcal{G})^{-1}}{1+{\bf 1}^\top\diag(\mathcal{F}/\mathcal{G})^{-1}{\bf 1}}\right){\bf 1}\\
&=\sigma^2\alpha^2-\frac{(\sigma^2\alpha^2(1-\lambda m))^2}{\mathcal{G}}\left(\frac{k\mathcal{G}}{\mathcal{F}}-\frac{k^2\mathcal{G}^2/\mathcal{F}^2}{1+k\mathcal{G}/\mathcal{F}}\right)\\
&=\sigma^2\alpha^2\left(1-\frac{\sigma^2\alpha^2(1-\lambda m)^2k}{\mathcal{F}+k\mathcal{G}}\right).
\end{align*}
Similarly, we can express the optimal weights $\mathcal{W}_k^*$ as
\begin{align*}
\mathcal{W}_k^*&=\frac{1}{\mathcal{G}}({\bf 1\cdot1^\top}+\diag(\mathcal{F}/\mathcal{G}))^{-1}V\\
&=\frac{\sigma^2\alpha^2(1-\lambda m)}{\mathcal{G}}({\bf 1\cdot1^\top}+\diag(\mathcal{F}/\mathcal{G}))^{-1}{\bf 1}\\
&=\frac{\sigma^2\alpha^2(1-\lambda m)}{\mathcal{G}}\left(\diag(\mathcal{F}/\mathcal{G})^{-1}-\frac{\mathcal{G}^2/\mathcal{F}^2}{1+k\mathcal{G}/\mathcal{F}}{\bf 1\cdot1^\top}\right){\bf 1}\\
&=\frac{\sigma^2\alpha^2(1-\lambda m)}{\mathcal{F}+k\mathcal{G}}{\bf 1}.
\end{align*}

\subsection{Gaussian MLE for signal and noise components}
\label{signoise}

Recall that our model is $Y=X\beta+\ep$ where $\beta$ and $\ep$ are independent. Let $\theta=(\sigma^2,\alpha^2)$ and define the Gaussian log-likelihood,
\beqs
\ell(\theta)=-\frac{1}{2}\log(\sigma^2)-\frac{1}{2n}\log\det\left(\frac{\alpha^2}{p}XX^\top+I\right)-\frac{1}{2\sigma^2n}Y^\top \left(\frac{\alpha^2}{p}XX^\top+I\right)^{-1}Y.
\eeqs
Note that $\ell(\theta)$ is the log-likelihood for $\theta$ under the Gaussian assumption of $\beta\sim\mathcal{N}(0,(\sigma^2\alpha^2/p)I)$ and $\ep\sim\mathcal{N}(0, \sigma^2I)$. For the MLE
\beqs
\htheta=(\hsigma^2,\hat{\alpha}^2)=\underset{\sigma^2,\alpha^2\geq 0}\argmax ~\ell(\theta),
\eeqs
we have the following result from \cite{dicker2017variance}.
\begin{theorem}[Consistency and asymptotic normality of the MLE, \cite{dicker2017variance}]
\label{mle}
Suppose $\theta=(\sigma^2, \alpha^2)$ are the true parameters, then $\htheta\to\theta$ in probability as $p/n\to\gamma$. Furthermore, define the Fisher information matrix for $\theta$ under the Gaussian assumption model
$$
\mathcal{I}_n(\theta)=\begin{bmatrix}
    I_2(\theta) & I_3(\theta) \\
    I_3(\theta) & I_4(\theta)
\end{bmatrix},
$$
where 
$$
I_k(\theta)=\frac{1}{2n\sigma^{8-2k}}\tr\left[\left(\frac{1}{p}XX^\top\right)^{k-2}\left(\frac{\alpha^2}{p}XX^\top+I\right)^{2-k}\right],~~k=2,3,4.
$$
Then $n^{1/2}\mathcal{I}_n(\theta)^{1/2}(\htheta-\theta)\to\mathcal{N}(0,I_2)$ in distribution as $p/n\to\gamma$.
\end{theorem}

In addition, if we put some assumptions on $X$ as we did in Theorem \ref{mp} and denote the limiting spectral distribution of $p^{-1}XX^\top$ by $F_\gamma$, then the entries of the Fisher information matrix $\mathcal{I}_n(\theta)$ have limits
$$
I_k(\theta)\to_{a.s.} \mathcal{J}_k(\theta)=\frac{1}{2\sigma^{8-2k}}\int\left(\frac{s}{\alpha^2s+1}\right)^{k-2}dF_\gamma(s),~~k=2,3,4.
$$
Thus $\mathcal{I}_n(\theta)$ converges almost surely to a limiting information matrix $\mathcal{I}(\theta)$ which characterizes the asymptotic variance of the MLE $\htheta$.

\subsection{Proof of Theorem \ref{AREsimple}}
\label{pf:AREsimple}

The proof for $v$ and $R$ is clear by Theorem \ref{ARE}. For the limit of $A$, the diagonal case is also direct. When $i\neq j$, recall that
\beqs E_{ij}=p^{-1}\tr \{(\hSigma_i+\lambda_i I_p)^{-1}(\hSigma_j+\lambda_j I_p)^{-1}\} \to \E_H \frac{1}{(x_iT+\lambda_i)(x_jT+\lambda_j)}.\eeqs
For $H=\delta_1$, the expectation decouples, we find
\beqs E_{ij} \to  \frac{1}{x_i+\lambda_i} \cdot \frac{1}{x_i+\lambda_j}= m_{\gamma_i}(-\lambda_i) m_{\gamma_j}(-\lambda_j).\eeqs 
Therefore, 
\beqs A_{ij} \to  
\sigma^2\alpha^2[1- \lambda_im_{\gamma_i}(-\lambda_i)]
\cdot 
[1- \lambda_jm_{\gamma_j}(-\lambda_j)].\eeqs 
Now let us put everything together. Recall that the optimal risk has the form $\textnormal{MSE}^*_{dist} = \|\beta\|^2- v^\top (A+R)^{-1} v$. Based on the above discussion, we have 
\beqs
\sigma^2\alpha^2(A+R) \to \sigma^2\alpha^2(\mathcal{A}+\mathcal{R})=V V^\top+D,
\eeqs
where $D$ is a diagonal matrix with $i$-th diagonal entry $\sigma^2\alpha^2(\mathcal{R}_{ii} + \mathcal{A}_{ii})-V_i^2$. Then, by using the Sherman$-$Morrison formula, we have 
\beqs
V^\top(VV^\top+D)^{-1}V
=\frac{V^\top D^{-1}V}{1+V^\top D^{-1}V}.
\eeqs
So the limiting distributed risk is
\begin{align*}
\mathcal{M}_k= \sigma^2\alpha^2 - \sigma^2\alpha^2\frac{V^\top D^{-1}V}{1+V^\top D^{-1}V}
= \frac{\sigma^2\alpha^2}{1+V^\top D^{-1}V}=\frac{\sigma^2\alpha^2}{1+\sum_{i=1}^k\frac{V_i^2}{D_i}},
\end{align*}
which finishes the proof.

\subsection{Explaining decoupling via free probability theory}
\label{free_prob_exp}

In this section, we provide an explanation via free probability theory for why the limiting distributed risk decouples. Specifically, we explain why the limit of the quantities $\beta^\top Q_i\beta\cdot\beta^\top Q_j\beta$ for $i\neq j$ becomes a product of terms depending on $i,j$.

We will use some basic notions from free probability theory \citep{voiculescu1992free,hiai2006semicircle,nica2006lectures,anderson2010introduction,couillet2011random}. Let us define our non-commutative probability space as
$$\left(\mathcal A=(L^{\infty-}\otimes M_p(\mathbb R)), \tau=\frac{1}{p}\tr \right),$$
where $L^{\infty-}$ denotes the collection of random variables with all moments finite and $M_p(\mathbb R)$ is the space of $p\times p$ real matrices. Recall that, a sequence of random variables $\{a_{1,p},a_{2,p}, \dots, a_{k,p}\}\subset \mathcal A$ is said to be asymptotically free almost surely if 
$$
\tau[\prod_{j=1}^mP_j(a_{i_j,p}-\tau(P_j(a_{i_j,p})))]\to_{a.s.}0,
$$
for any positive integer $m$, any polynomials $P_1,\dots, P_m$ and any indices $i_1,\dots,i_m\in[k]$ with no two adjacent $i_j$ equal. Suppose $A_p, B_p$ are two sequences of independent random matrices and at least one of them is orthogonally invariant, then it is well-known that $\{A_p, B_p\}\subset \mathcal A$ is asymptotically free almost surely. 

Now, let us assume that $X^\top X$ is orthogonally invariant, which is the case when $X^\top X$ follows the white Wishart distribution. Then clearly $X_i^\top X_i$ and $X_j^\top X_j$ are asymptotically free almost surely. It follows that $Q_i$ and $Q_j$ are also asymptotically free almost surely. By using the definition of asymptotic freeness, we have for $i\neq j$
$$\tau[(Q_i-\frac{1}{p}\tr(Q_i)I)(Q_j-\frac{1}{p}\tr(Q_j)I)]\to_{a.s.}0,$$
which is equivalent to 
$$\frac{1}{p}\tr(Q_iQ_j)-\frac{1}{p}\tr(Q_i)\frac{1}{p}\tr(Q_j)\to_{a.s.}0.$$
Hence, under the random-effects assumption for $\beta$, the limit of $\beta^\top\beta\cdot\beta^\top Q_iQ_j\beta$ ($i\neq j$) will decouple and is the same as the limit of $\beta^\top Q_i\beta\cdot\beta^\top Q_j\beta$.

\subsection{Proof of Proposition \ref{parameter}}
\label{pf:parameter}

Recall that
\begin{align*}
\frac{V_i^2}{D_i}
 &=\frac{\sigma^4\alpha^4(1-\lambda_im_{\gamma_i}(-\lambda_i))^2}
{\sigma^4\alpha^4\lambda_i^2[m'_{\gamma_i}(-\lambda_i)-m^2_{\gamma_i}(-\lambda_i)]
+\sigma^4\alpha^2\gamma_i[m_{\gamma_i}(-\lambda_i)-\lambda_im'_{\gamma_i}(-\lambda_i)]}\\
 &=\frac{\alpha^2(1-\lambda_im_{\gamma_i}(-\lambda_i))^2}
{\alpha^2\lambda_i^2[m'_{\gamma_i}(-\lambda_i)-m^2_{\gamma_i}(-\lambda_i)]
+\gamma_i[m_{\gamma_i}(-\lambda_i)-\lambda_im'_{\gamma_i}(-\lambda_i)]},
\end{align*}
and our goal is to find $\lambda_i$ that maximizes $V_i^2/D_i$. Luckily, from \cite{dobriban2018high} it follows that for $k=1$, i.e. when there is only one machine, the optimal choice of the tuning parameter $\lambda$ is $\gamma/\alpha^2$. This means that the maximizer of $V^2/D$ is $\lambda=\gamma/\alpha^2$. Now, due to the decoupled structure of $\mathcal{M}_k$, the optimal tuning parameters are $\lambda_i=\gamma_i/\alpha^2$. Plugging in the parameters, we have 
\begin{align*}
\frac{V_i^2}{D_i} &
=\frac{\alpha^2}{\gamma_im_{\gamma_i}(-\gamma_i/\alpha^2)}-1.
\end{align*}
Then the optimal risk can be simplified to
\beqs
\mathcal{M}_k=\frac{\sigma^2\alpha^2}
{1+\sum_{i=1}^k\left[\frac{\alpha^2}{\gamma_im_{\gamma_i}(-\gamma_i/\alpha^2)}-1\right]}.
\eeqs
When $k=1$, this equals to $\sigma^2\gamma m_{\gamma}(-\gamma/\alpha^2)$ which matches the known result from \cite{dobriban2018high}.

\subsection{Proof of Proposition \ref{optimalrisk}}
\label{pf:optimalrisk}

The explicit form is easy to derive by plugging $z=-\gamma/\alpha^2$ into the formula of $m_\gamma(z)$. Next, we can check monotonicity by computing $\phi'(\gamma)$:
\beqs
\phi'(\gamma)=\frac{\alpha^2}{2\gamma^2}\left(1+\frac{(1-1/\alpha^2)\gamma-1}{\sqrt{[(1-1/\alpha^2)\gamma-1]^2+4\gamma^2/\alpha^2}}\right)>0.
\eeqs
Finally, for the convexity, let us consider the two cases separately.
\benum
\item $\alpha\leq 1$:  With some effort, we find the second derivative of $\phi$
\begin{align*}
\phi''(\gamma)=\frac{\alpha^2\left(\frac{2\gamma^2}{\alpha^2}-\left(((1-\frac{1}{\alpha^2})\gamma-1)^2+\frac{4\gamma^2}{\alpha^2}\right)\left((1-\frac{1}{\alpha^2})\gamma-1)\right)-\left(((1-\frac{1}{\alpha^2})\gamma-1)^2+\frac{4\gamma^2}{\alpha^2}\right)^{3/2}\right)}{\gamma^3[((1-1/\alpha^2)\gamma-1)^2+4\gamma^2/\alpha^2]^{3/2}}.
\end{align*}
To analyze the second derivative, it is helpful to denote $1-(1-\frac{1}{\alpha^2})\gamma$ by $\Delta$. Clearly, in this case, $\Delta\geq1$. Then 
we can rewrite $\phi''$ as 
\begin{align*}
\phi''(\gamma)&=\frac{\alpha^2}{\gamma^3[\Delta^2+4\gamma^2/\alpha^2]^{3/2}}\left(\frac{2\gamma^2}{\alpha^2}+(\Delta^2+\frac{4\gamma^2}{\alpha^2})\Delta-(\Delta^2+\frac{4\gamma^2}{\alpha^2})^{3/2}\right)\\
&=\frac{\alpha^2}{\gamma^3[\Delta^2+4\gamma^2/\alpha^2]^{3/2}}\left(\frac{2\gamma^2}{\alpha^2}+(\Delta^2+\frac{4\gamma^2}{\alpha^2})\left(\Delta-\sqrt{\Delta^2+\frac{4\gamma^2}{\alpha^2}}\right)\right)\\
&=\frac{\alpha^2}{\gamma^3[\Delta^2+4\gamma^2/\alpha^2]^{3/2}}\left(\frac{2\gamma^2}{\alpha^2}-\frac{4\gamma^2}{\alpha^2}\cdot\frac{\Delta^2+4\gamma^2/\alpha^2}{\Delta+\sqrt{\Delta^2+\frac{4\gamma^2}{\alpha^2}}}\right)\\
&=\frac{\alpha^2}{\gamma^3[\Delta^2+4\gamma^2/\alpha^2]^{3/2}}\left(\frac{2\gamma^2}{\alpha^2}-\frac{4\gamma^2}{\alpha^2}\cdot\frac{\sqrt{\Delta^2+4\gamma^2/\alpha^2}}{\Delta+\sqrt{\Delta^2+\frac{4\gamma^2}{\alpha^2}}}\right)\\
&\leq\frac{\alpha^2}{\gamma^3[\Delta^2+4\gamma^2/\alpha^2]^{3/2}}\left(\frac{2\gamma^2}{\alpha^2}-\frac{4\gamma^2}{\alpha^2}\cdot\frac{1}{2}\right)=0.
\end{align*}
Thus, $\phi(\gamma)$ is always concave in this case.

\item $\alpha>1$: Here we can consider the Taylor expansion of $\phi''$ near the origin. We can check that $\phi''(\gamma)=2(1-1/\alpha^2)\gamma^3+o(\gamma^3)$ as $\gamma\to 0$, which means $\phi''(\gamma)>0$ for small $\gamma$. When $\gamma$ is very large, we can immediately see that $\phi''(\gamma)<0$, since the leading order in the numerator of $\phi''(\gamma)$ is $-\gamma^3$. Then the desired result follows.

\eenum

\subsection{Proof of Theorem \ref{AREprop}}
\label{pf:AREprop}

For the first property, minimizing the ARE is equivalent to maximizing the following quantity
\beqs
\sum_{i=1}^k \frac{\gamma_im_{\gamma_i}(-\gamma_i/\alpha^2)}{\alpha^2}=\sum_{i=1}^k\frac{\phi(\gamma_i)}{\alpha^2}.
\eeqs
It is helpful to introduce $r(t)=\phi(\gamma)$, where $t=1/\gamma$. We can easily compute that 
$$
r'(t)=\frac{\alpha^2}{2}\left(-1+\frac{t-(1-1/\alpha^2)}{\sqrt{(t-(1-1/\alpha^2))^2+4/\alpha^2}}\right)<0~~,r''(t)=\frac{2}{[(t-(1-1/\alpha^2))^2+4/\alpha^2]^{3/2}}>0.
$$
Thus, $r(t)$ is a decreasing and convex function. We can show the ARE achieves minimum when the samples are equally distributed by considering the following optimization problem
\begin{equation*}
\begin{aligned}
& \underset{t_i}{\text{max}}
& & \sum_{i=1}^k \frac{r(t_i)}{\alpha^2}\\
& \text{subject to}
& & \sum_{i=1}^kt_i=\frac{1}{\gamma},\\
& & & t_i\geq0, i=1,2,\dots,k.
\end{aligned}
\end{equation*}
We denote the objective by $R(t_1,\dots,t_k)$, 
and the corresponding Lagrangian by $R_\xi=R-\xi(\sum_it_i-1/\gamma)$. Then it is easy to check that the condition $\frac{\partial R_\xi}{\partial t_i}=0$ reduces to
\beqs
\frac{r'(t_i)}{\alpha^2}-\xi=0,~i=1,2,\dots,k.
\eeqs
Since $r'(t)$ is also monotone, the unique solution to the stationary condition is $t_1=t_2=\cdots=t_k=1/(k\gamma)$. If some $t_i$ equals to $0$, then it reduces to a problem with $k-1$ machines. So it remains to check the boundary case where only one $t_i$ is non-zero and equals to $1/\gamma$. Obviously, this is the trivial case where the ARE is $1$. Therefore, we have shown that the ARE attains its minimum when the samples are equally distributed across $k$ machines.

Next, for fixed $\alpha^2$ and $\gamma$, we can check
\beqs
\frac{\partial \psi}{\partial k}=\frac{\gamma m_{\gamma}(-\gamma/\alpha^2)}{\alpha^2}\left(\frac{\alpha^2}{2\gamma}\cdot\frac{\left(\gamma/\alpha^2+\gamma\right)^2k+\gamma/\alpha^2-\gamma}{\sqrt{\left(\gamma/\alpha^2+\gamma\right)^2k^2+2\left(\gamma/\alpha^2-\gamma\right)k+1}}-\frac{1+\alpha^2}{2}\right)\leq 0.
\eeqs
Moreover, the limit of $\psi$ is
\begin{align*}
h(\alpha^2, \gamma) = \lim_{k\to\infty}\psi(k,\gamma,\alpha^2)&=\frac{\gamma m_{\gamma}(-\gamma/\alpha^2)}{\alpha^2}\left(1+\frac{\alpha^2}{\gamma(1+\alpha^2)}\right)\\
&=\frac{-\gamma/\alpha^2+\gamma-1+\sqrt{(-\gamma/\alpha^2+\gamma-1)^2+4\gamma^2/\alpha^2}}{2\gamma}\left(1+\frac{\alpha^2}{\gamma(1+\alpha^2)}\right).
\end{align*}
Then for fixed $\alpha^2$, we can differentiate $h(\alpha^2, \gamma)$ with respect to $\gamma$:
\begin{align*}
\frac{\partial h}{\partial\gamma}&=-\frac{\alpha^2}{\gamma^2(1+\alpha^2)}\cdot\frac{-\gamma/\alpha^2+\gamma-1+\sqrt{(-\gamma/\alpha^2+\gamma-1)^2+4\gamma^2/\alpha^2}}{2\gamma}\\
&+\left(1+\frac{\alpha^2}{\gamma(1+\alpha^2)}\right)\cdot\frac{1-1/\alpha^2+\frac{(-\gamma/\alpha^2+\gamma-1)(1-1/\alpha^2)+4\gamma/\alpha^2}{\sqrt{(-\gamma/\alpha^2+\gamma-1)^2+4\gamma^2/\alpha^2}}}{2\gamma}\\
&-\left(1+\frac{\alpha^2}{\gamma(1+\alpha^2)}\right)\cdot\frac{-\gamma/\alpha^2+\gamma-1+\sqrt{(-\gamma/\alpha^2+\gamma-1)^2+4\gamma^2/\alpha^2}}{2\gamma^2}.
\end{align*}
After tedious calculation, we find $\frac{\partial h}{\partial\gamma}\geq 0$.
Finally, we can evaluate the limit of $h$ as $\gamma\to 0$ and $\gamma\to\infty$
\beqs
\lim_{\gamma\to0}h(\alpha^2, \gamma)=\frac{1}{1+\alpha^2},~~~\lim_{\gamma\to\infty}h(\alpha^2, \gamma)=1.
\eeqs
On the other hand, for fixed $\gamma$, we can check that $h$ is a decreasing function of $\alpha^2$  and
\beqs
\lim_{\alpha^2\to0}h(\alpha^2, \gamma)=1,~~
\lim_{\alpha^2\rightarrow\infty}h(\alpha^2, \gamma)
=
\begin{cases}
1-\frac{1}{\gamma^2},~~\gamma>1,\\
0,~~0<\gamma\leq1.
\end{cases}
\eeqs

\subsection{Proof of Theorem \ref{weights}}
\label{pf:weights}

Recall that the optimal weights are $w^*=(A+R)^{-1}v$ and $\sigma^2\alpha^2(A+R)\to VV^\top+D$. Denote the limit of the optimal weights by $W$, so that we have 
$$
W=\sigma^2\alpha^2(VV^\top+D)^{-1}V=\frac{\sigma^2\alpha^2D^{-1}V}{1+V^\top D^{-1}V}.
$$
When we choose $\lambda_i=\gamma_i/\alpha^2$ for each $i$, we can write the limiting optimal weights as
$$
W=\mathcal{M}_k\cdot D^{-1}V.$$
So, it follows from the formulas of $\mathcal{M}_k, D$ and $V$ that
$$
W_i=\left(\frac{\alpha^2}{\gamma_im_{\gamma_i}(-\gamma_i/\alpha^2)}\right)\cdot \left(\frac{1}{1+\sum_{i=1}^k\left[\frac{\alpha^2}{\gamma_im_{\gamma_i}(-\gamma_i/\alpha^2)}-1\right]}\right).
$$
For the sum of the coordinates, we have
$$
1^\top W=\frac{\sum_{i=1}^k\left(\frac{\alpha^2}{\gamma_im_{\gamma_i}(-\gamma_i/\alpha^2)}\right)}{1+\sum_{i=1}^k\left[\frac{\alpha^2}{\gamma_im_{\gamma_i}(-\gamma_i/\alpha^2)}-1\right]}=\frac{\sum_{i=1}^k\left(\frac{\alpha^2}{\gamma_im_{\gamma_i}(-\gamma_i/\alpha^2)}\right)}{1-k+\sum_{i=1}^k\left(\frac{\alpha^2}{\gamma_im_{\gamma_i}(-\gamma_i/\alpha^2)}\right)}\geq 1.
$$

In the special case where all $\gamma_i$ are equal, i.e., $\gamma_i = k\gamma$, we have all $W_i$ equal to
$$
W_i=\frac{\frac{\alpha^2}{k\gamma \cdot m_{k\gamma}(-k\gamma/\alpha^2)}}
{1-k+\frac{\alpha^2}{\gamma \cdot m_{k\gamma}(-k\gamma/\alpha^2)}}
=\frac{1}{k+ (1-k)\cdot k\gamma/\alpha^2 \cdot m_{k\gamma}(-k\gamma/\alpha^2)}
.
$$
In terms of the optimal risk function $\phi(\gamma)=\phi(\gamma,\alpha)=\gamma m_{\gamma}(-\gamma/\alpha^2)$ defined before, this can also be written as the following optimal weight function
$$
\mathcal{W}(k,\gamma, \alpha)=
\frac{1}{k- (k-1)\cdot \phi(k\gamma)/\alpha^2}.
$$
The monotonicity and the limits of $\mathcal{W}$ can be checked directly.

\subsection{Intuitive explanation for the need to use weights summing to greater than unity}
\label{int_exp}
Consider a much more simplified problem, where we are estimating a scalar parameter $\theta$. We have an estimator $\htheta$, which is generally biased, and we would like to find the scale multiple $c\cdot \htheta$ that minimizes the mean squared error. A calculation reveals that 

$$M(c) = \E(c\cdot \htheta-\theta)^2 
= c^2\E(\htheta^2)-2c\cdot\E\htheta\cdot\theta+\theta^2$$

Hence the optimal scale factor is $c = \E\htheta\cdot\theta/ \E(\htheta^2)$. 

We can achieve a better understanding of this optimal scale if we consider the bias-variance decomposition of $\htheta$. Let us define the bias and the variance as
\begin{align*}
B &= \E\htheta-\theta\\
V &= \E(\htheta-\E\htheta)^2
\end{align*}
We then see that the optimal scale factor is
\begin{align*}
c &= \frac{B+\theta}{V+(B+\theta)^2}\theta 
= 1 - \frac{V+B(B+\theta)}{V+(B+\theta)^2}.
\end{align*}
This quantity is an ``inflation factor", i.e., greater than or equal to unity, if $V+B(B+\theta)\le 0$. This can be written as
\begin{align*}
V+B^2\le -B\theta
\end{align*}
Hence, this condition can only hold if the bias $B$ has opposite sign with $\theta$. This would be the case for a \emph{shrinkage estimator $\theta$}. In that case, the condition could hold if the parameter $\theta$ has a large magnitude. 

Returning to our main problem, ridge regression is a shrinkage estimator, and averages of ridge regression estimators are still shrinkage estimators. Therefore, it makes sense that their weighted average should be inflated to minimize mean squared error. This provides an intuitive explanation for why the weights sum to greater than one. 

\subsection{Proof of Proposition \ref{oe_ridge}}
\label{pf:oe}

Recall the definition of the out-of-sample prediction error is $\E\|y_t-x_t^\top\hbeta\|^2$. So for any estimator $\hbeta$, under the assumption $\Sigma=I$, we have
\begin{align*}
\E\|y_t-x_t^\top\hbeta\|^2&=\E\|x_t^\top(\hbeta-\beta)+\ep_t\|^2=\E\|x_t^\top(\hbeta-\beta)\|^2+\sigma^2\\
&=\E[(\hbeta-\beta)^\top x_t\cdot x_t^\top(\hbeta-\beta)]+\sigma^2\\
&=\E[(\hbeta-\beta)^\top \Sigma(\hbeta-\beta)]+\sigma^2\\
&=\E\|\hbeta-\beta\|^2+\sigma^2.
\end{align*}
When we consider the distributed estimator and take the limit, we obtain
$$
\mathcal{O}_k=\sigma^2+\mathcal{M}_k,
$$
and the formula for OE. For the inequality between OE and ARE, it is sufficient to notice that $ARE\leq 1$. Finally, the explicit formulas follow easily from previous results.

\subsection{Proof of Theorem \ref{minimax}}
\label{pf:minimax}

It is equivalent to show that the ARE is always greater than or equal to $1/(1+\alpha^2)$. To do this, we need to use Theorem \ref{AREprop}. From the first property, we have ARE$\geq \psi(k,\gamma,\alpha^2)$. Then, since $\psi$ is a decreasing function of $k$, it is lower bounded by its limit at infinity, which is $h(\alpha^2,\gamma)$. Finally, $h(\alpha^2,\gamma)$ is an increasing function of $\gamma$, so it is lower bounded by the limit at $0$, which is $1/(1+\alpha^2)$. When $\gamma>1$, $h(\alpha^2,\gamma)$ is a decreasing function of $\alpha^2$, so it is lower bounded by the limit at infinity, which is $1-1/\gamma^2$. The desired result follows.

{\small
\setlength{\bibsep}{0.2pt plus 0.3ex}
\bibliographystyle{plainnat-abbrev}
\bibliography{references}

\begin{thebibliography}{64}
\providecommand{\natexlab}[1]{#1}
\providecommand{\url}[1]{\texttt{#1}}
\expandafter\ifx\csname urlstyle\endcsname\relax
  \providecommand{\doi}[1]{doi: #1}\else
  \providecommand{\doi}{doi: \begingroup \urlstyle{rm}\Url}\fi

\bibitem[Anderson et~al.(2010)Anderson, Guionnet, and
  Zeitouni]{anderson2010introduction}
G.~W. Anderson, A.~Guionnet, and O.~Zeitouni.
\newblock \emph{An Introduction to Random Matrices}.
\newblock Number 118. Cambridge University Press, 2010.

\bibitem[Anderson(2003)]{anderson1958introduction}
T.~W. Anderson.
\newblock \emph{An Introduction to Multivariate Statistical Analysis}.
\newblock Wiley New York, 2003.

\bibitem[Bai and Silverstein(2009)]{bai2009spectral}
Z.~Bai and J.~W. Silverstein.
\newblock \emph{Spectral analysis of large dimensional random matrices}.
\newblock Springer Series in Statistics. Springer, 2009.

\bibitem[Banerjee et~al.(2016)Banerjee, Durot, and Sen]{banerjee2016divide}
M.~Banerjee, C.~Durot, and B.~Sen.
\newblock Divide and conquer in non-standard problems and the super-efficiency
  phenomenon.
\newblock \emph{arXiv preprint arXiv:1605.04446}, 2016.

\bibitem[Battey et~al.(2018)Battey, Fan, Liu, Lu, and
  Zhu]{battey2018distributed}
H.~Battey, J.~Fan, H.~Liu, J.~Lu, and Z.~Zhu.
\newblock Distributed testing and estimation under sparse high dimensional
  models.
\newblock \emph{The Annals of Statistics}, 46\penalty0 (3):\penalty0
  1352--1382, 2018.

\bibitem[Bertin-Mahieux et~al.(2011)Bertin-Mahieux, Ellis, Whitman, and
  Lamere]{Bertin-Mahieux2011}
T.~Bertin-Mahieux, D.~P. Ellis, B.~Whitman, and P.~Lamere.
\newblock The million song dataset.
\newblock In \emph{{Proceedings of the 12th International Conference on Music
  Information Retrieval ({ISMIR} 2011)}}, 2011.

\bibitem[Bertsekas and Tsitsiklis(1989)]{bertsekas1989parallel}
D.~P. Bertsekas and J.~N. Tsitsiklis.
\newblock \emph{Parallel and distributed computation: numerical methods},
  volume~23.
\newblock Prentice hall Englewood Cliffs, NJ, 1989.

\bibitem[Blelloch and Maggs(2010)]{blelloch2010parallel}
G.~E. Blelloch and B.~M. Maggs.
\newblock Parallel algorithms.
\newblock In \emph{Algorithms and theory of computation handbook}, pages
  25--25. Chapman \& Hall/CRC, 2010.

\bibitem[Boyd et~al.(2011)Boyd, Parikh, Chu, Peleato, Eckstein,
  et~al.]{boyd2011distributed}
S.~Boyd, N.~Parikh, E.~Chu, B.~Peleato, J.~Eckstein, et~al.
\newblock Distributed optimization and statistical learning via the alternating
  direction method of multipliers.
\newblock \emph{Foundations and Trends{\textregistered} in Machine learning},
  3\penalty0 (1):\penalty0 1--122, 2011.

\bibitem[Braverman et~al.(2016)Braverman, Garg, Ma, Nguyen, and
  Woodruff]{braverman2016communication}
M.~Braverman, A.~Garg, T.~Ma, H.~L. Nguyen, and D.~P. Woodruff.
\newblock Communication lower bounds for statistical estimation problems via a
  distributed data processing inequality.
\newblock In \emph{Proceedings of the forty-eighth annual ACM symposium on
  Theory of Computing}, pages 1011--1020. ACM, 2016.

\bibitem[Cai and Wei(2020)]{cai2020distributed}
T.~T. Cai and H.~Wei.
\newblock Distributed gaussian mean estimation under communication constraints:
  Optimal rates and communication-efficient algorithms.
\newblock \emph{arXiv preprint arXiv:2001.08877}, 2020.

\bibitem[Chen et~al.(2018{\natexlab{a}})Chen, Liu, and
  Zhang]{xichen2018Quantile}
X.~Chen, W.~Liu, and Y.~Zhang.
\newblock Quantile regression under memory constraint.
\newblock \emph{arXiv preprint arxiv:1810.08264}, 2018{\natexlab{a}}.

\bibitem[Chen et~al.(2018{\natexlab{b}})Chen, Liu, and Zhang]{xichen2018fone}
X.~Chen, W.~Liu, and Y.~Zhang.
\newblock First-order newton-type estimator for distributed estimation and
  inference.
\newblock \emph{arXiv preprint arxiv:1811.11368}, 2018{\natexlab{b}}.

\bibitem[Chen and Xie(2014)]{chen2014split}
X.~Chen and M.-g. Xie.
\newblock A split-and-conquer approach for analysis of extraordinarily large
  data.
\newblock \emph{Statistica Sinica}, pages 1655--1684, 2014.

\bibitem[Couillet and Debbah(2011)]{couillet2011random}
R.~Couillet and M.~Debbah.
\newblock \emph{Random {M}atrix {M}ethods for {W}ireless {C}ommunications}.
\newblock Cambridge University Press, 2011.

\bibitem[Couillet et~al.(2011)Couillet, Debbah, and
  Silverstein]{couillet2011deterministic}
R.~Couillet, M.~Debbah, and J.~W. Silverstein.
\newblock A deterministic equivalent for the analysis of correlated mimo
  multiple access channels.
\newblock \emph{IEEE Trans. Inform. Theory}, 57\penalty0 (6):\penalty0
  3493--3514, 2011.

\bibitem[Dicker and Erdogdu(2017)]{dicker2017variance}
L.~Dicker and M.~Erdogdu.
\newblock Flexible results for quadratic forms with applications to variance
  components estimation.
\newblock \emph{The Annals of Statistics}, 45\penalty0 (1):\penalty0 386--414,
  2017.

\bibitem[Dicker(2016)]{dicker2014ridge}
L.~H. Dicker.
\newblock Ridge regression and asymptotic minimax estimation over spheres of
  growing dimension.
\newblock \emph{Bernoulli}, 22\penalty0 (1):\penalty0 1--37, 2016.

\bibitem[Dicker(2014)]{dicker2014variance}
L.~H. Dicker.
\newblock Variance estimation in high-dimensional linear models.
\newblock \emph{Biometrika}, 101\penalty0 (2):\penalty0 269--284, 2014.

\bibitem[Dicker and Erdogdu(2016)]{dicker2016variance}
L.~H. Dicker and M.~A. Erdogdu.
\newblock Maximum likelihood for variance estimation in high-dimensional linear
  models.
\newblock In \emph{Proceedings of the 19th International Conference on
  Artificial Intelligence and Statistics}, volume~51 of \emph{Proceedings of
  Machine Learning Research}, pages 159--167. PMLR, 2016.

\bibitem[Dobriban and Sheng(2018)]{dobriban2018Distributed}
E.~Dobriban and Y.~Sheng.
\newblock Distributed linear regression by averaging.
\newblock \emph{arXiv preprint arxiv:1810.00412}, 2018.

\bibitem[Dobriban and Wager(2018)]{dobriban2018high}
E.~Dobriban and S.~Wager.
\newblock High-dimensional asymptotics of prediction: Ridge regression and
  classification.
\newblock \emph{The Annals of Statistics}, 46\penalty0 (1):\penalty0 247--279,
  2018.

\bibitem[Duan et~al.(2018)Duan, Qiao, and Cheng]{duan2018distributed}
J.~Duan, X.~Qiao, and G.~Cheng.
\newblock Distributed nearest neighbor classification.
\newblock \emph{arXiv preprint arXiv:1812.05005}, 2018.

\bibitem[Duchi et~al.(2014)Duchi, Jordan, Wainwright, and
  Zhang]{duchi2014optimality}
J.~C. Duchi, M.~I. Jordan, M.~J. Wainwright, and Y.~Zhang.
\newblock Optimality guarantees for distributed statistical estimation.
\newblock \emph{arXiv preprint arXiv:1405.0782}, 2014.

\bibitem[Fan et~al.(2017)Fan, Wang, Wang, and Zhu]{fan2017distributed}
J.~Fan, D.~Wang, K.~Wang, and Z.~Zhu.
\newblock Distributed estimation of principal eigenspaces.
\newblock \emph{arXiv preprint arXiv:1702.06488}, 2017.

\bibitem[Guo et~al.(2017)Guo, Shi, and Wu]{guo2017learning}
Z.-C. Guo, L.~Shi, and Q.~Wu.
\newblock Learning theory of distributed regression with bias corrected
  regularization kernel network.
\newblock \emph{The Journal of Machine Learning Research}, 18\penalty0
  (1):\penalty0 4237--4261, 2017.

\bibitem[Hachem et~al.(2007)Hachem, Loubaton, and
  Najim]{hachem2007deterministic}
W.~Hachem, P.~Loubaton, and J.~Najim.
\newblock Deterministic equivalents for certain functionals of large random
  matrices.
\newblock \emph{The Annals of Applied Probability}, 17\penalty0 (3):\penalty0
  875--930, 2007.

\bibitem[Hiai and Petz(2006)]{hiai2006semicircle}
F.~Hiai and D.~Petz.
\newblock \emph{The semicircle law, free random variables and entropy}.
\newblock Number~77. American Mathematical Soc., 2006.

\bibitem[Huo and Cao(2018)]{huo2018aggregated}
X.~Huo and S.~Cao.
\newblock Aggregated inference.
\newblock \emph{Wiley Interdisciplinary Reviews: Computational Statistics},
  page e1451, 2018.

\bibitem[Jiang(1996)]{jiang1996reml}
J.~Jiang.
\newblock Reml estimation: asymptotic behavior and related topics.
\newblock \emph{The Annals of Statistics}, 24\penalty0 (1):\penalty0 255--286,
  1996.

\bibitem[Jiang et~al.(2016)Jiang, Li, Paul, Yang, and Zhao]{jiang2016high}
J.~Jiang, C.~Li, D.~Paul, C.~Yang, and H.~Zhao.
\newblock On high-dimensional misspecified mixed model analysis in genome-wide
  association study.
\newblock \emph{The Annals of Statistics}, 44\penalty0 (5):\penalty0
  2127--2160, 2016.

\bibitem[Jordan et~al.(2016)Jordan, Lee, and Yang]{jordan2016communication}
M.~I. Jordan, J.~D. Lee, and Y.~Yang.
\newblock Communication-efficient distributed statistical inference.
\newblock \emph{arXiv preprint arXiv:1605.07689}, 2016.

\bibitem[Koutris et~al.(2018)Koutris, Salihoglu, Suciu,
  et~al.]{koutris2018algorithmic}
P.~Koutris, S.~Salihoglu, D.~Suciu, et~al.
\newblock Algorithmic aspects of parallel data processing.
\newblock \emph{Foundations and Trends{\textregistered} in Databases},
  8\penalty0 (4):\penalty0 239--370, 2018.

\bibitem[Lee et~al.(2017)Lee, Liu, Sun, and Taylor]{lee2017communication}
J.~D. Lee, Q.~Liu, Y.~Sun, and J.~E. Taylor.
\newblock Communication-efficient sparse regression.
\newblock \emph{Journal of Machine Learning Research}, 18\penalty0
  (5):\penalty0 1--30, 2017.

\bibitem[Li et~al.(2013)Li, Lin, and Li]{li2013statistical}
R.~Li, D.~K. Lin, and B.~Li.
\newblock Statistical inference in massive data sets.
\newblock \emph{Applied Stochastic Models in Business and Industry},
  29\penalty0 (5):\penalty0 399--409, 2013.

\bibitem[Lin et~al.(2017)Lin, Guo, and Zhou]{lin2017distributed}
S.-B. Lin, X.~Guo, and D.-X. Zhou.
\newblock Distributed learning with regularized least squares.
\newblock \emph{The Journal of Machine Learning Research}, 18\penalty0
  (1):\penalty0 3202--3232, 2017.

\bibitem[Liu et~al.(2018)Liu, Shang, and Cheng]{liu2018many}
M.~Liu, Z.~Shang, and G.~Cheng.
\newblock How many machines can we use in parallel computing for kernel ridge
  regression?
\newblock \emph{arXiv preprint arXiv:1805.09948}, 2018.

\bibitem[Lynch(1996)]{lynch1996distributed}
N.~A. Lynch.
\newblock \emph{Distributed algorithms}.
\newblock Elsevier, 1996.

\bibitem[Mackey et~al.(2011)Mackey, Jordan, and Talwalkar]{mackey2011divide}
L.~W. Mackey, M.~I. Jordan, and A.~Talwalkar.
\newblock Divide-and-conquer matrix factorization.
\newblock In \emph{Advances in neural information processing systems}, pages
  1134--1142, 2011.

\bibitem[Marchenko and Pastur(1967)]{marchenko1967distribution}
V.~A. Marchenko and L.~A. Pastur.
\newblock Distribution of eigenvalues for some sets of random matrices.
\newblock \emph{Mat. Sb.}, 114\penalty0 (4):\penalty0 507--536, 1967.

\bibitem[Mcdonald et~al.(2009)Mcdonald, Mohri, Silberman, Walker, and
  Mann]{mcdonald2009efficient}
R.~Mcdonald, M.~Mohri, N.~Silberman, D.~Walker, and G.~S. Mann.
\newblock Efficient large-scale distributed training of conditional maximum
  entropy models.
\newblock In \emph{Advances in Neural Information Processing Systems}, pages
  1231--1239, 2009.

\bibitem[Nica and Speicher(2006)]{nica2006lectures}
A.~Nica and R.~Speicher.
\newblock \emph{Lectures on the combinatorics of free probability}, volume~13.
\newblock Cambridge University Press, 2006.

\bibitem[Paul and Aue(2014)]{paul2014random}
D.~Paul and A.~Aue.
\newblock Random matrix theory in statistics: A review.
\newblock \emph{Journal of Statistical Planning and Inference}, 150:\penalty0
  1--29, 2014.

\bibitem[Pourshafeie et~al.(2018)Pourshafeie, Bustamante, and
  Prabhu]{pourshafeie2018caring}
A.~Pourshafeie, C.~D. Bustamante, and S.~Prabhu.
\newblock Caring without sharing: Meta-analysis 2.0 for massive genome-wide
  association studies.
\newblock \emph{bioRxiv}, page 436766, 2018.

\bibitem[Rauber and R{\"u}nger(2013)]{rauber2013parallel}
T.~Rauber and G.~R{\"u}nger.
\newblock \emph{Parallel programming: For multicore and cluster systems}.
\newblock Springer Science \& Business Media, 2013.

\bibitem[Rosenblatt and Nadler(2016)]{rosenblatt2016optimality}
J.~D. Rosenblatt and B.~Nadler.
\newblock On the optimality of averaging in distributed statistical learning.
\newblock \emph{Information and Inference: A Journal of the IMA}, 5\penalty0
  (4):\penalty0 379--404, 2016.

\bibitem[Rubio and Mestre(2011)]{rubio2011spectral}
F.~Rubio and X.~Mestre.
\newblock Spectral convergence for a general class of random matrices.
\newblock \emph{Statistics \& Probability Letters}, 81\penalty0 (5):\penalty0
  592--602, 2011.

\bibitem[Searle et~al.(2009)Searle, Casella, and McCulloch]{searle2009variance}
S.~R. Searle, G.~Casella, and C.~E. McCulloch.
\newblock \emph{Variance components}, volume 391.
\newblock John Wiley \& Sons, 2009.

\bibitem[Serdobolskii(2007)]{serdobolskii2007multiparametric}
V.~I. Serdobolskii.
\newblock \emph{Multiparametric {S}tatistics}.
\newblock Elsevier, 2007.

\bibitem[Shang and Cheng(2017)]{shang2017computational}
Z.~Shang and G.~Cheng.
\newblock Computational limits of a distributed algorithm for smoothing spline.
\newblock \emph{The Journal of Machine Learning Research}, 18\penalty0
  (1):\penalty0 3809--3845, 2017.

\bibitem[Shi et~al.(2018)Shi, Lu, and Song]{shi2018massive}
C.~Shi, W.~Lu, and R.~Song.
\newblock A massive data framework for m-estimators with cubic-rate.
\newblock \emph{Journal of the American Statistical Association}, pages 1--12,
  2018.

\bibitem[Smith et~al.(2016)Smith, Forte, Ma, Tak{\'a}c, Jordan, and
  Jaggi]{smith2016cocoa}
V.~Smith, S.~Forte, C.~Ma, M.~Tak{\'a}c, M.~I. Jordan, and M.~Jaggi.
\newblock Cocoa: A general framework for communication-efficient distributed
  optimization.
\newblock \emph{arXiv preprint arXiv:1611.02189}, 2016.

\bibitem[Smola(2012)]{SML}
A.~Smola.
\newblock Course notes on scalable machine learning, 2012.

\bibitem[Voiculescu et~al.(1992)Voiculescu, Dykema, and
  Nica]{voiculescu1992free}
D.~V. Voiculescu, K.~J. Dykema, and A.~Nica.
\newblock \emph{Free random variables}.
\newblock Number~1. American Mathematical Soc., 1992.

\bibitem[Volgushev et~al.(2017)Volgushev, Chao, and
  Cheng]{volgushev2017distributed}
S.~Volgushev, S.-K. Chao, and G.~Cheng.
\newblock Distributed inference for quantile regression processes.
\newblock \emph{arXiv preprint arXiv:1701.06088}, 2017.

\bibitem[Wang et~al.(2018)Wang, Yang, Chen, and Liu]{xichen2018svm}
X.~Wang, Z.~Yang, X.~Chen, and W.~Liu.
\newblock Distributed inference for linear support vector machine.
\newblock \emph{arXiv preprint arxiv:1811.11922}, 2018.

\bibitem[Xu et~al.(2016)Xu, Shang, and Cheng]{xu2016optimal}
G.~Xu, Z.~Shang, and G.~Cheng.
\newblock Optimal tuning for divide-and-conquer kernel ridge regression with
  massive data.
\newblock \emph{arXiv preprint arXiv:1612.05907}, 2016.

\bibitem[Yao et~al.(2015)Yao, Bai, and Zheng]{yao2015large}
J.~Yao, Z.~Bai, and S.~Zheng.
\newblock \emph{Large Sample Covariance Matrices and High-Dimensional Data
  Analysis}.
\newblock Cambridge University Press, 2015.

\bibitem[Zhang et~al.(2012)Zhang, Wainwright, and
  Duchi]{zhang2012communication}
Y.~Zhang, M.~J. Wainwright, and J.~C. Duchi.
\newblock Communication-efficient algorithms for statistical optimization.
\newblock In \emph{Advances in Neural Information Processing Systems}, pages
  1502--1510, 2012.

\bibitem[Zhang et~al.(2013{\natexlab{a}})Zhang, Duchi, and
  Wainwright]{zhang2013divide}
Y.~Zhang, J.~Duchi, and M.~Wainwright.
\newblock Divide and conquer kernel ridge regression.
\newblock In \emph{Conference on Learning Theory}, pages 592--617,
  2013{\natexlab{a}}.

\bibitem[Zhang et~al.(2013{\natexlab{b}})Zhang, Duchi, and
  Wainwright]{zhang2013communication}
Y.~Zhang, J.~C. Duchi, and M.~J. Wainwright.
\newblock Communication-efficient algorithms for statistical optimization.
\newblock \emph{Journal of Machine Learning Research}, 14:\penalty0 3321--3363,
  2013{\natexlab{b}}.

\bibitem[Zhang et~al.(2015)Zhang, Duchi, and Wainwright]{zhang2015divide}
Y.~Zhang, J.~Duchi, and M.~Wainwright.
\newblock Divide and conquer kernel ridge regression: A distributed algorithm
  with minimax optimal rates.
\newblock \emph{The Journal of Machine Learning Research}, 16\penalty0
  (1):\penalty0 3299--3340, 2015.

\bibitem[Zhao et~al.(2016)Zhao, Cheng, and Liu]{zhao2016partially}
T.~Zhao, G.~Cheng, and H.~Liu.
\newblock A partially linear framework for massive heterogeneous data.
\newblock \emph{Annals of statistics}, 44\penalty0 (4):\penalty0 1400, 2016.

\bibitem[Zhu and Lafferty(2018)]{zhu2018distributed}
Y.~Zhu and J.~Lafferty.
\newblock Distributed nonparametric regression under communication constraints.
\newblock \emph{arXiv preprint arXiv:1803.01302}, 2018.

\end{thebibliography}
}

\end{document}